\pdfoutput=1
\RequirePackage{ifpdf}
\ifpdf 
\documentclass[pdftex]{sigma}
\else
\documentclass{sigma}
\fi

\numberwithin{equation}{section}
\usepackage{mathtools}

\newcommand{\nco}{\newcommand}
\def\R{\mathbb{R}}
\nco{\red}{\color{red}}
\nco{\blue}{\color{blue}}
\def\inv#1{\frac{1}{#1}}
\def\tr{\operatorname{tr}}
\def\({\left(}
\def\){\right)}

\nco{\rnc}{\renewcommand}

\def\p{{\tt p}}

\def\bbeta{\underline{\beta}}
\def\diag{\operatorname{diag}}
\def\ii{\mathrm{i}}
\def\oh{\frac{1}{2}}
\def\su{{\rm su}}
\def\SU{{\rm SU}}\def\U{{\rm U}}\def\SO{{\rm SO}}\def\O{{\rm O}}\def\USp{{\rm USp}}
\def\inv#1{\frac{1}{#1}}

\def\eq=#1{\buildrel #1 \over{=}}
\def\CH{{\mathcal H}} \def\CI{{\mathcal I}} \def\CJ{{\mathcal J}}  

\def\bH{\widetilde{\mathbf{H}}}
\def\GO{\Omega}
\def\E{\mathbb{E}}
\def\I{\mathbb{I}}

\def\R{\mathbb{R}}

\begin{document}

\allowdisplaybreaks

\newcommand{\arXivNumber}{1809.03394}

\renewcommand{\PaperNumber}{029}

\FirstPageHeading

\ShortArticleName{The Horn Problem for Real Symmetric and Quaternionic Self-Dual Matrices}

\ArticleName{The Horn Problem for Real Symmetric\\ and Quaternionic Self-Dual Matrices}

\Author{Robert COQUEREAUX~$^\dag$ and Jean-Bernard ZUBER~$^\ddag$}

\AuthorNameForHeading{R.~Coquereaux and J.-B.~Zuber}

\Address{$^\dag$~Aix Marseille Univ, Universit\'e de Toulon, CNRS, CPT, Marseille, France}
\EmailD{\href{mailto:robert.coquereaux@cpt.univ-mrs.fr}{robert.coquereaux@cpt.univ-mrs.fr}}
\URLaddressD{\url{http://www.cpt.univ-mrs.fr/~coque/}}

\Address{$^\ddag$~Sorbonne Universit\'e, UMR 7589, LPTHE, F-75005, Paris, France\\
\hphantom{$^\ddag$}~CNRS, UMR 7589, LPTHE, F-75005, Paris, France}
\EmailD{\href{mailto:jean-bernard.zuber@upmc.fr}{jean-bernard.zuber@upmc.fr}}
\URLaddressD{\url{http://www.lpthe.jussieu.fr/~zuber/}}

\ArticleDates{Received December 20, 2018, in final form April 06, 2019; Published online April 16, 2019}

\Abstract{Horn's problem, i.e., the study of the eigenvalues of the sum $C=A+B$ of two matrices, given the spectrum of $A$ and
of $B$, is re-examined, comparing the case of real symmetric, complex Hermitian and self-dual quaternionic $3\times 3$ matrices.
In particular, what can be said on the probability distribution function (PDF) of the eigenvalues of $C$ if $A$ and $B$
are independently and uniformly distributed on their orbit under the action of, respectively, the orthogonal, unitary and
symplectic group?
While the two latter cases (Hermitian and quaternionic)
may be studied by use of explicit formulae for the relevant orbital integrals, the case of
real symmetric matrices is much harder. It is also quite intriguing, since numerical experiments reveal the
occurrence of singularities where the PDF of the eigenvalues diverges.
Here we show that the computation of the PDF
of the symmetric functions of the eigenvalues for traceless $3\times 3$ matrices may be carried out in terms of
algebraic functions~-- roots of quartic polynomials~-- and their integrals. The computation is carried out in detail in a
particular case, and reproduces the expected singular patterns. The divergences are of logarithmic or inverse power type.
We also relate this PDF to the (rescaled) structure constants of zonal polynomials and introduce a
zonal analogue of the Weyl $\SU(n)$ characters.}


\Classification{17B08; 17B10; 22E46; 43A75; 52Bxx}

\section{Introduction}
Recall what Horn's problem is: given two $n\times n$ matrices $A$ and $B$ of given spectrum of eigenvalues, what can be said about the spectrum of their sum $C=A+B$? The problem has been addressed by many authors, see~\cite{Fu} for a review and references. For Hermitian matrices, the support of the spectrum of $C$ has been completely determined after some crucial work by Klyashko~\cite{Kl} and by Knutson and Tao~\cite{KT, KT00, KTW01}.

In a recent paper \cite{Z1}, a more specific question was considered: given $\alpha=\{\alpha_1\ge \alpha_2\ge \cdots \ge \alpha_n\}$ and $\beta=\{\beta_1\ge \beta_2\ge \cdots \ge \beta_n\}$, take the Hermitian matrices $A$ and $B$ uniformly and independently distributed on the orbit of $\diag(\alpha)$ and $\diag(\beta)$ under the action of the $\SU(n)$ group. The probability distribution function (PDF) of the eigenvalues $\gamma$ of $C=A+B$ may then be computed, see also \cite{Fa, FG}.

The same question may, however, be raised if instead of Hermitian matrices, one considers other classes of matrices and the action of an appropriate group. The case of skew-symmetric real matrices was considered in \cite{Z1}, but more interesting is that of real symmetric matrices and the action of the orthogonal group $\SO(n)$.

Let us start with a numerical experiment with specific $3\times3$ matrices. Since $\tr C=\tr A +\tr B$, only two eigenvalues of $C$, say $\gamma_1$ and $\gamma_2$ are linearly independent. Then compare the distribution of points in the $(\gamma_1, \gamma_2)$-plane for the three cases of
\begin{enumerate}\itemsep=0pt
\item[(a)] orbits of real symmetric matrices $A$ and $B$ equivalent to $J_z:=\diag(1,0,-1)$ under conjugation by the orthogonal group $\SO(3)$;
\item[(b)] orbits of such matrices, regarded as Hermitian, under conjugation by the unitary group $\SU(3)$;
\item[(c)] orbits of such matrices, regarded as quaternionic self-dual (i.e., $A\otimes \I_2$, $B\otimes \I_2$ as $6\times6$ matrices), under conjugation by the unitary symplectic group $\USp(3)$.
\end{enumerate}

Following the nomenclature introduced by Dyson, we label these three cases by the index $\bbeta= 1,2$ or~4, respectively.\footnote{The reader should not confuse this index $\bbeta$ with the multiplet $\beta$ of eigenvalues of the $B$ matrix.}

\begin{figure}[t!] \centering
 \includegraphics[width=12pc]{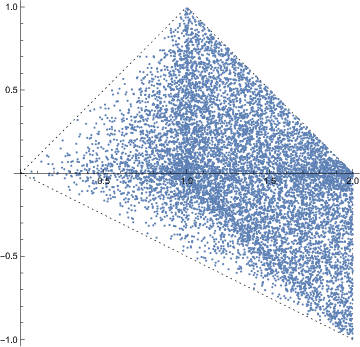}\ \includegraphics[width=12pc]{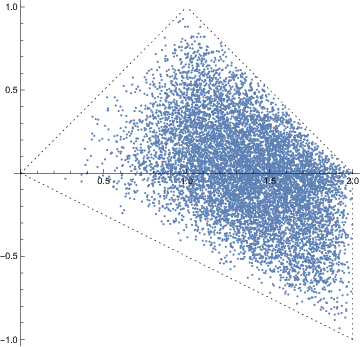} \ \includegraphics[width=12pc]{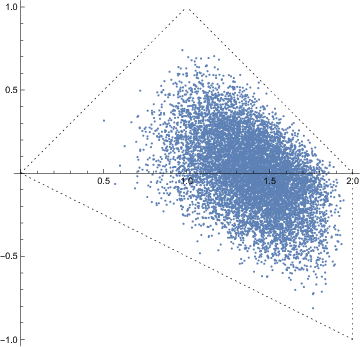} \\
 \includegraphics[width=12pc]{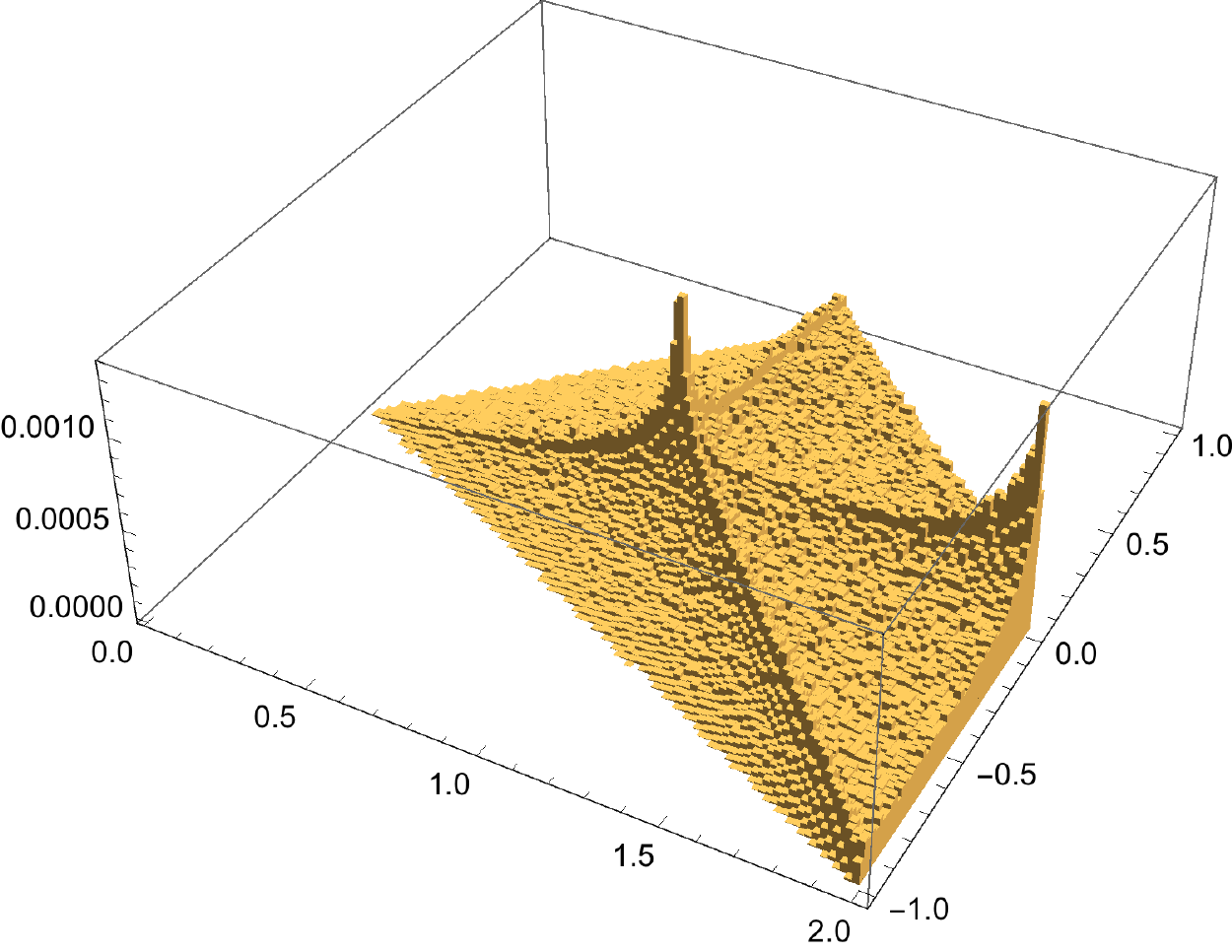} \ \includegraphics[width=12pc]{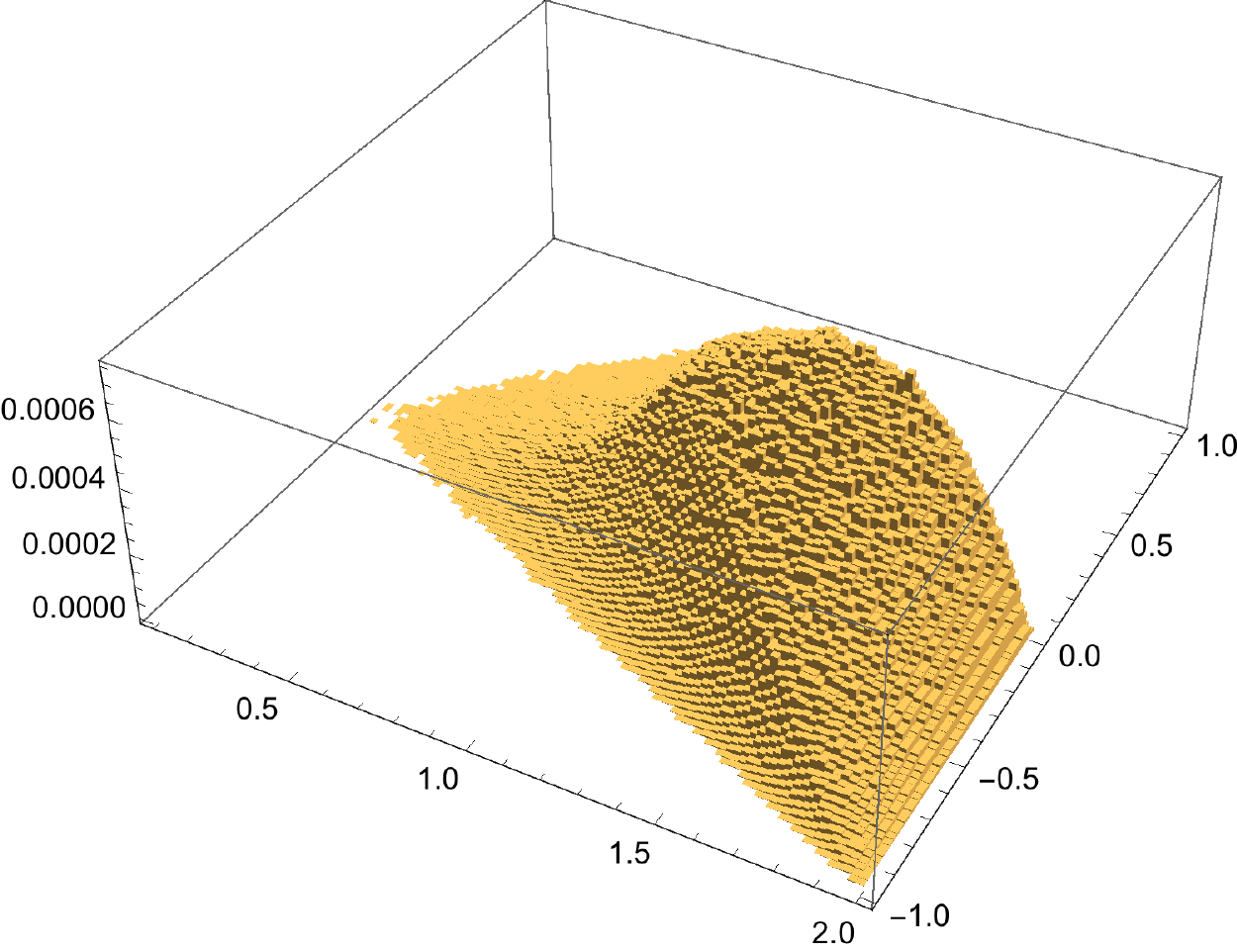}\ \includegraphics[width=12pc]{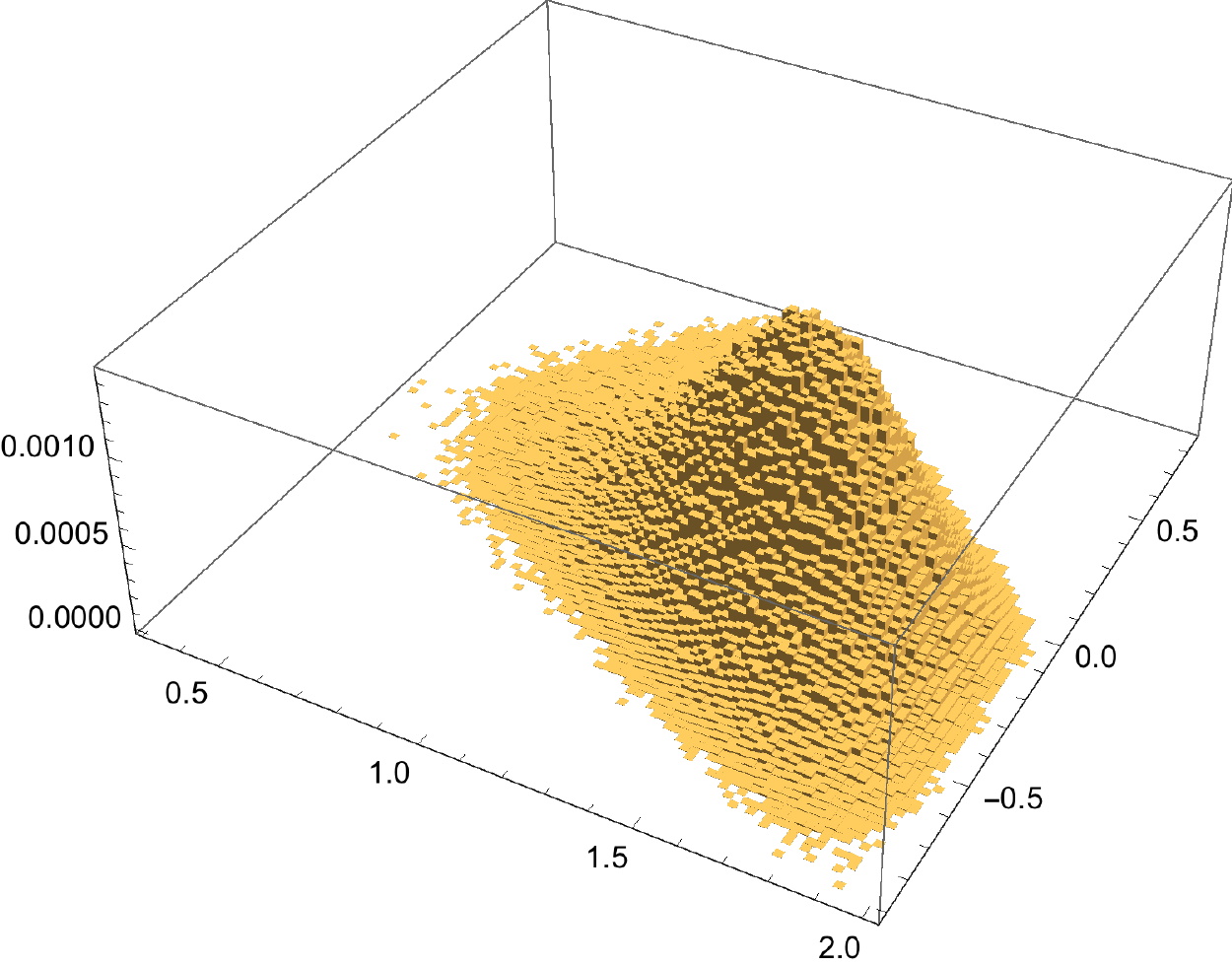}\\[5pt]
 \qquad \qquad (a) \qquad \qquad \qquad \qquad \qquad \qquad(b) \qquad \qquad \qquad \qquad \qquad \qquad \qquad(c)
 \caption{Plots and histograms of $(\gamma_1,\gamma_2)$ for $A\sim B\sim J_z$ in the three cases: (a) orbits under SO(3), (b) under SU(3), (c) USp(3). Plots made of $10^4$ points, histograms of $10^6$ points.}\label{threecases}
\end{figure}

Some features appear clearly on the plots and histograms of Fig.~\ref{threecases}:
\begin{enumerate}\itemsep=0pt
\item[(i)] the PDF vanishes faster on the boundaries of the Horn domain as $\bbeta$ increases;
\item[(ii)] the non-analyticities are stronger and stronger as $\bbeta$ decreases;
\item[(iii)] these singularities seem to appear at the same place in the $(\gamma_1,\gamma_2)$ plane (for $\alpha$ and $\beta$ fixed).
\end{enumerate}
In the Hermitian case (and in the quaternionic case as well, see below), it is known that the PDF is a~piece-wise polynomial function. The plots of Fig.~\ref{threecases}(a) suggest that this cannot be true for real symmetric matrices. It should be emphasized that these general features do not depend on the explicit case we have chosen. Similar singular patterns have been observed in numerical experiments with other matrices $A$ and $B$, see \cite[Fig.~7]{Z1}.

\looseness=-1 The aim of this paper is to compare the three cases, to reproduce analytically the previous empirical observations and in particular to analyse the location and nature of the singularities that occur in the symmetric case. After a brief review of known results on the relevant orbital integrals (Section~\ref{sectionorbitalintegrals}), we treat rapidly the easy case of quaternionic self-dual matrices (Section~\ref{quaternion}), before turning to the more challenging case of real symmetric matrices in Sections~\ref{PDF} and~\ref{sectionparticularcase}. In Section~\ref{PDF}, it is shown that for $3\times 3$ traceless matrices, the introduction of the two symmetric functions $p$ and $q$ of the three eigenvalues (of vanishing sum) simplifies matters: one may express the PDF $\rho(p,q)$ in terms of roots of some polynomial equations and integral thereof, see~(\ref{pdf4}) below. The actual computation is carried out again for our pet example of $\alpha=\beta=(1,0,-1)$ in Section~\ref{sectionparticularcase}. In particular we reproduce and analyze in Section~\ref{gen-sing} and Appendix~\ref{detail-sing} the singularities that are apparent in Fig.~\ref{threecases}(a). Finally in Section~\ref{zonalpolynomials}, we show that this function $\rho(p,q)$ is related to the distribution of the (rescaled) structure constants of zonal polynomials, thus elaborating on a claim of~\cite{FG2}.

\section[The orbital integrals for $\bbeta=1,2,4$]{The orbital integrals for $\boldsymbol{\bbeta=1,2,4}$}\label{sectionorbitalintegrals}
Let us consider the {\it orbital integrals}
\begin{gather*}
\CI_{\bbeta}(X,A)=\int_{G_{\bbeta}} {\rm d}\GO \exp\big[\tr X \GO A \GO^\dagger\big]\end{gather*}
in the following three cases:
\begin{enumerate}\itemsep=0pt
\item[1)] $\bbeta=1$, $\GO\in G_1:=\SO(n)$, $X$, $A$ real symmetric matrices;
\item[2)] $\bbeta=2$, $\GO\in G_2:=\U(n)$, $X$, $A$ complex Hermitian matrices;
\item[3)] $\bbeta=4$, $\GO\in G_4:=\USp(n)$, $X$, $A$ real quaternionic self-dual matrices.
\end{enumerate}
In fact each of these integrals depends only on the eigenvalues $\alpha_i$, resp.~$x_i$, of $A$, resp.~$X$, and by a small abuse of notations, we write them $\CI_{\bbeta}(x,\alpha)$ in the following.

For $\bbeta=2$, this orbital integral is well known~\cite{HC,IZ}
\begin{gather*}
\CI_2(x,\alpha)=\prod_{p=1}^{n-1} p! \sum_{P\in S_n} \inv{\Delta(x)\Delta(\alpha_P)} {\rm e}^{\sum x_j \alpha_{Pj}}=
\prod_{p=1}^{n-1} p! \frac{\det {\rm e}^{x_i \alpha_j}}{\Delta(x)\Delta(\alpha)}. \end{gather*}

Here and below, $\Delta$ with no subscript is the Vandermonde determinant $\Delta(x)=\prod\limits_{i<j}(x_i-x_j)$.

For $\bbeta=4$ and for a given $n$, as shown by \cite{BH1, BH2}, we have closed formulae
\begin{gather}\label{orbintquater}
\CI_4(x,\alpha)=\hat\kappa_n \sum_{P\in S_n} \inv{\Delta(x)^3\Delta(\alpha_P)^3} {\rm e}^{\sum x_j \alpha_{Pj}} f_n(\tau(\alpha,P)),\end{gather}
where $f_n$ a polynomial in the variables
\begin{gather} \label{tau}\tau(\alpha,P)=\{\tau_{ij}:=(x_i-x_j)(\alpha_{Pi}-\alpha_{Pj}) \,\vert\, 1\le i<j\le n\},\end{gather}
and $\hat\kappa_n$ is an $\alpha$- and $x$-independent constant. For example for $n=2$, $\hat\kappa_2=-12$ and $f_2(\tau(\alpha,P)) =1-\oh \tau_{12}$ while for $n=3$, $\hat\kappa_3=-2\times 3! \times 6!$ and{\samepage
\begin{gather*}
f_3(\tau(\alpha,P)) =1-\inv{3} (\tau_{12}+\tau_{13}+\tau_{23})+\inv{6}( \tau_{12}\tau_{13}+\tau_{13}\tau_{23}+\tau_{23}\tau_{12}) -\inv{12} \tau_{12}\tau_{13}\tau_{23}.
\end{gather*}
These closed formulae allow an explicit computation of the PDF, see \cite{Z1} and below, Section~\ref{quaternion}.}

In contrast, for $\bbeta=1$, i.e., for real symmetric matrices under the action of the orthogonal group~$\SO(n)$, the best that can be achieved is an expansion in zonal polynomials \cite{BH1, FG2, BH2, JamesZonalPoly2, OO}
 \begin{gather} \label{CI1}\CI_1(x,\alpha)= \sum_{m=0}^\infty \inv{m! \prod\limits_{q=0}^{m-1}(1+2q)} \sum_{\kappa\vdash m}c(\kappa) \frac{Z(\kappa)(x) Z(\kappa)(\alpha)}{Z(\kappa)(I)},\end{gather}
 where the second sum runs over partitions $\kappa$ of $m$ with no more than~$n$ parts, and~$c(\kappa)$ is a constant which depends on the normalization of the $Z(\kappa)$'s, see Section~\ref{zonalpolynomials}.
 If the zonal polynomials in the above formula are written with the so-called James normalization, one can use the values of~$Z(\kappa)(I)$ and of the coefficient $c(\kappa)$ tabulated by James up to $m=4$, for all $n$ (the dimension of $I$) in the Appendix of \cite[p.~157]{JamesZonalPoly2}; one can also use commands contained in the Mathematica package~\cite{MathematicaProgramsRC}. The infinite sum in the previous expansion, however, makes the computation of the PDF intractable (as far as we can see), and we will have to follow another route, see below Section~\ref{PDF}.

In Section~\ref{zonalpolynomials}, however, we return to this formulation in terms of zonal polynomials and show a connection between the PDF and the distribution of (rescaled) ``zonal multiplicities'' i.e., appropriate structure constants of zonal polynomials.

\section[Quaternionic case for $n=3$]{Quaternionic case for $\boldsymbol{n=3}$}\label{quaternion}
We start with a warming up exercise: compute the PDF of the eigenvalues $\gamma$ for a sum of two self-dual quaternionic matrices that are independently and uniformly distributed on their orbit under the action of the unitary symplectic group. The computation of the PDF of the $\gamma$'s follows the same lines as that in the Hermitian case. We will therefore be a bit sketchy in its derivation, referring the reader to \cite{Z1} for details of the computation.

Up to an overall factor, this PDF is given by the integral of three orbital integrals of the type~$\CI_4$ in~(\ref{orbintquater}),
\begin{gather*}
\p(\gamma|\alpha,\beta)=\frac{\kappa_n^2}{(2\pi)^N} \Delta(\gamma)^4 \int \prod_{i=1}^{n} {\rm d}x_i \Delta(x)^4 \CI_4(x,\ii \alpha) \CI_4(x,\ii \beta) \CI_4(x,-\ii \gamma), \end{gather*}
where $N=n(2n-1)$ is the number of independent matrix elements of a self-dual quaternion matrix and $\kappa_n=\frac{(2\pi)^{(N-n)/2}}{\prod\limits_{j=1}^n ( (2j)!/2)}$ stems from the change from those~$N$ variables to the $n$ eigenvalues. Thus
\begin{gather} \label{pUSpn} \p(\gamma|\alpha,\beta) = \frac{2\kappa_n^2 \hat\kappa_n^3}{(2\pi)^N} { (-1)^{n(n-1)/2}} \pi^n \delta\bigg(\sum_k(\alpha_k+\beta_k-\gamma_k)\bigg)
\frac{ \Delta(\gamma)}{\Delta(\alpha)^3\Delta(\beta)^3} \CJ_n,\\ \nonumber
 \CJ_n =\frac{ n! }{{\ii^{n(n-1)/2}}\pi^{n-1}} \sum_{P,P'\in S_n} \varepsilon_P \varepsilon_{P'}
\int \frac{{\rm d}^{n-1}u}{\Delta'^5(u)} \\
\hphantom{\CJ_n =}{} \times \prod_{j=1}^{n-1} {\rm e}^{\ii u_j A_j(P,P',I)} f_n(\tau(\ii \alpha,P)) f_n(\tau(\ii \beta,P')) f_n(\tau(-\ii \gamma,I)) \label{CJn}
\end{gather}
with $u_j:=x_j-x_{j+1}$, $\Delta'(u):=\prod\limits_{1\le i < j\le n} (u_i+u_{i+1}+\cdots +u_{j-1})$ and
\begin{gather} \label{Aj} A_j(P,P',P'')= \sum_{k=1}^j (\alpha_{P(k)}+\beta_{P'(k)}-\gamma_{P''(k)}).\end{gather}
The power~5 of $\Delta'(u)$ results from $3\times 3$ (three denominators $\Delta^3(x)$), minus~$4$ from the Jaco\-bian~$\Delta^4(x)$.

\begin{figure}[t!] \centering\includegraphics[width=14pc]{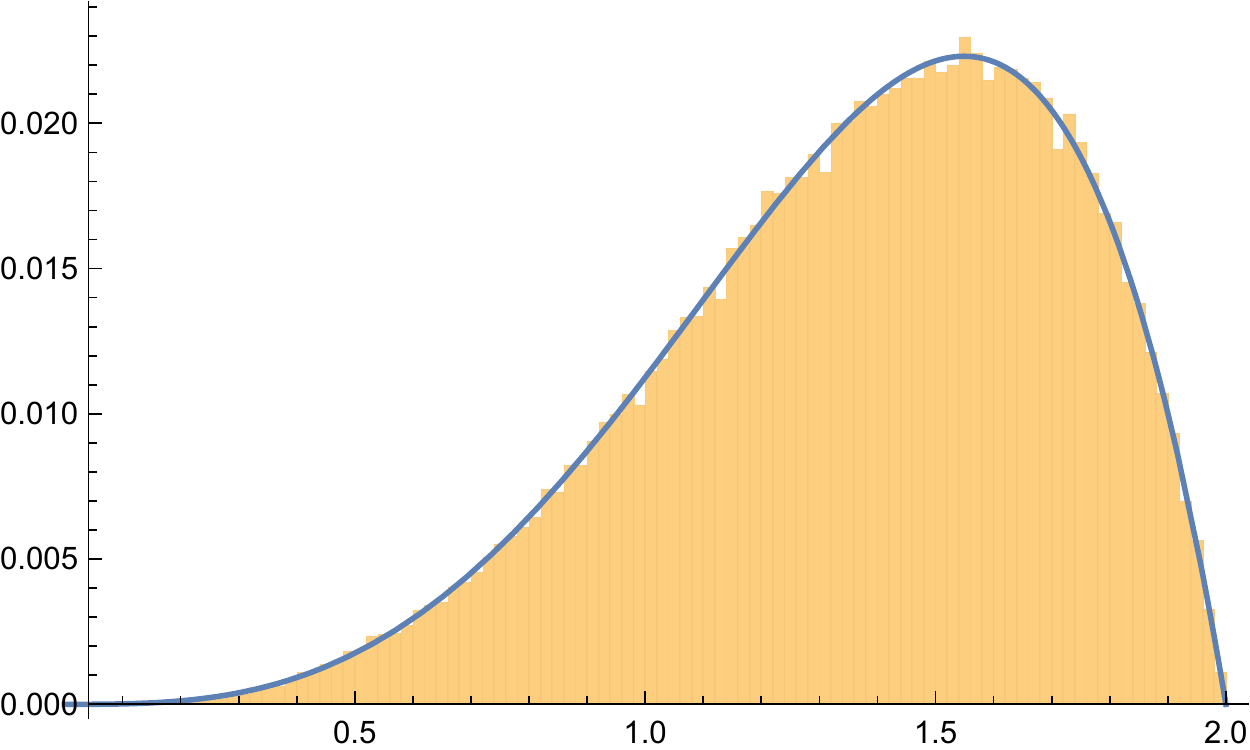}\qquad \includegraphics[width=14pc]{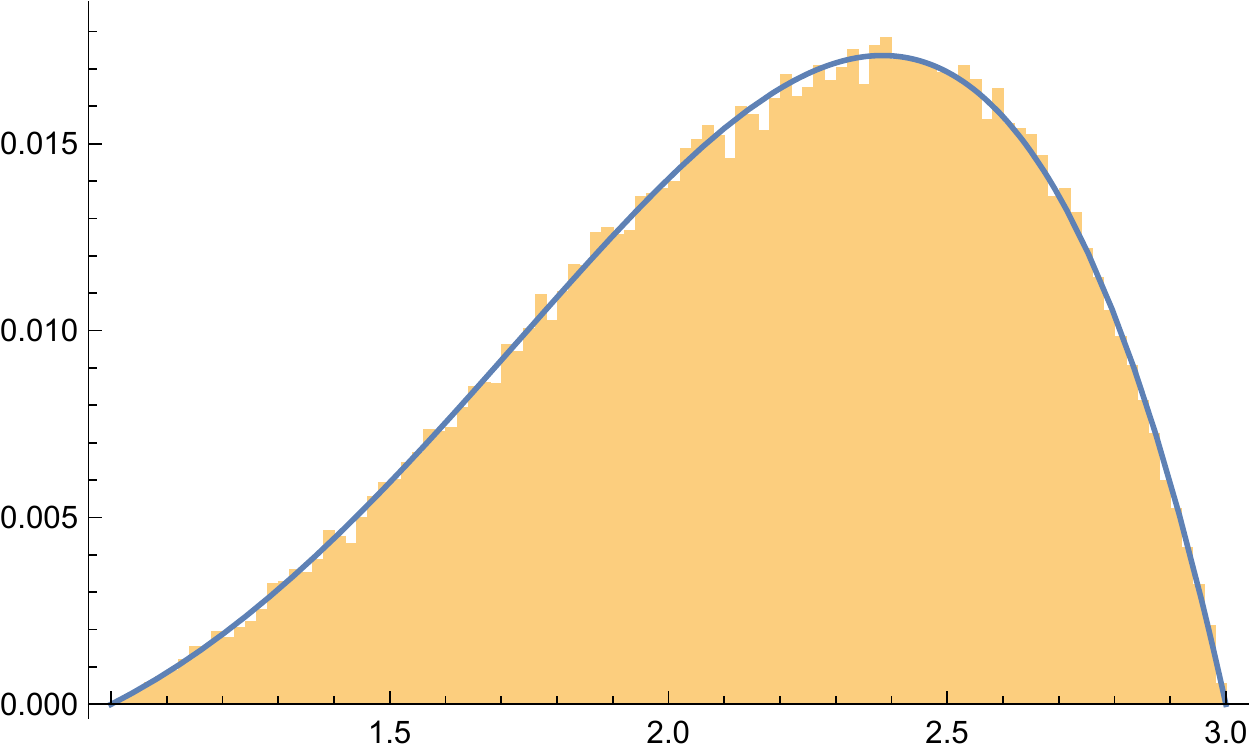}
 \caption{Comparing the histogram of $10^4$ points with the PDF of~(\ref{pdfqu2}) for (left) $\alpha=\beta=(1,-1)$ and (right) $\alpha=(1,-1)$, $\beta=(2,-2)$. }\label{histo-pdf-n=2}
\end{figure}

With the little extra complication caused by the $f_n$ factors of (\ref{orbintquater}), the integration may be carried out as in~\cite{Z1}, namely by partial fraction decomposition of the integrand and repeated use of the Dirichlet formula for the Cauchy principal value
\begin{gather*}
P\int_\Bbb{R} \frac{{\rm d}u}{u^r} {\rm e}^{\ii u A} =\ii \pi \frac{(\ii A)^{r-1}}{(r-1)!} \epsilon(A),
\end{gather*}
with $\epsilon$ the sign function. The result, though cumbersome, is completely explicit. For $n=2$, we leave as an exercise for the reader to check that for $\alpha_2=-\alpha_1$, $\beta_2=-\beta_1$, we have
\begin{gather*} \p(\gamma|\alpha,\beta)= 6 \frac{\Delta(\gamma)}{\Delta^3(\alpha) \Delta^3(\beta)} \delta(\gamma_1+\gamma_2) \CJ_2\end{gather*} with
\begin{gather}
\CJ_2=\frac{1}{4} \big({-}\big({\alpha_1^2-\beta_1^2}\big)^2+2 \big({\alpha_1}^2 +\beta_1^2\big){\gamma_1}^2-{\gamma_1}^4\big)\nonumber\\
\hphantom{\CJ_2=}{}\times
 \big({\epsilon}(\gamma_1-{\alpha_1}+{\beta_1}) +{\epsilon}({\gamma_1}+{\alpha_1}-{\beta_1})
 -{\epsilon}(\gamma_1-{\alpha_1}-{\beta_1}) -{\epsilon}({\gamma_1}+{\alpha_1}+{\beta_1})\big)\nonumber\\
\hphantom{\CJ_2}{}= \frac{1}{2} \big({-}\big({\alpha_1^2-\beta_1^2}\big)^2+2 \big({\alpha_1}^2 +\beta_1^2\big){\gamma_1}^2-{\gamma_1}^4\big)
 \big( {\bf 1}_I(\gamma_1)-{\bf 1}_{-I}(\gamma_1)\big), \label{pdfqu2} \end{gather}
 showing that the spectrum of $\gamma_1$ is supported by two segments $I:=[ |\alpha_1-\beta_1|, \alpha_1+\beta_1]$ and $-I$, or only by the former if one imposes $\gamma_2\le \gamma_1$. Compare with the analogous formulae obtained for Hermitian and real symmetric matrices in~\cite{Z1}. In Fig.~\ref{histo-pdf-n=2}, we show the agreement with the histogram made of $10^4$ for two pairs $\alpha=\beta=(1,-1)$ and $\alpha=(1,-1)$, $\beta=(2,-2)$.

For $n=3$,
\begin{gather} \label{p3} \p(\gamma|\alpha,\beta)= 8640 \frac{\Delta(\gamma)}{\Delta^3(\alpha)\Delta^3(\beta)} \delta\bigg(\sum_k \gamma_k-\alpha_k-\beta_k\bigg) \CJ_3,\end{gather}
where $\CJ_3$ is a piecewise polynomial of degree 13 in $\gamma_1$ and $\gamma_2$. Indeed
\begin{gather}\label{CJ3} \CJ_3 =\frac{3!}{\ii^3 \pi^2}\sum_{P,P'\in S_3} \varepsilon_P \varepsilon_{P'} \int \frac{{\rm d}u_1 {\rm d}u_2}{(u_1 u_2 (u_1+u_2))^5} {\rm e}^{\ii (u_1 A_1+u_2 A_2)} F(u_1, u_2),\end{gather}
where $F(u_1, u_2)$ stands for the product {$ f_3\big(\tau(\ii \alpha,P)\big) f_3\big(\tau(\ii \beta,P')\big) f_3\big(\tau(-\ii \gamma,I)\big)$ in~(\ref{CJn}), see also (\ref{tau})};
one therefore finds, by homogeneity, that $\CJ_3$ behaves as $[\gamma]^{15-2=13}$. The function $\CJ_3$ is of differentiability class $C^2$: indeed as $u_1$, say, goes to infinity, the integrand behaves as $ \frac{1}{u_1^4} {\rm e}^{\ii A_1 u_1}$ (up to subdominant terms), so that the integral is twice continuously differentiable with respect to~$A_1$ (or~$\gamma$), hence of class $C^2$. Non-analyticities are expected (and do occur) along the boundaries of the Horn polygon and across the same singular lines (or half-lines) as in the Hermitian case, namely
\begin{gather}
\gamma_2=\alpha_2+\beta_2;\nonumber \\
 \gamma_1= \alpha_1+\beta_2,\qquad \gamma_2\ge \alpha_3+\beta_1;\qquad \gamma_2= \alpha_3+\beta_1,\qquad \gamma_1\le \alpha_1+\beta_2;\nonumber\\
\gamma_3=\alpha_2+\beta_3,\qquad \gamma_1\ge \alpha_1+\beta_2, \qquad \text{the same with} \ \alpha\leftrightarrow \beta. \label{loci-sing}
\end{gather}

 \begin{figure}[t!] \centering\includegraphics[width=14pc]{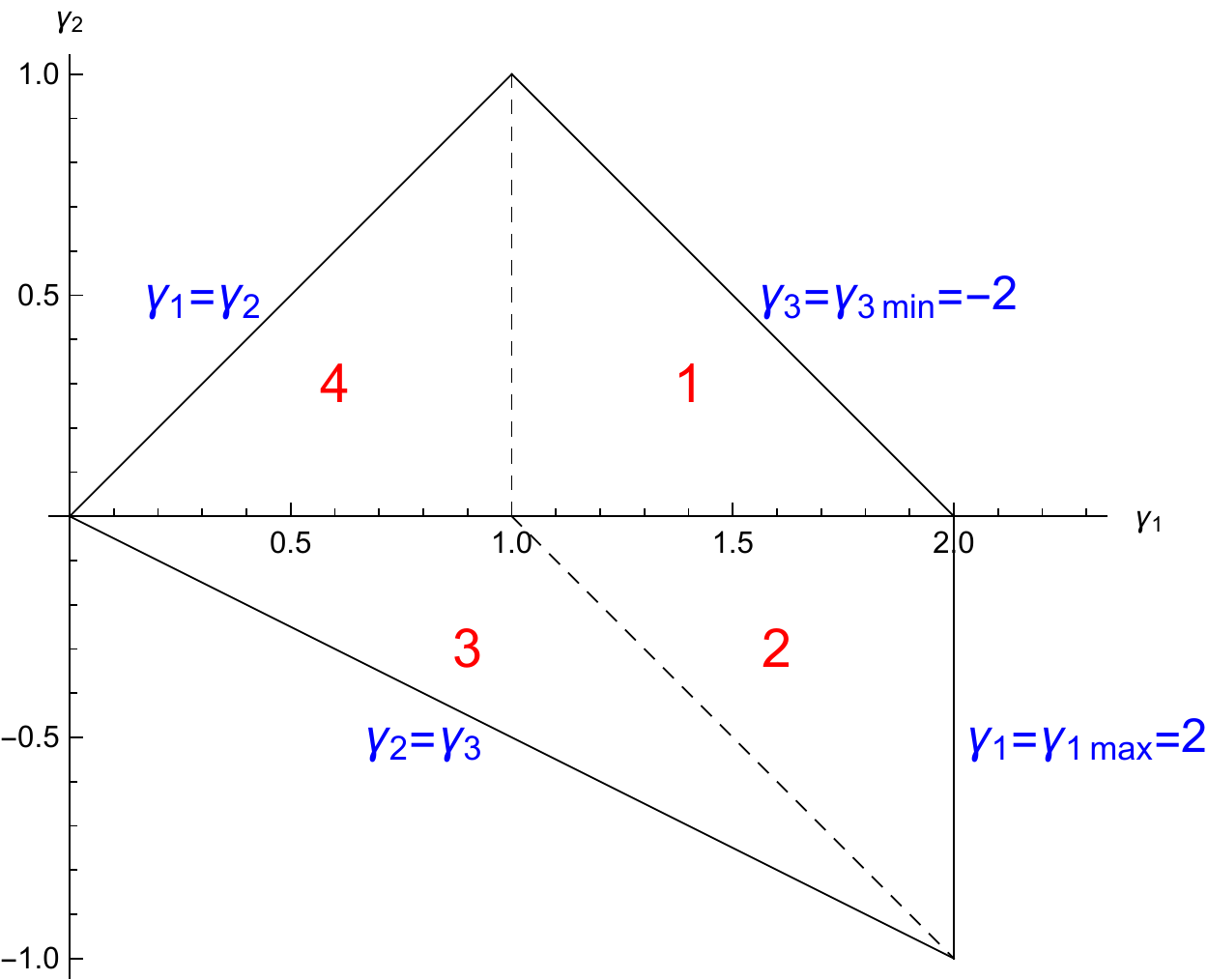}
 \caption{The four sectors in the Horn polygon of the quaternionic case for $\alpha=\beta=\{1,0,-1\}$.}\label{HornPolygonQuaternionic}
\end{figure}

For our favorite example of $A\sim B\sim J_z=\diag(1,0,-1)$, the function has four sectors of piecewise polynomiality, labelled by $i=1,\ldots, 4$ according to Fig.~\ref{HornPolygonQuaternionic}. In each of these four sectors the function takes the form
\begin{gather}\label{J3} \CJ_3|_{\mathrm{sector }\ i} =-\frac{2\times 5!}{13!} {\mathcal P}_i(\gamma_1, \gamma_2),\end{gather}
where
\begin{gather} \label{form-qu} {\mathcal P}_i(\gamma_1, \gamma_2) = p_i(\gamma_1, \gamma_2)^3 \tilde p_i(\gamma_1, \gamma_2), \\
 p_i(\gamma_1, \gamma_2)=\begin{cases}(2-\gamma_1- \gamma_2)&\text{in sector} \ 1,\ \text{i.e.}, \gamma_1\ge {1},\ \gamma_2\ge 0,\ 2-\gamma_1- \gamma_2\ge 0,\\
(2-\gamma_1) & \text{in sector} \ 2,\ \text{i.e.}, \gamma_1\le 2,\ \gamma_2\le 0, \ \gamma_1+\gamma_2\ge 1, \\
(\gamma_2-\gamma_3)=(\gamma_1+2\gamma_2)& \text{in sector} \ 3,\ \text{i.e.}, \gamma_2\le 0, \ \gamma_1+\gamma_2\le 1,\ \gamma_1+2\gamma_2\ge 0, \\
(\gamma_1-\gamma_2)& \text{in sector} \ 4,\ \text{i.e.}, \gamma_1\le 1,\ \gamma_2\ge 0,\ \gamma_1-\gamma_2\ge 0,
\end{cases} \nonumber
\end{gather}
and each of the $\tilde p_i(\gamma_1, \gamma_2)=40 \gamma_1^{10}+\cdots$ is too cumbersome to be given here.\footnote{This may be found on the web-site \url{http://www.lpthe.jussieu.fr/~zuber/Z_Unpub2.html}.
Other data that are not given explicitly in this paper may be found either on this web-site or on \url{http://www.cpt.univ-mrs.fr/~coque/Varia/JackZonalSchurResults.html}\label{repfootnote}}

Note that the form (\ref{form-qu}) guarantees that $\CJ_3$ vanishes ``cubically'' on each boundary of the polygon. This explains one of the features observed on the plots of Fig.~\ref{threecases}(c). As a side remark, we note that $\tilde p_4(\gamma_1, \gamma_2)$ is symmetric under $\gamma_1\leftrightarrow \gamma_2$, and $\tilde p_3(\gamma_1, \gamma_2)$ is symmetric under $\gamma_2\leftrightarrow \gamma_3=-\gamma_1-\gamma_2$.

One may also compute the transition functions between adjacent domains:
\begin{alignat*}{3}
& {\mathcal P}_1- {\mathcal P}_2 = \gamma_2^3 (\cdots),\qquad && {\mathcal P}_2- {\mathcal P}_3= (1-\gamma_1-\gamma_2)^3 (\cdots),& \\
& {\mathcal P}_3- {\mathcal P}_4= \gamma_2^3 (\cdots),\qquad && {\mathcal P}_4- {\mathcal P}_1= (1 - \gamma_1)^3 (\cdots),&
\end{alignat*}
where in each case, the ellipsis stands for a polynomial of degree 10 that we refrain from dis\-playing. Just like on the boundary, the transition functions vanish cubically along the non-analyticity lines. This is in agreement with the differentiability argument above.

\begin{figure}[t!] \centering
 \includegraphics[width=16pc]{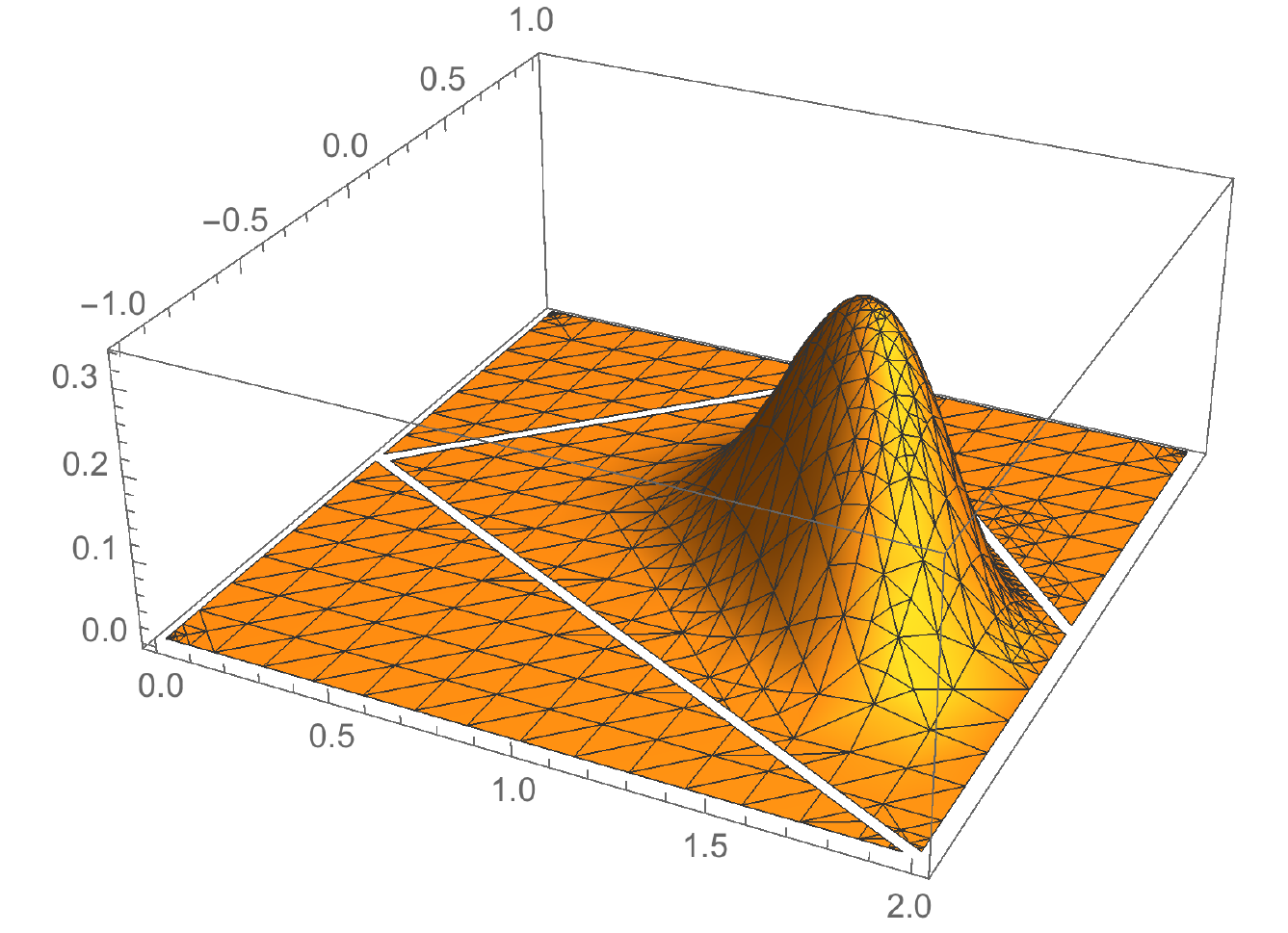}\qquad
 \includegraphics[width=18pc]{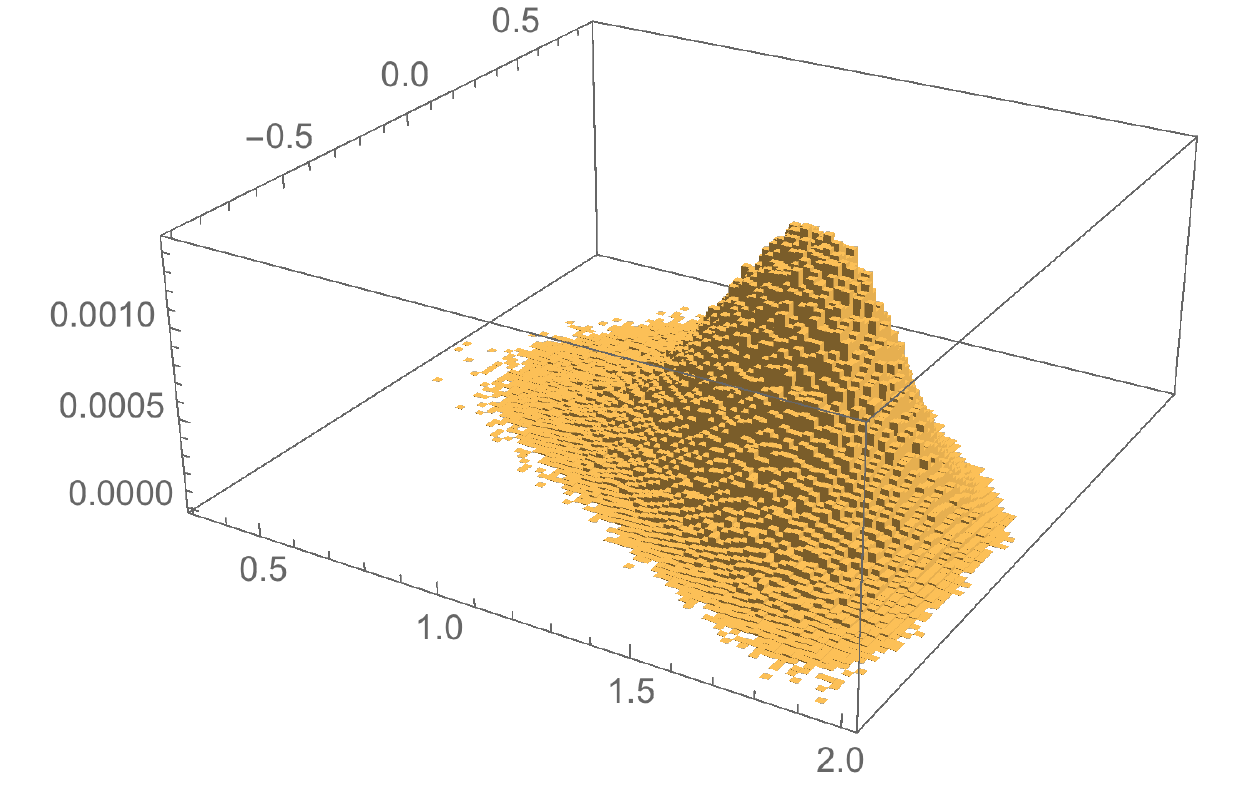}\\[5pt]
 \qquad (a) \qquad \qquad \qquad \qquad \qquad \qquad \qquad \qquad \qquad \qquad(b)
 \caption{Comparing the result of (\ref{pUSpn})--(\ref{form-qu}) with the ``experimental'' histogram.}\label{PDFQuaternionic}
\end{figure}

The resulting PDF is plotted in Fig.~\ref{PDFQuaternionic}(a) and compared to the ``experimental'' histogram of Fig.~\ref{PDFQuaternionic}(b).

A good check of our calculation is that the integral of the PDF over the whole domain of~$\gamma$'s, made of $3!$ copies of the Horn polygon $\bH_{\alpha,\beta}$ (see below in~(\ref{inequ3})), is indeed~1, thanks to
\begin{gather*} \int_{\bH_{\alpha,\beta}} {\rm d}\gamma_1 {\rm d}\gamma_2 \Delta(\gamma) \CJ_3= \inv{810}.\end{gather*}

\section[Computing the PDF for a sum for two real symmetric matrices]{Computing the PDF for a sum\\ for two real symmetric matrices}\label{PDF}

As explained above, in contrast with the unitary or the symplectic groups, there is no closed expression for the orbital integral on the orthogonal group and the formula
\begin{gather}\label{pdfsym} \p(\gamma|\alpha,\beta)={\rm const} |\Delta(\gamma)|
\int \prod_{i=1}^{n} {\rm d}x_i |\Delta(x)| \CI_1(x,\ii \alpha) \CI_1(x,\ii \beta) \CI_1(x,-\ii \gamma) \end{gather}
remains {intractable}. We thus have to resort to a different approach. We first trade the PDF of eigenvalues $\gamma_i$ of the sum $C=A+B$ for the PDF {$\rho(p,q)$} of symmetric functions of the $\gamma$'s, see~(\ref{pdf1}) below. The relation between the two is given below in~(\ref{PDFgamma}).

\subsection[The support of the symmetric functions $P$ and $Q$]{The support of the symmetric functions $\boldsymbol{P}$ and $\boldsymbol{Q}$}\label{PandQ}

For $3\times 3 $ real symmetric {\it traceless} matrices $A$ and $B$, the characteristic polynomial of their sum~$C$ reads\footnote{The assumption of tracelessness is harmless: one may always translate $A$ and $B$ by a multiple of the identity matrix so as to enforce it, at the price of shifting the eigenvalues of their sum $C$ by a common real number.}
\begin{gather*} 
\det(z \I - C)=\det(z \I - A-B)= z^3 + Pz +Q. \end{gather*}
The support of the eigenvalues $\gamma_i$ of $C$ is a convex polygon $\bH_{\alpha,\beta}$ in $\R^2$, defined in terms of the eigenvalues $\alpha_i$ and $\beta_i$ of $A$ and $B$ by the same Horn inequalities as in the Hermitian case~\cite{Fu}:
\begin{gather}
\gamma_{3\min}:=\alpha_3+\beta_3 \le \gamma_3 \le \min(\alpha_1+\beta_3,\alpha_2+\beta_2,\alpha_3+\beta_1)=:\gamma_{3\max},\nonumber\\
 \gamma_{2\min}:=\max(\alpha_2+\beta_3,\alpha_3+\beta_2) \le \gamma_2 \le \min(\alpha_1+\beta_2,\alpha_2+\beta_1)=:\gamma_{2\max}, \nonumber\\
\gamma_{1\min}:=\max(\alpha_1+\beta_3,\alpha_2+\beta_2,\alpha_3+\beta_1) \le \gamma_1 \le \alpha_1+\beta_1=:\gamma_{1\max}, \nonumber\\
\text{as well as (by convention)}\quad \gamma_3\le \gamma_2 \le \gamma_1.\label{inequ3}
\end{gather}

\begin{figure}[!tbp] \centering
 \includegraphics[width=16pc]{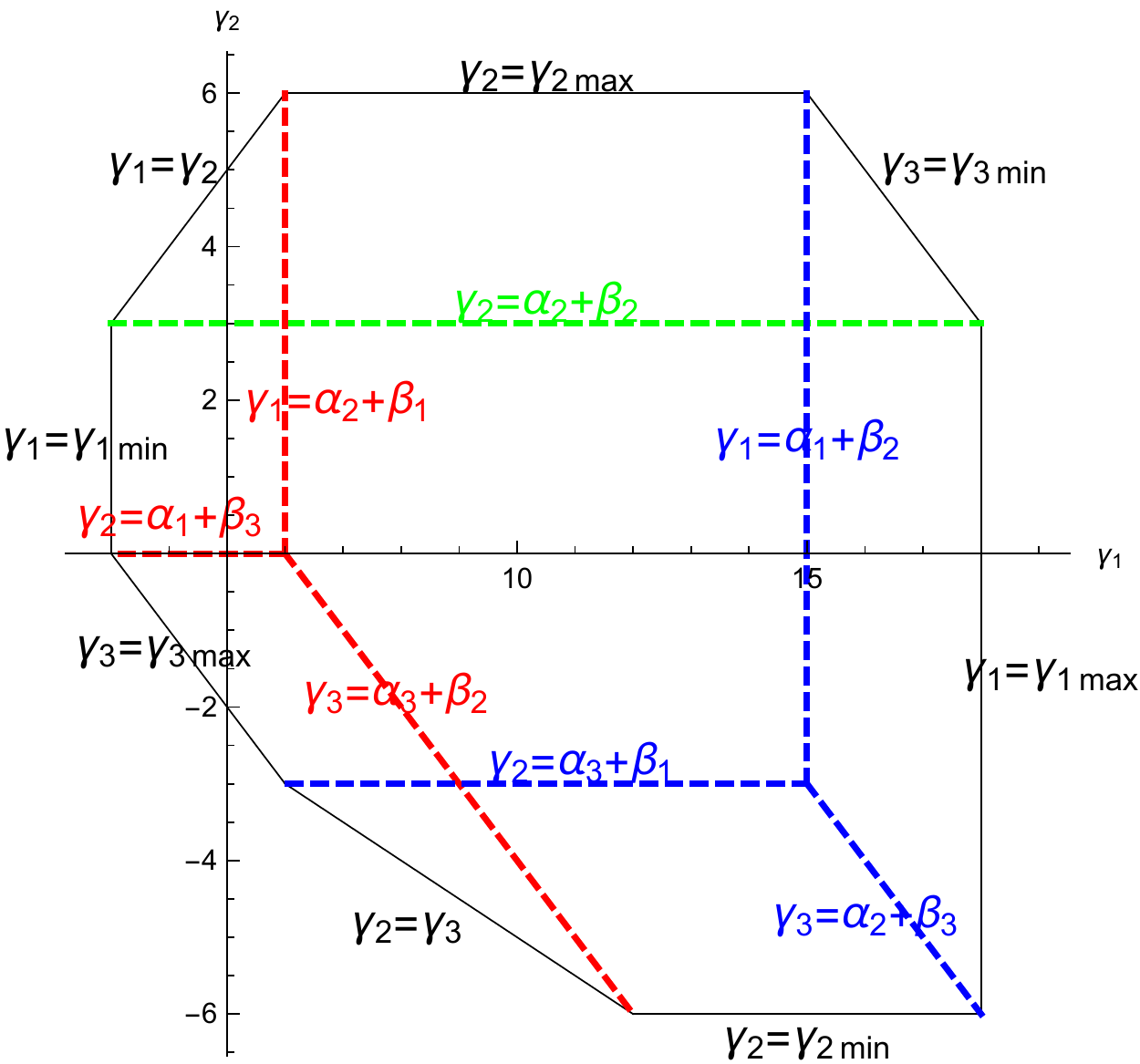} \qquad \includegraphics[width=16pc]{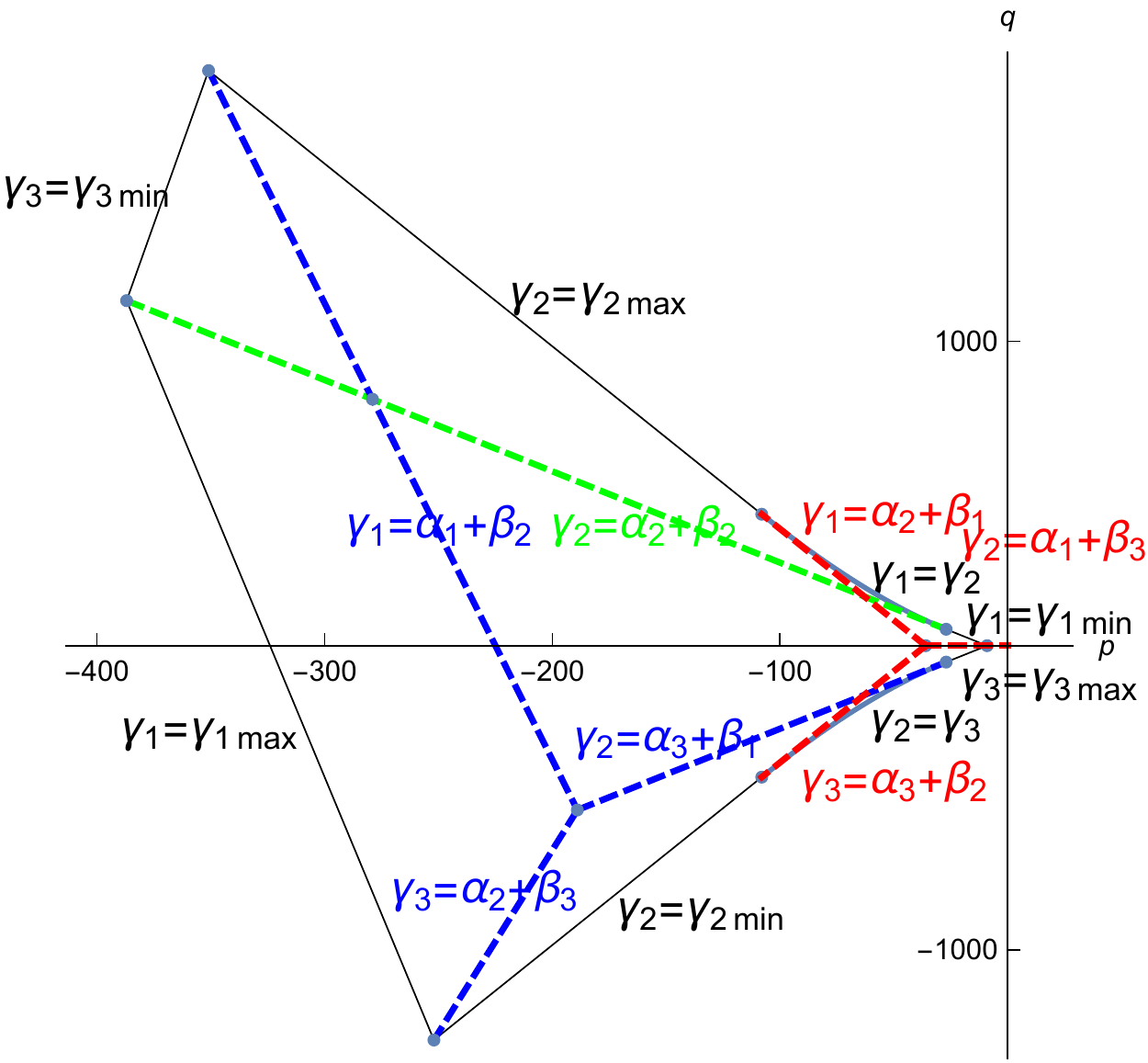}
 \caption{The Horn-tensor polygon $\bH_{\alpha\beta}$ in the $(\gamma_1,\gamma_2)$ plane, and the curvilinear polygon in the $(p,q)$ plane, drawn here for $\alpha=\{11, -1, -10\}$, $\beta=\{7, 4, -11\}$. The dashed (blue, red and green) lines are the expected loci of singularities of the PDF. (A histogram of eigenvalues $\gamma$ for an equivalent configuration of $(\alpha, \beta)$ has
 appeared in \cite[Fig.~7]{Z1}.)}\label{polygongamma}
\end{figure}

This translates for the symmetric functions $P=\gamma_1\gamma_2+\gamma_2\gamma_3+\gamma_3\gamma_1$ and $Q=-\gamma_1 \gamma_2 \gamma_3$ into a curvilinear polygon, whose sides are either segments of lines $cp+q=-c^3$ (for sides of $\bH_{\alpha,\beta}$ where some $\gamma_i=c$) or arcs of the cubic $4p^3+27q^2=0$ for sides of $\bH_{\alpha,\beta}$ where $\gamma_i=\gamma_j$. See Fig.~\ref{polygongamma} for an example with $\alpha=\{11, -1, -10\}$, $\beta=\{7, 4, -11\}$.

It has been suggested~\cite{Z1} that the Horn problem for real symmetric matrices not only shares the same support as in the Hermitian case (for given $\alpha$ and $\beta$), as proved by Fulton~\cite{Fu}, but has also singularities at the same locations, although the former singularities look much stronger. Thus we expect singularities to occur for $3\times 3 $ matrices along the same lines (or half-lines) as in~(\ref{loci-sing}).\footnote{That non $C^\infty$ differentiability occurs only on these lines has now been established in full generality, following an argument of M.~Vergne, see~\cite{CMcSZ}.} These lines are illustrated in Fig.~\ref{polygongamma} for $\alpha=\{11, -1, -10\}$, $\beta=\{7, 4, -11\}$ in the $(\gamma_1,\gamma_2)$ and in the $(p,q)$ planes.

\def\Ro{{\mathcal R}}
\subsection[The statistics of the symmetric functions $P$ and $Q$]{The statistics of the symmetric functions $\boldsymbol{P}$ and $\boldsymbol{Q}$}
For $3\times 3 $ real {\it traceless} symmetric matrices $A$ and $B$, of respective eigenvalues $\alpha$ and $\beta$, the characteristic polynomial of their sum $C$ reads
\begin{gather*} \det(z \I_3 - C)=\det\big(z \I_3- \diag(\alpha)- \Ro \diag(\beta) \Ro^T\big)= z^3 + P(\Ro) z +Q(\Ro).\end{gather*}
For given $\alpha$'s and $\beta$'s, and $\Ro$ regarded as a random variable uniformly distributed in SO(3) (in the sense of the Haar measure), $P(\Ro)$ and $Q(\Ro)$ are also random variables, whose PDF may be written as
\begin{gather} \label{pdf1} \rho(p,q) =\E\big(\delta(P-p) \delta(Q-q)\big) =\int D\Ro \delta(P(\Ro)-p) \delta(Q(\Ro)-q). \end{gather}
We parametrize $\Ro$ in terms of Euler angles
\begin{gather*} \Ro=\Ro_z(\phi) \Ro_y(\theta) \Ro_z(\psi)\end{gather*}
with
\begin{gather*} 0\le \phi\le 2\pi, \qquad 0\le \theta \le \pi,\qquad 0\le \psi\le 2\pi\end{gather*}
and the normalized Haar measure is then
\begin{gather*}
D\Ro= \inv{8\pi^2} \sin \theta {\rm d}\theta {\rm d} \phi {\rm d} \psi.\end{gather*}
The differences $P_p:=P-p$ and $Q_q:=Q-q$ are then degree 2 polynomials of the variable $c:=\cos \theta$, of the form
\begin{gather}
P_p(c) = -{\alpha_1}^2-\alpha_1 \alpha_2-{\alpha_2}^2-{\beta_1}^2-{\beta_1} {\beta_2}-{\beta_2}^2+ \alpha_1 (\beta_2-\beta_1) +\alpha_2(\beta_1+2\beta_2)\nonumber \\
\hphantom{P_p(c) =}{} + x^2 (\alpha_1 - \alpha_2) (\beta_2 + 2 \beta_1) + y^2 (\alpha_2 + 2 \alpha_1) (\beta_1 - \beta_2) - x^2 y^2 ( \alpha_1 - \alpha_2) (\beta_1 - \beta_2) -\blue{p}\nonumber
\\
\hphantom{P_p(c) =}{} + 2 x y ({\alpha_1}-{\alpha_2}) ({\beta_1}-{\beta_2}) \sin \phi \sin \psi \red{c}\nonumber \\
\hphantom{P_p(c) =}{} - \big(\alpha_1+2 {\alpha_2}+{\alpha_1} x^2-{\alpha_2} x^2\big) \big({\beta_1}+2 {\beta_2}+{\beta_1} y^2-{\beta_2} y^2\big){\red{c^2}},\label{PpQq}
\\
 Q_q(c) = (\alpha_1 + \beta_1) (\alpha_1 + \alpha_2 - \beta_2) (\alpha_2 - \beta_1 - \beta_2)+ (\alpha_1 - \alpha_2) (\alpha_1 + \alpha_2 - \beta_2) (2 \beta_1 + \beta_2) x^2 \nonumber\\
\hphantom{Q_q(c) =}{} -(2 \alpha_1 + \alpha_2) (\alpha_2 - \beta_1 - \beta_2) (\beta_1 - \beta_2) y^2\nonumber \\
\hphantom{Q_q(c) =}{} -(\alpha_1 - \alpha_2) (\beta_1 - \beta_2) (\alpha_1 + \alpha_2 + \beta_1 + \beta_2) x^2 y^2 -\blue{q} \nonumber \\
\hphantom{Q_q(c) =}{} + 2 x y (\alpha_1 - \alpha_2) (\beta_1 - \beta_2) (\alpha_1 + \alpha_2 + \beta_1 + \beta_2) \sin \phi \sin\psi \red{c}\nonumber \\
\hphantom{Q_q(c) =}{} +\big((\alpha_1 + 2 \alpha_2) (\alpha_1 + \beta_1) (\beta_1 + 2 \beta_2)
-(\alpha_1 - \alpha_2) (\alpha_1 + \alpha_2 - \beta_1) (\beta_1 + 2 \beta_2) x^2 \nonumber \\
\hphantom{Q_q(c) =}{} + (\alpha_1 + 2 \alpha_2) (\alpha_1 - \beta_1 - \beta_2) (\beta_1 -\beta_2) y^2 \nonumber \\
\hphantom{Q_q(c) =}{} -(\alpha_1 - \alpha_2) (\beta_1 - \beta_2) (\alpha_1 + \alpha_2 + \beta_1 + \beta_2)x^2 y^2\big) {\red{c^2}},\nonumber
 \end{gather}
 where $x:=\cos \phi$ and $y:=\cos \psi$.

Note the peculiar feature of these polynomials in $c$: their degree 0 and 2 terms are (degree 1) polynomials in the variables $x^2$ and $y^2$, while their degree 1 term in $c$ is of the form $x y\sin \phi \sin \psi$, up to a ($\alpha$- and $\beta$-dependent) factor.

Note also that these expressions are $\pi$-periodic in $\phi$ and $\psi$, making it possible to restrict these angles to the interval $(0,\pi)$ where their sine is non negative. This will be implicit in the following.
Thus
\begin{gather*} 
\rho(p,q)=\inv{2\pi^2} \int_0^{\pi} {\rm d}\phi \int_0^{\pi} {\rm d}{\psi} \int_{-1}^1 {\rm d}c \delta(Q_q) \delta(P_p).
\end{gather*}
The PDF for the independent variables $\gamma_1$, $\gamma_2$ then follows simply:
\begin{gather} \label{PDFgamma} \p(\gamma_1,\gamma_2)= |\Delta(\gamma)| \rho(p,q).\end{gather}

\subsection[Reducing $\delta(Q_q) \delta(P_p)$ to $\delta(R)$]{Reducing $\boldsymbol{\delta(Q_q) \delta(P_p)}$ to $\boldsymbol{\delta(R)}$}\label{reducingtoR}

In this subsection, we show that $\int {\rm d}c \delta(Q_q) \delta(P_p)$ may be reduced, up a factor, to a single $\delta(R)$, where $R$ is the {\it resultant} of the two polynomials $P_p(c)$ and $Q_q(c)$. We shall make repeated use of two classical identities \cite{GelfandShilov}:
\begin{enumerate}\itemsep=0pt
\item[--] for $f(t)$ a function with a finite number of ``simple'' zeros $t_i$ (i.e., such that $f'(t_i)\ne 0$),
\begin{gather}\label{ident1}\delta(f(t))=\sum_{t_i} \frac{\delta(t-t_i)}{|f'(t_i)|} ;\end{gather}
\item[--] for $f$ and $g$ two functions with no common zero,
\begin{gather} \label{ident2} \delta(f(t) g(t)) =\frac{\delta(f(t))}{|g(t)|} +\frac{\delta(g(t))}{|f(t)|},\end{gather}
to which we may then apply the previous identity.
\end{enumerate}

Then starting from the product $\delta(Q_q) \delta(P_p)$, we assume that the discriminant $\Delta_Q$ of $Q_q$ is positive, in such a way that the roots $c_{1,2}$ are real and distinct, and we may write
\begin{gather*} \delta(Q_q) =\inv{\sqrt{\Delta_Q}} \big(\delta(c-c_1) +\delta(c-c_2)\big)\end{gather*}
and
\begin{gather*}\int {\rm d}c \delta(Q_q) \delta(P_p)= \inv{\sqrt{\Delta_Q}}\big(\delta(P_p(c_1)) +\delta(P_p(c_2))\big)\end{gather*}
(where it is understood that the deltas act on functions of the remaining variables $\phi$ and $\psi$ or~$x$ and~$y$).

We want to compare this expression with $\delta(R)$ where, as said above, $R$ is the {resultant} of the two polynomials $P_p(c)$ and $Q_q(c)$. If $a$ and $a'$ are the coefficients of terms of degree 2 of the polynomials $Q_q$ and~$P_p$, respectively, and $c'_1,\ c'_2$ the roots of the latter, thus
\begin{gather*} P_p(c)=a' (c-c'_1)(c-c'_2),\qquad Q_q(c)=a(c-c_1)(c-c_2),\end{gather*} the resultant $R$ defined as $a^2 a'^2 \prod\limits_{i,j=1,2}(c_i-c'_j)$ may also be written as
\begin{gather*} R=a^2 P_p(c_1) P_p(c_2).\end{gather*}
(For the polynomials of (\ref{PpQq}), this is a quite cumbersome polynomial of degree {4} in $u=x^2$ and in $z=y^2$, with 4089 $\alpha$- and $\beta$-dependent terms${}^{\ref{repfootnote}}$.)

According to (\ref{ident2}), one writes
\begin{gather} \label{deltaR} \delta(R)= \inv{a^2} \left(\frac{\delta(P_p(c_1))}{|P_p(c_2)|}+\frac{\delta(P_p(c_2))}{|P_p(c_1)|}\right).\end{gather}
But $P_p(c_{1,2})$ has the general form
\begin{gather*} P_p(c_{1,2})= A \pm B \sqrt{\Delta_Q} \end{gather*}
 with $A$ and $B$ functions of $x$ and $y$. Thus whenever $P_p(c_1)$ vanishes for some $(x,y)$, i.e., $A=-B \sqrt{\Delta_Q} $, we have $P_p(c_2)=-2 B \sqrt{\Delta_Q} $ for those values. And vice versa, if $P_p(c_2)$ vanishes, then $P_p(c_1)=2 B \sqrt{\Delta_Q} $. One may thus rewrite (\ref{deltaR}) as
\begin{gather} \label{deltaR2} \int {\rm d}c \delta(Q_q) \delta(P_p) = \inv{\sqrt{\Delta_Q}}(\delta(P_p(c_1))+\delta(P_p(c_2))) = {2a^2 |B|} \delta(R),\end{gather}
a non-trivial and useful identity.\footnote{This identity may be generalized to higher degree polynomials with real and distinct roots~\cite{BZ19}.}
For the polynomials $P_p$ and $Q_q$ of~(\ref{PpQq}), it is easy to compute that $2 a^2 B=(x y\sin \phi \sin \psi ) \tilde B$, with
\begin{gather*}
\tilde B = 2(\alpha_1-\alpha_2)(\beta_1-\beta_2) [(2 \alpha_1 + \alpha_2)( \alpha_1 + 2\alpha_2) ((\beta_1 - \beta_2) z + (\beta_1 + 2 \beta_2))\nonumber\\
\hphantom{\tilde B =}{} + (2 \beta_1+ \beta_2) ( \beta_1+ 2\beta_2) ((\alpha_1 - \alpha_2) u + (\alpha_1 + 2 \alpha_2) )],
\end{gather*}
while the prefactor $|x| |y|\sin \phi \sin \psi $ (recalling that $x=\cos \phi$, $y=\cos \psi$) enables us to change variables to $u=x^2$, $z=y^2$, with the result that
\begin{gather}\label{pdf2a} \rho(p,q)= \inv{{2}\pi^2} \int_0^1 {\rm d}z \int_0^1 {\rm d}u |\tilde B| \delta(R).\end{gather}
The final transformation of this expression follows from (\ref{ident1}). If $u_i(z)$ denote {those} roots of the polynomial $R(u,z)$ that belong to the interval $[0,1]$, we may write
\begin{gather}\label{pdf2b} \rho(p,q)=\inv{2\pi^2} \int_0^1 {\rm d}z \sum_{i} \frac{|\tilde B(u_i,z)|}{|R'_u(u_i,z)|}.\end{gather}
This integral will be studied more explicitly in a particular example in the next section.

Three remarks are in order:
\begin{enumerate}\itemsep=0pt
\item[--] we have assumed from the start that the discriminant $\Delta_Q$ is positive, and this led us to~(\ref{pdf2a}). Conversely, the vanishing of $R$ for real values of the variables $u$ and $z$ encompassed in~(\ref{pdf2a}) implies that $P_p(c)$ and $Q_q(c)$ have a common root, and this may only be possible if that root is a root of $a P_p -a' Q_q$ which is a degree 1 polynomial in $c$ with real coefficients. Thus the common roots of $P_p(c)$ and $Q_q(c)$ are necessarily real, justifying our assumption that $\Delta_Q{\ge} 0$;
\item[--] the roots $c_1$ and $c_2$ have to lie in $[-1,1]$ for the consistency of the derivation;
\item[--] the two functions $P_p(c_1)$ and $P_p(c_2)$ have to have no common zero. Otherwise if that happened at some values $(x,y)$, both $A$ and $B$ would vanish.
\end{enumerate}
The two latter points will be verified in the particular case that we discuss now.

\section[The particular case $A\sim B\sim J_z$]{The particular case $\boldsymbol{A\sim B\sim J_z}$} \label{sectionparticularcase}
\subsection{Symmetries and reduction to algebraic equations}
From now on, we restrict ourselves to the particular case of $\alpha=\beta=\{1,0,-1\}$. The polyno\-mials~$P_p$ and~$Q_q$ then reduce to
\begin{gather*}
P_p(c)= P(R)-p = - c^2 \big(1+x^2 \big) \big(1+y^2 \big)+2c x y \sin \phi \sin \psi - \big(x^2 y^2 -2x^2 -2 y^2 +3\big) -p, \\
Q_q(c)=Q(R)-q= 2c^2 \big(1- x^2 y^2\big)+4 c x y \sin \phi \sin \psi -2 \big(1-x^2\big) \big(1-y^2\big) -q,
\end{gather*}
where we recall that $c=\cos \theta$, $x=\cos \phi$, $y=\cos \psi$, and we could substitute $\sin \phi=\sqrt{1-x^2}$, $\sin\psi=\sqrt{1-y^2}$ since $\phi$ and $\psi$ are restricted to $(0,\pi)$. Note that these expressions are invariant under the exchange of $\phi$ and $\psi$, or of $x$ and $y$.

\begin{figure}[t] \centering
\includegraphics[width=15pc]{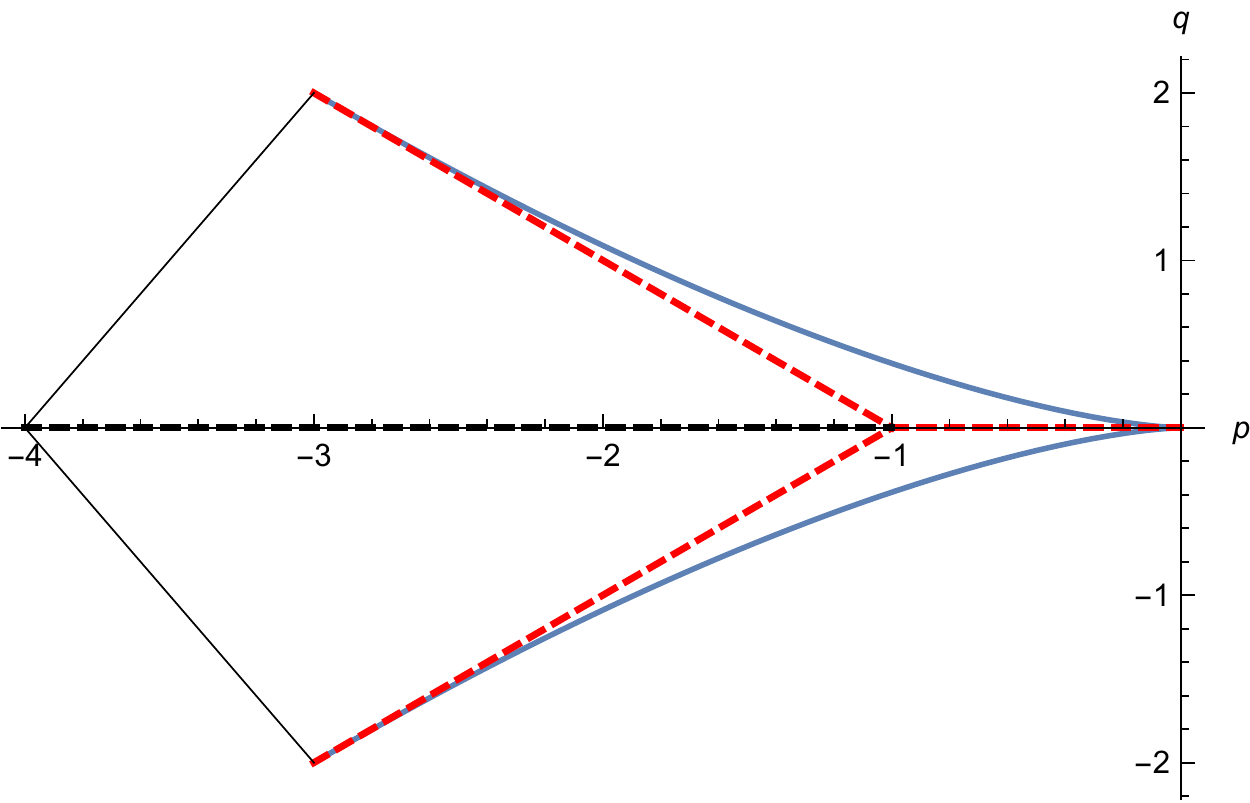} \qquad {\includegraphics[width=15pc]{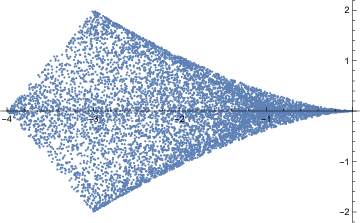}}
 \caption{Left, the curvilinear polygon, bounded by the lines $2p+8\pm q=0$ and two arcs of the cubic $(p/3)^3+(q/2)^2=0$; (right) the distribution of 8000 points in the $(p,q)$ plane for $A=B=\diag(1, 0, -1)$. The (red and black) dashed lines along which singularities are expected are $q=0$ and $p=-1\pm q$.}\label{polygonJzJz}
\end{figure}

Another symmetry of the problem will be quite useful. Because $-J_z$ is {conjugate to $J_z$ by the action of $R_y(\pi/2)$}, it is clear that the distribution of the $\gamma$'s will be invariant under change of sign, and the distribution of the $(p,q)$ variables will be invariant under change of sign of~$q$. This is apparent in Fig.~\ref{polygonJzJz} where we see that both the support of $\rho$, here a curvilinear quadrangle, and the distribution of points (here a simulation with $8000$ points), are symmetric under $q\to -q$. Our formalism, however, does not have that manifest symmetry, and we will find it useful in the following to choose~$q$ {\it negative}, as we see shortly. The PDF for all values of $q$ will then be reconstructed by symmetry.

The discriminant of $Q_q$ is $\Delta_Q=16 \big(1 - x^2\big) \big(1 - y^2\big) +8 q \big(1 - x^2 y^2\big)$. As explained previously, its positivity follows by consistency from the vanishing of the resultant $R$. For $-2\le q\le 0$, one may check that, if $\Delta_Q\ge 0$, the real roots obey $|c_{1,2}|\le 1$ for all $x, y \in [-1,1]$. {\it We shall hereafter assume that $-2\le q\le 0$.}

Following the discussion of Section~\ref{reducingtoR}, we then determine the resultant of the two polyno\-mials~$P_p(c)$ and $Q_q(c)$, a degree 4 polynomial in~$u$ and~$z$
\begin{gather*}
R=\operatorname{Res}(P_p,Q_q; c) = 4 p^2 (1-u z)^2+4p q(1+u)(1+z)(1-u z) \nonumber\\
\hphantom{R=}{} -8p \big({-}4 + 2 (u + z) + (u + z)^2 - 2 u z (u + z) + u z (u - z)^2\big) \nonumber\\
\hphantom{R=}{} +q^2 \big( (1 + u z)^2 + (u + z)^2 + 2 u z (u + z) + 2 (u + z)\big) \\ 
\hphantom{R=}{} +4 q \big( \big(4 - u^3 + 2 u^2 z^2 - z^3 + 2 (u + z) - u z (u + z) - 3 (u + z)^2 + 3 u z \big(u^2 + z^2\big)\big) \big)\nonumber \\
\hphantom{R=}{} +4 \big(u^4-8 u^3 z-2 u^2 z^2-8 u z^3+z^4\big)+16 (u+z)^3-16 (u-z)^2+64 (-u-z+1). \nonumber
 \end{gather*}
Of course, this polynomial $R$ is also symmetric under the swapping of $u$ and $z$, since $P$ and $Q$ were under $x\leftrightarrow y$. The factor $\tilde B$ appearing in (\ref{pdf2a}) is $\tilde B=4 \big(2+x^2+y^2\big)$. Then, according to~(\ref{pdf2b}), the PDF reads
\begin{gather}\label{pdf4} \rho(p,q)= \frac{2}{\pi^2}\int_{0}^1 {\rm d}z \sum_{\mathrm{roots}\ u_i \ {\rm of}\ R\atop 0\le u_i \le 1} \frac{(2+u+ z)}{|R'_u| }\Big|_{u=u_i}.\end{gather}

\subsection[Roots $u_i$ of the resultant $R$ and their singularities]{Roots $\boldsymbol{u_i}$ of the resultant $\boldsymbol{R}$ and their singularities}\label{subsectionroots}

Let us start with some general features of the roots $u_i(z)$ of $R$:
\begin{enumerate}\itemsep=0pt
\item[--] The polynomial $R(u,z)$ being symmetric in $u$, $z$, its roots $u_i(z)$ or $z_i(u)$ are built by the same function: $z\buildrel {\varsigma_i} \over{\mapsto}u_i(z)$, $u\buildrel {\varsigma_i} \over{\mapsto}z_i(u)$. This means that their graph is symmetric with respect to the first diagonal, see Fig.~\ref{gallery} below.
 \item[--] Within the full domain $q\le 0$, there are either 0, two or four roots $u_i(z)\in [0,1]$ for $ z\in [0,1]$. When $z$ varies in $(0,1)$, this number may change: either some of these roots may evade the interval $[0,1]$ through one of its end points 0 or 1, but they always do it pairwise; or a pair of real roots ``pops out'' of the complex plane or disappears into it, but this may occur only when they coalesce. In both cases, the discriminant $\Delta_R$ of $R$ with respect to $u$ vanishes. We compute
\begin{gather*}
\Delta_R(z,p,q)= 8192 (8 + 2 p - q)^2 (1 - z)^2 z^2 \Delta^{(1)}_R(z,p,q),\end{gather*}
where $ \Delta^{(1)}_R(z,p,q) $ is a fairly cumbersome polynomial of degree 8 in $z$ that we refrain from writing here${}^{\ref{repfootnote}}$. The vanishing of the first factor does not occur for $q<0$ and $p> -4$. In the following, we denote the ordered roots $z_s$ of $\Delta_R$ belonging to $[0,1]$ as
\begin{gather*} 
z_{s_0}=0 \le z_{s_1} \le \cdots \le z_{s_r} \le 1. \end{gather*}
It turns out there are up to five roots $z_s$ of $ \Delta^{(1)}_R(z,p,q) $ in the open interval $]0,1[$, hence seven in $[0,1]$. Some of these roots $z_s$ may be {\it irrelevant}, in the sense that they are associated with the merging of {\it irrelevant roots} $u_i(z)$ of $R$, i.e., roots that {\it do not} satisfy $u_i\in[0,1]$.
\end{enumerate}

\begin{figure}[t] \centering
\includegraphics[width=20pc]{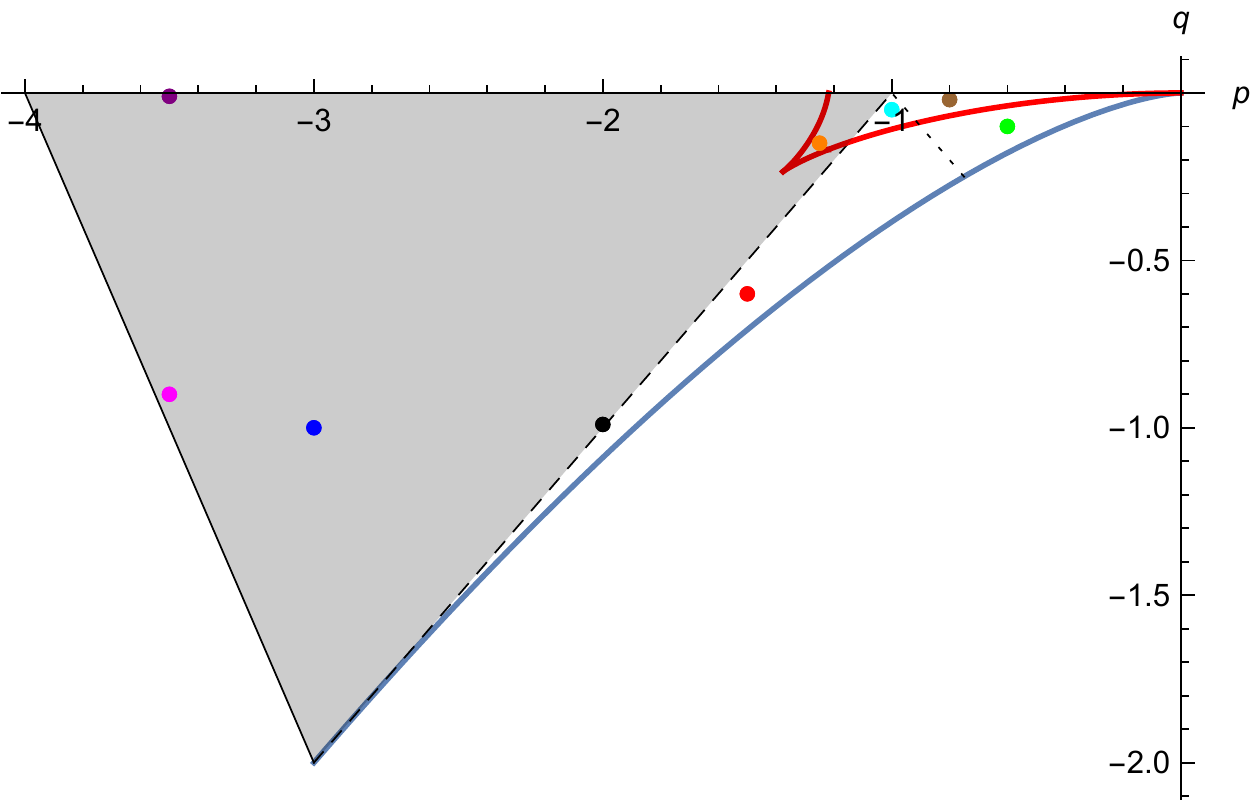}
 \caption{The $q\le 0$ part of the domain in the $(p,q)$-plane. The special curves are the lines $p-q+1=0$ (dashed) and $p+q+1=0$ (dotted), and an arc of the ``horned'' curve where $T=0$ (red). The colored dots refer to the plots of Fig.~\ref{gallery}.}\label{horneddomain}
\end{figure}

At some particular values of $(p,q)$, the number of $z_s$ roots may change. Either some $z_s$ evade the interval $[0,1]$ or enter it, through 0 or 1: computing $ \Delta^{(1)}_R(0,p,q) $ and $ \Delta^{(1)}_R(1,p,q) $, one finds that this happens along the lines $p\pm q +1=0$ and $q=0$ which are the dashed lines in Fig.~\ref{polygonJzJz}. Or two roots of $ \Delta^{(1)}_R(z,p,q) $ coalesce, and this occurs for values of $(p,q)$ that are
roots of the discriminant $\Delta_{R}^{(2)}$ of $\Delta^{(1)}_R(z,p,q) $ with respect to~$z$,
\begin{gather*} \Delta_{R}^{(2)}(p,q) = q^4 (2-q)^2 (1+p-q)^2 (4+p-q)^2 (5+p-q)^2 (8+2 p-q)^2 \big(4 p^3+27 q^2\big)\nonumber \\
\hphantom{\Delta_{R}^{(2)}(p,q) =}{} \times\big(64+16 p-48 q - 8 p q+13 q^2+p q^2-q^3\big)^2 \nonumber\\
\hphantom{\Delta_{R}^{(2)}(p,q) =}{} \times \big(8+10 p+2 p^2-19 q-10 p q-p^2 q+8 q^2+2 p q^2-q^3\big)^2 (T(p,q))^3,
\end{gather*}
where $T(p,q)$ is a horrendous polynomial of degree 21 in $p$ and 18 in $q$.${}^{\ref{repfootnote}}$ The relevant roots of~$\Delta_{R}^{(2)}$ (in fact of~$T$) for our discussion define the ``horned'' (red) curve in the upper right part of Fig.~\ref{horneddomain}. The cusp of that curve occurs at $(p_c,q_c)=(-1.37657\ldots,-0.234765\ldots)$. Note that this $p_c$ is a root of the ``third generation discriminant" $\Delta^{(3)}_R(p):= \Delta_T$, in fact of its factor $\big(2 - 3 p - p^2 + 3 p^3 + p^4\big)$. The horned curve intersects the $q$-axis at $p=p_0:=-1.21891\ldots$, a~root of $1328 + 1325 p + 171 p^2 - 17 p^3 + p^4$, which is a factor of $T(p,0)$. Also, it intersects the dotted line $p+q+1=0$ at $p=-0.910988\ldots$, $q=-0.089012\ldots$; it intersects the dashed line $p-q+1=0$ at $p= -1.14617\ldots$, $q=-0.146174\ldots$.

\begin{figure}[t] \centering
 \includegraphics[width=20pc]{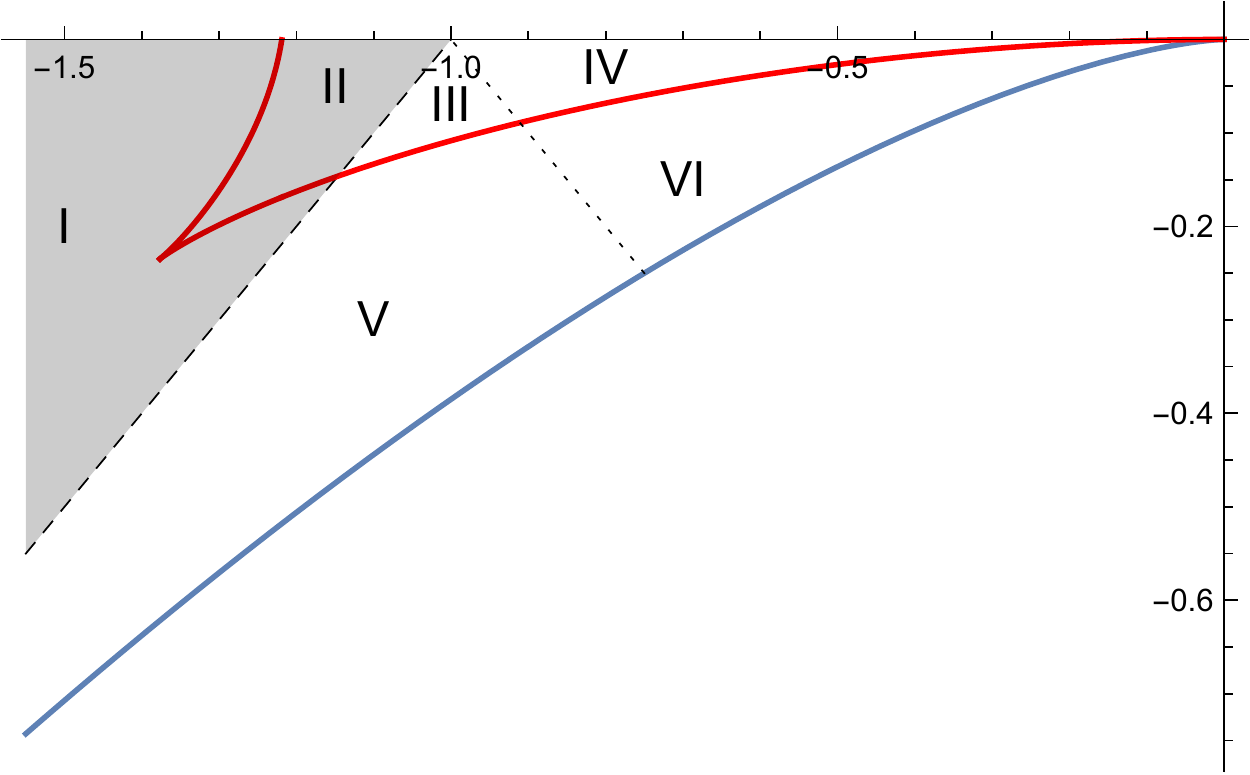}
 \caption{Zoom on the 6 regions.}\label{subdomains}
\end{figure}

\begin{figure}[t]
 \centering\includegraphics[width=8.7pc]{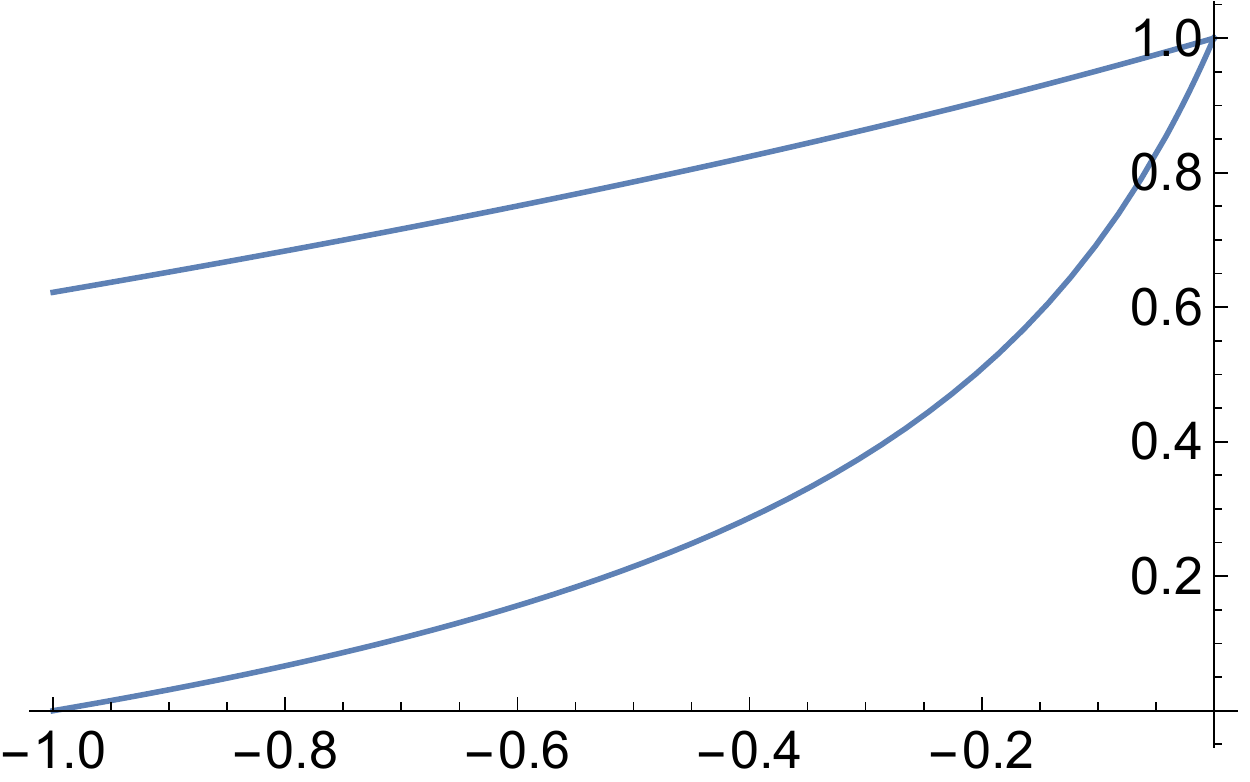}\quad\includegraphics[width=8.7pc]{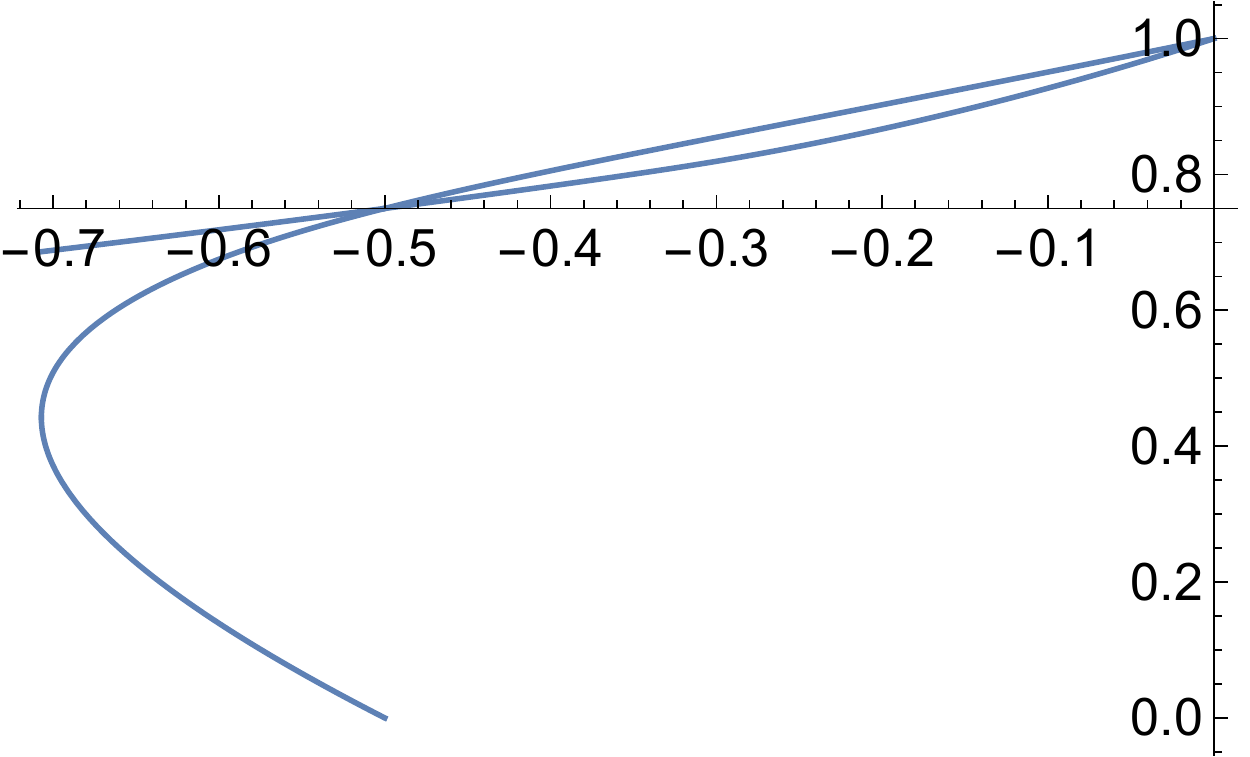}\quad \includegraphics[width=8.7pc]{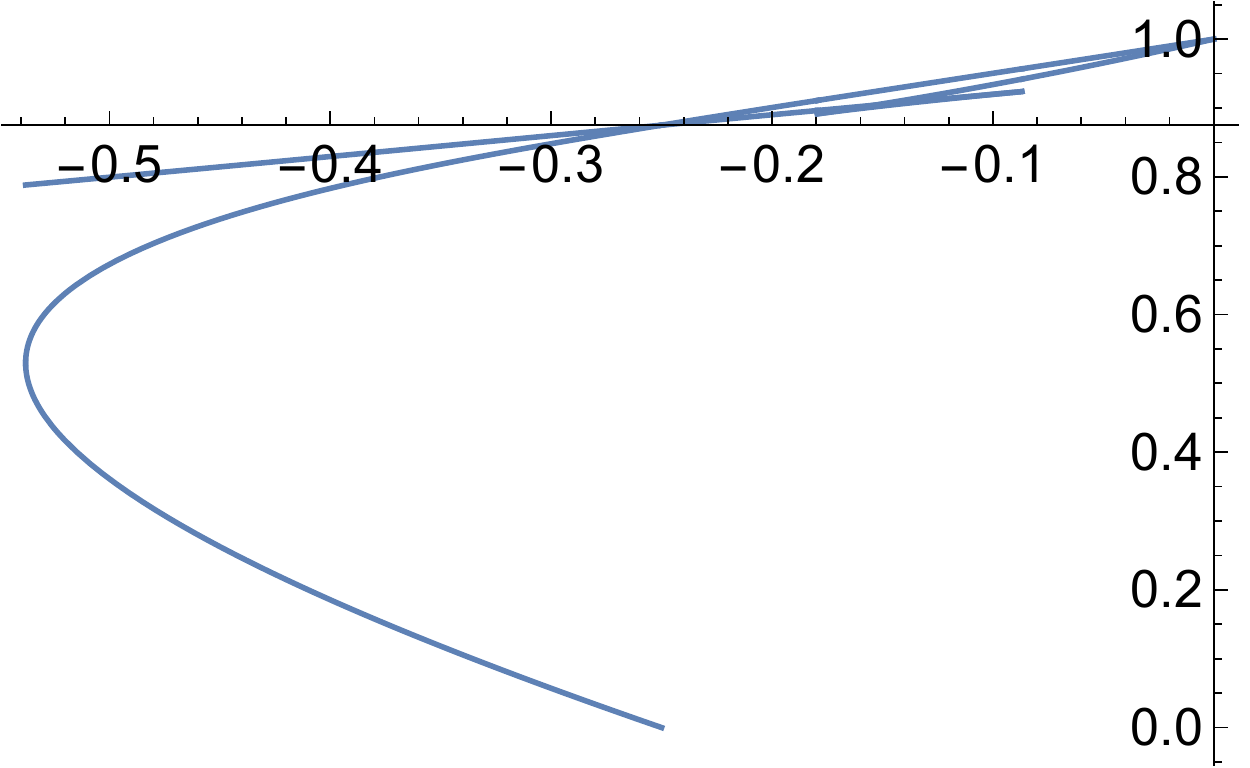}\quad \includegraphics[width=8.7pc]{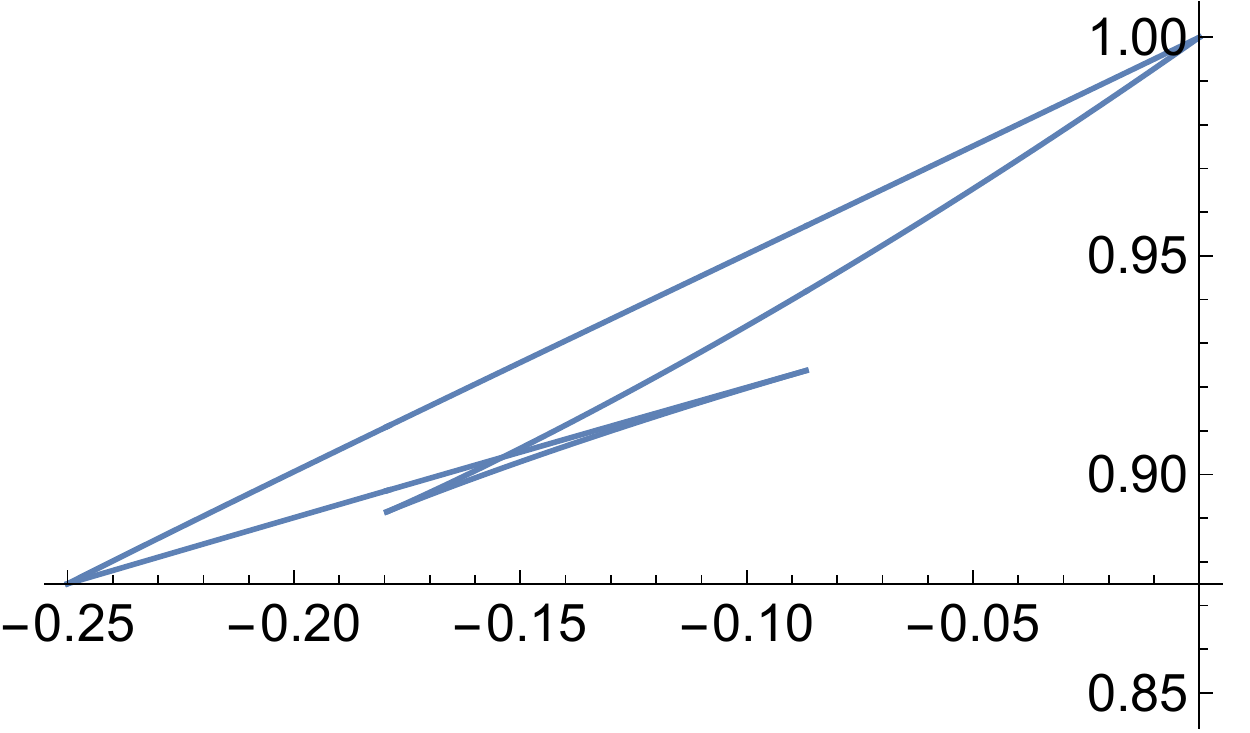}\\
\qquad \qquad \qquad (a)\hfill (b)\hfill (c) \hfill (d)\qquad \qquad\qquad
 \caption{Plot of the roots $z_s$ of $ \Delta^{(1)}_R$ as a function of $q$: (a) for $p=-3.5$; (b) for $p=-1.5$; (c) for $p=-1.25$, for $q_{{\rm min}}(p)\le q \le 0$; (d) zoom of the latter on $-.25\le q\le 0$. The cusps at $q=-0.179503$ and $q=-0.0868984$ lie on the boundary of the horned domain. In each case, the largest $z_s$ is irrelevant.}\label{plot-zs}
\end{figure}

A detailed analysis shows that one has to distinguish six regions in the domain of $(p,q)$, $q\le 0$, see Fig.~\ref{subdomains}. These regions differ by the subset of relevant values $z_s$. Note that we have to carry out the $z$-integration of (\ref{pdf4}) in each interval $(z_{s_i},z_{s_{i+1}})$, after one another, because the integrand is singular at each $z_{s_j}$, as we discuss in the next subsection. The properties of these six regions are summarized in Table 1; the pattern of $z_s$ when $q$ varies while $p$ is fixed at some value are displayed in Fig.~\ref{plot-zs}; and the various scenarios for the roots $u_i(z)$, which describe several branches of a closed curve in the $(z,u)$ plane, are illustrated in Fig.~\ref{gallery}, where colors refer to points of Fig.~\ref{horneddomain}.

\begin{table}[th!]\centering
\begin{tabular}{|c| c|c|c|c|c|c}
\hline
Region &Number of $z_s\in {} ]0,1[$ & Relevant $z_s$ &Intervals& Roots $u_i$ \\
\hline
I& 2 & $z_{s_0}=0, z_{s_1}$ & $0\le z\le z_{s_1}$& $u_1\le u_2$ \\
\hline
& & & $0\le z < z_{s_1}$& $u_1\le u_2$ \\
II & 4 & $ z_{s_0}=0, z_{s_1}, z_{s_2}, z_{s_3} $ & $ z_{s_1}< z < z_{s_2}$& $u_1\le u_3 \le u_4 \le u_2$ \\
 & & & $ z_{s_2}< z < z_{s_3}$& $ u_4 \le u_2$ \\
\hline
& & & $ z_{s_1}\le z < z_{s_2}$& $u_1\le u_2$ \\
 III & 5 & $ z_{s_1}, z_{s_2}, z_{s_3}, z_{s_4} $ & $ z_{s_2}< z < z_{s_3}$& $u_1\le u_3 \le u_4 \le u_2$ \\
 & & & $ z_{s_3}< z < z_{s_4}$& $ u_4 \le u_2$ \\
\hline
& && $ z_{s_2}\le z < z_{s_3}$& $u_1\le u_2$ \\
 IV &5 & $ z_{s_2}, z_{s_3}, z_{s_4}, z_{s_5} $ & $ z_{s_3}< z < z_{s_4}$& $u_1\le u_3 \le u_4 \le u_2$ \\
 & & & $ z_{s_4}< z < z_{s_5}$& $ u_1 \le u_3$ \\
\hline
V& 3 & $z_{s_1}, z_{s_2}$ & $z_{s_1}\le z\le z_{s_2}$& $u_1\le u_2$ \\
\hline
VI& 3 & $z_{s_2}, z_{s_3}$ & $z_{s_2}\le z\le z_{s_3}$& $u_1\le u_2$ \\
\hline
\end{tabular}
\caption{Pattern of roots $z_s$ of $\Delta_R$ and of roots $u_i(z)$ of $R$ in the various regions.}\label{Table1}
\end{table}

For example, in region II, (shaded triangle $\cap$ horned region), there are four roots $z_{s_j} $, $j=1,\dots,4$, of $\Delta_R^{(1)}$ but only the first three are relevant, and four relevant roots for $\Delta_R$, namely $ z_{s_0}=0, z_{s_1}, z_{s_2}, z_{s_3}$:
\begin{enumerate}\itemsep=0pt
\item[$*$] for $0 \le z \le z_{s_1}$, there are two roots $ 0 \le u_1\le u_2 \le 1$;
\item[$*$] at $z=z_{s_1}$, a new pair of roots $(u_3,u_4)$ pops out of the complex plane;
\item[$*$] for $z_{s_1}<z < z_{s_2}$, we have four roots $0 \le u_1 \le u_3 \le u_4 \le u_2 < 1$;
\item[$*$] at $z=z_{s_2}$, the pair $(u_1,u_3)$ merges and disappears into the complex plane;
\item[$*$] for $z_{s_2} \le z \le z_{s_3}$, we are left with two roots $u_4 \le u_2$, which merge at $z= z_{s_3}$ and disappear in the complex plane;
\item[$*$] for $ z> z_{s_3}$, there is no root in the interval $u\in [0,1]$. The $z$-integration must be carried out separately on the three intervals $(0,z_{s_1})$, $(z_{s_1},z_{s_2})$, and
 $(z_{s_2},z_{s_3})$.
 \end{enumerate}

\begin{figure}[h!] \centering
\includegraphics[width=8.7pc]{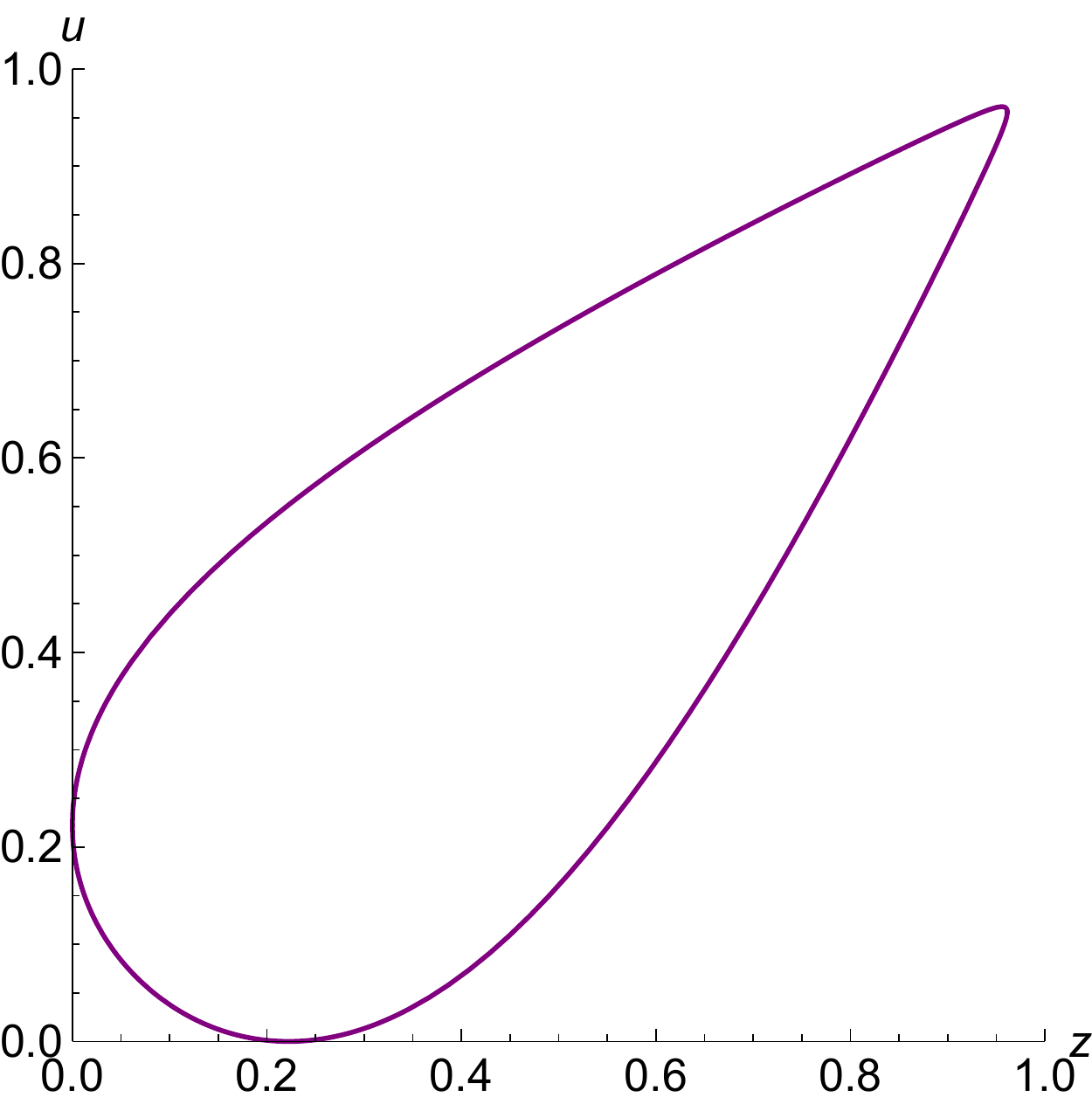}\quad \includegraphics[width=8.7pc]{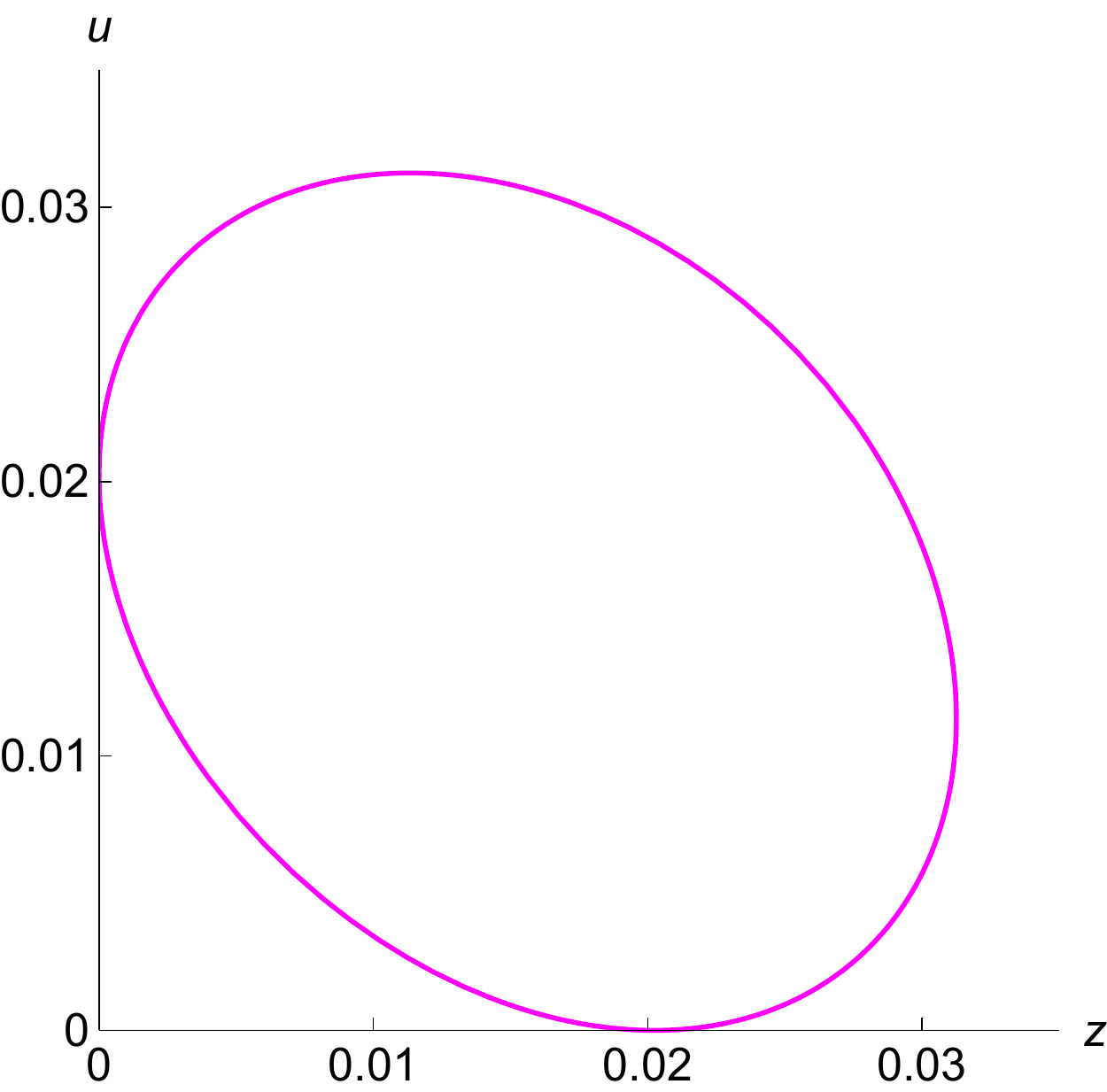}\quad\includegraphics[width=8.7pc]{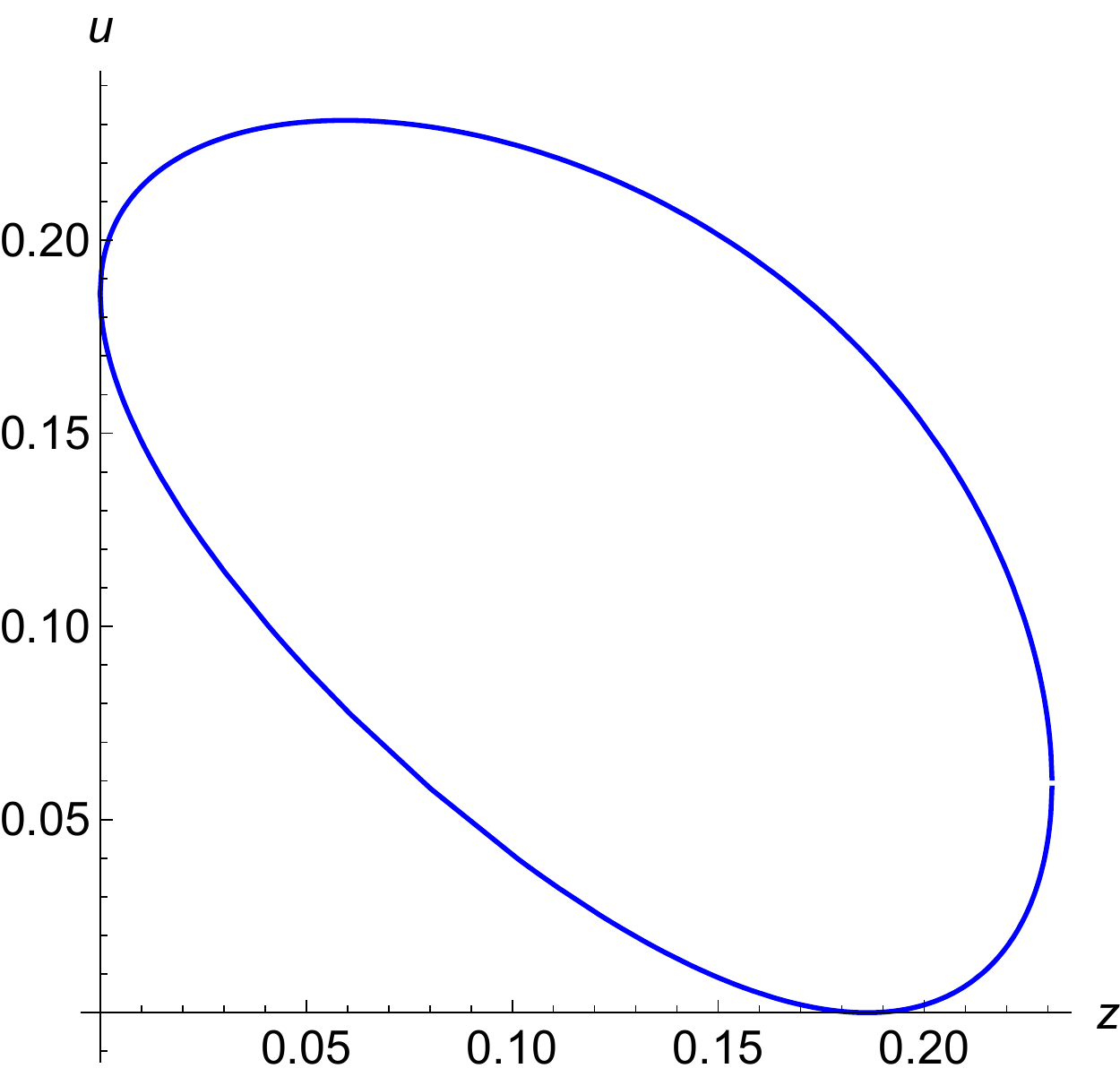}
\quad\includegraphics[width=8.7pc]{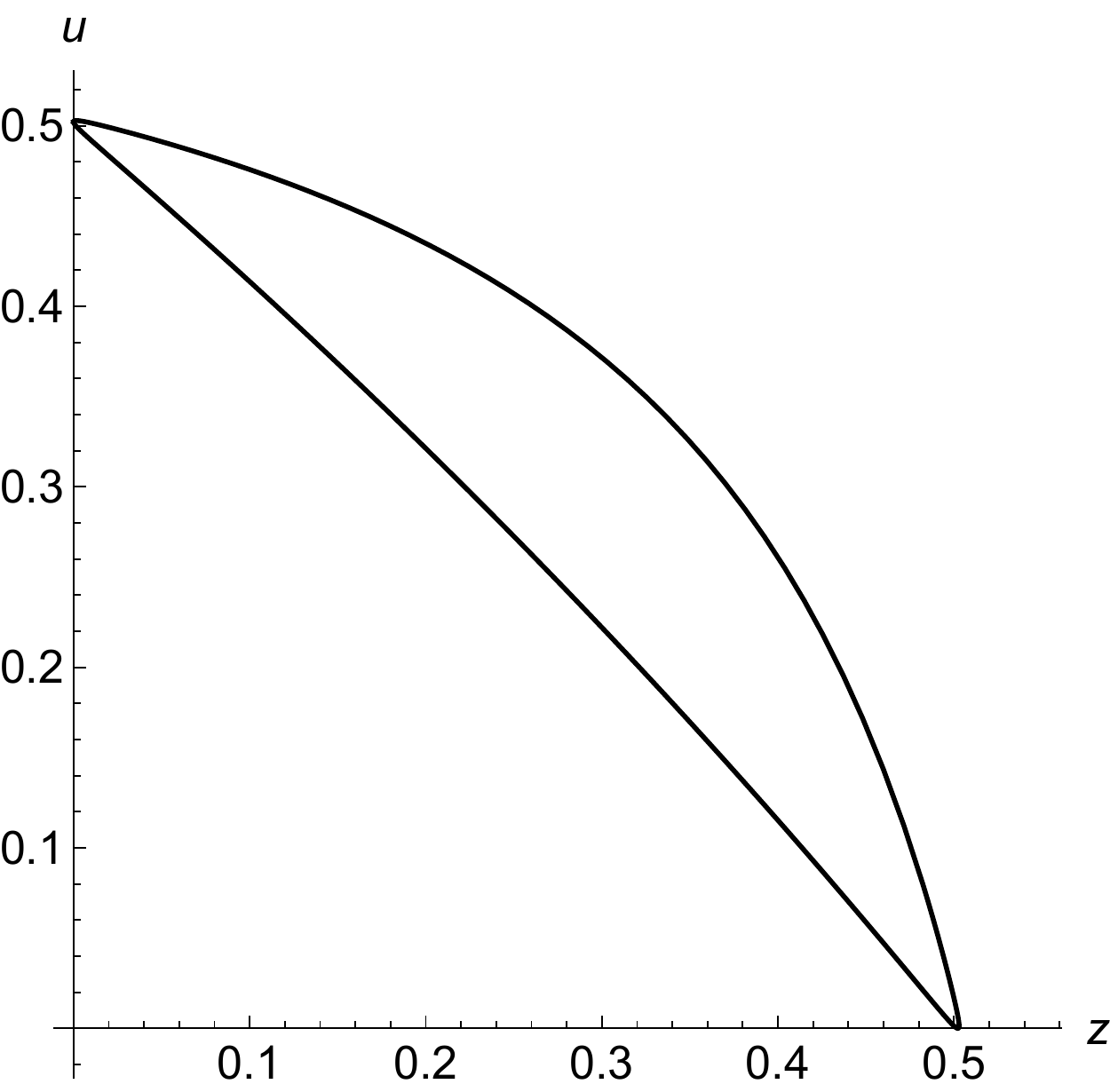}\\
(a): $(-3.5,-.01)$ {\bf (I)}\qquad (b): $(-3.5,-.9)$ {\bf (I)} \qquad (c): $(-3,-1)$ {\bf (I)}\qquad (d): $(-2,-.99)$ {\bf (I)}\\
{\small $z_{s_0}=0,\,z_{s_1}=0.961219$ \quad $z_{s_0}=0,\,z_{s_1}=0.031247$ \quad $z_{s_0}=0,\,z_{s_1}=0.2310119$ \quad $z_{s_0}=0,\,z_{s_1}=0.502881$}
\\[3pt]
\includegraphics[width=8.7pc]{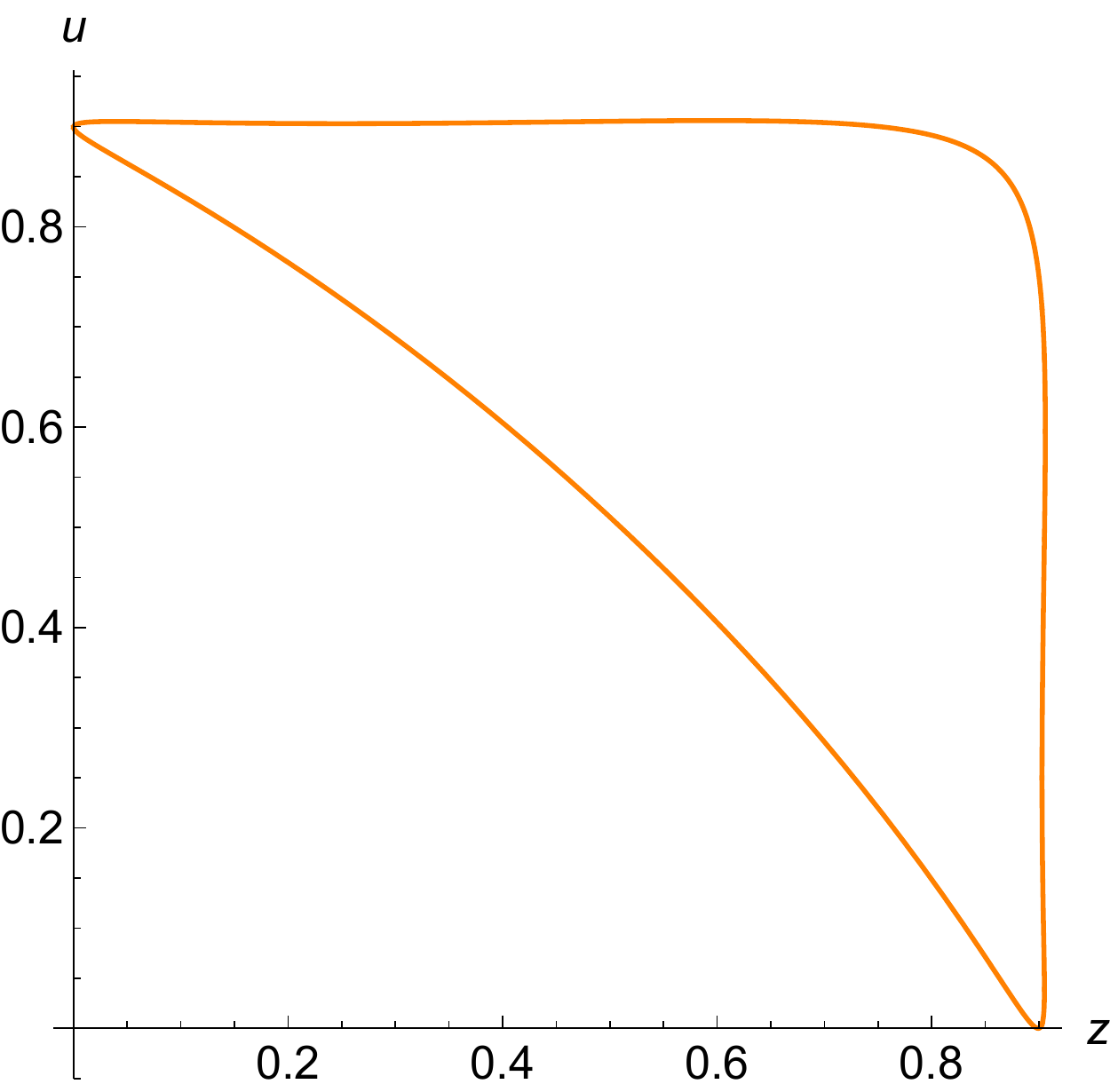}\quad \includegraphics[width=8.7pc]{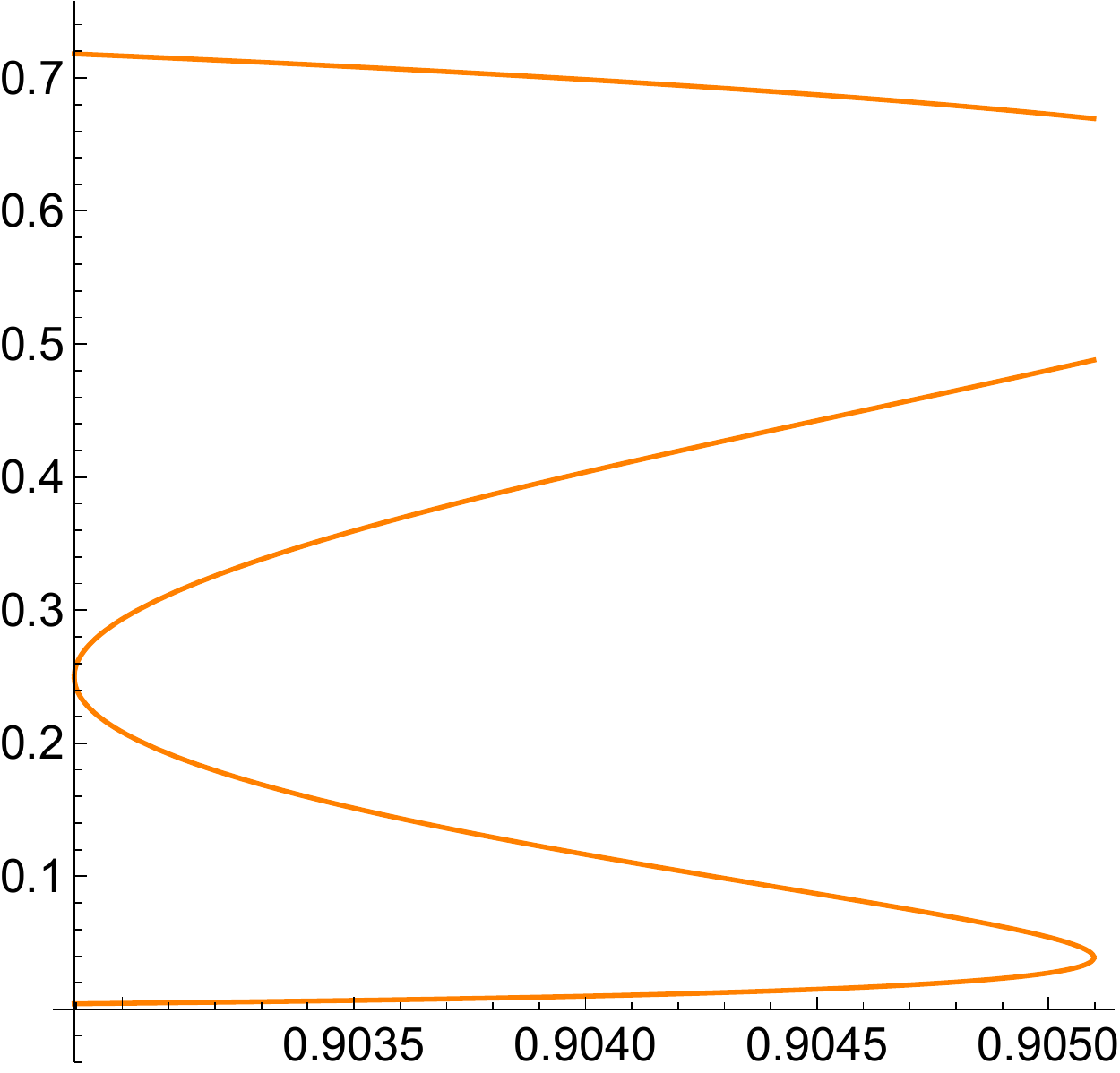}\quad\includegraphics[width=8.7pc]{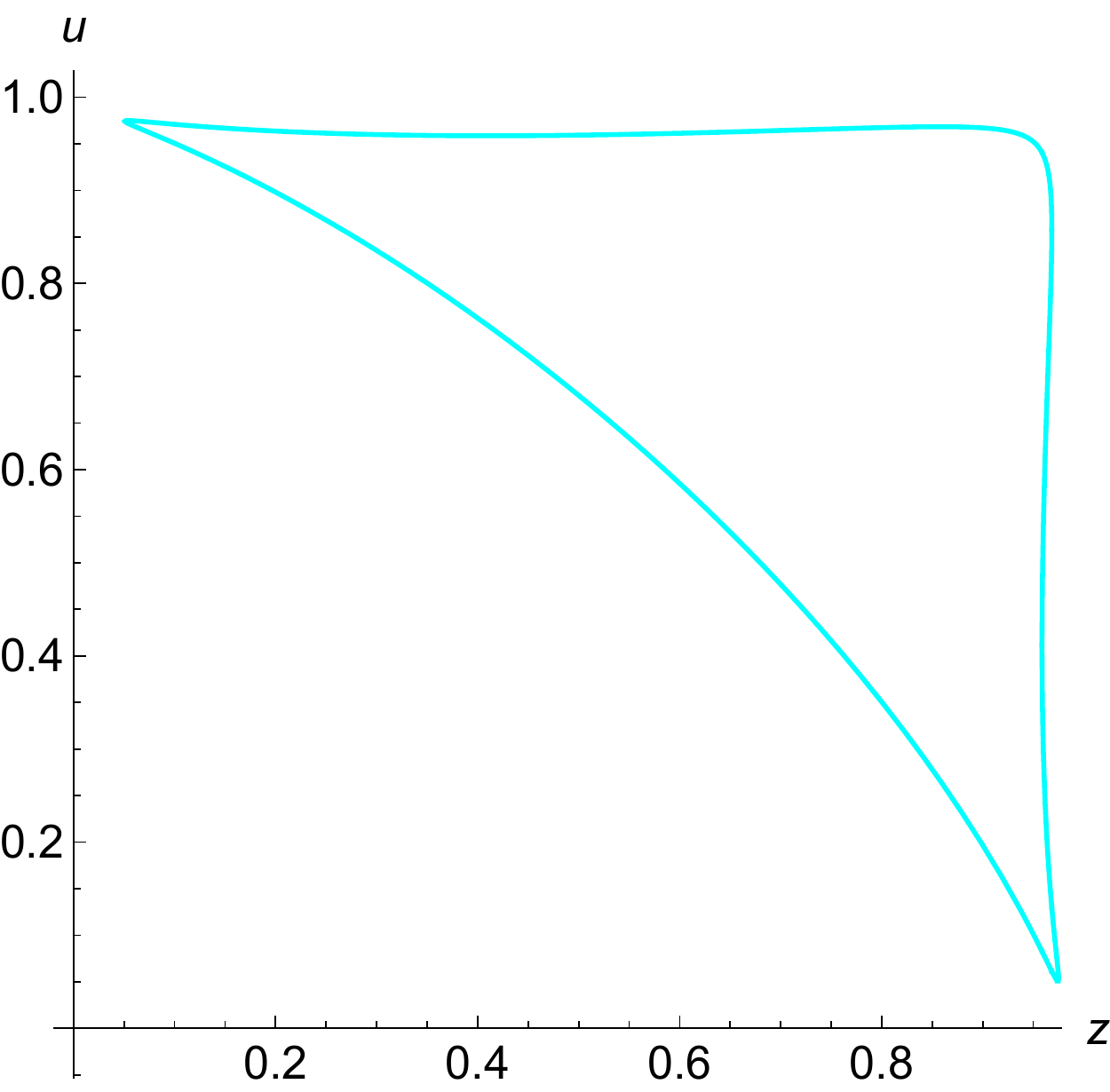}\quad \includegraphics[width=8.7pc]{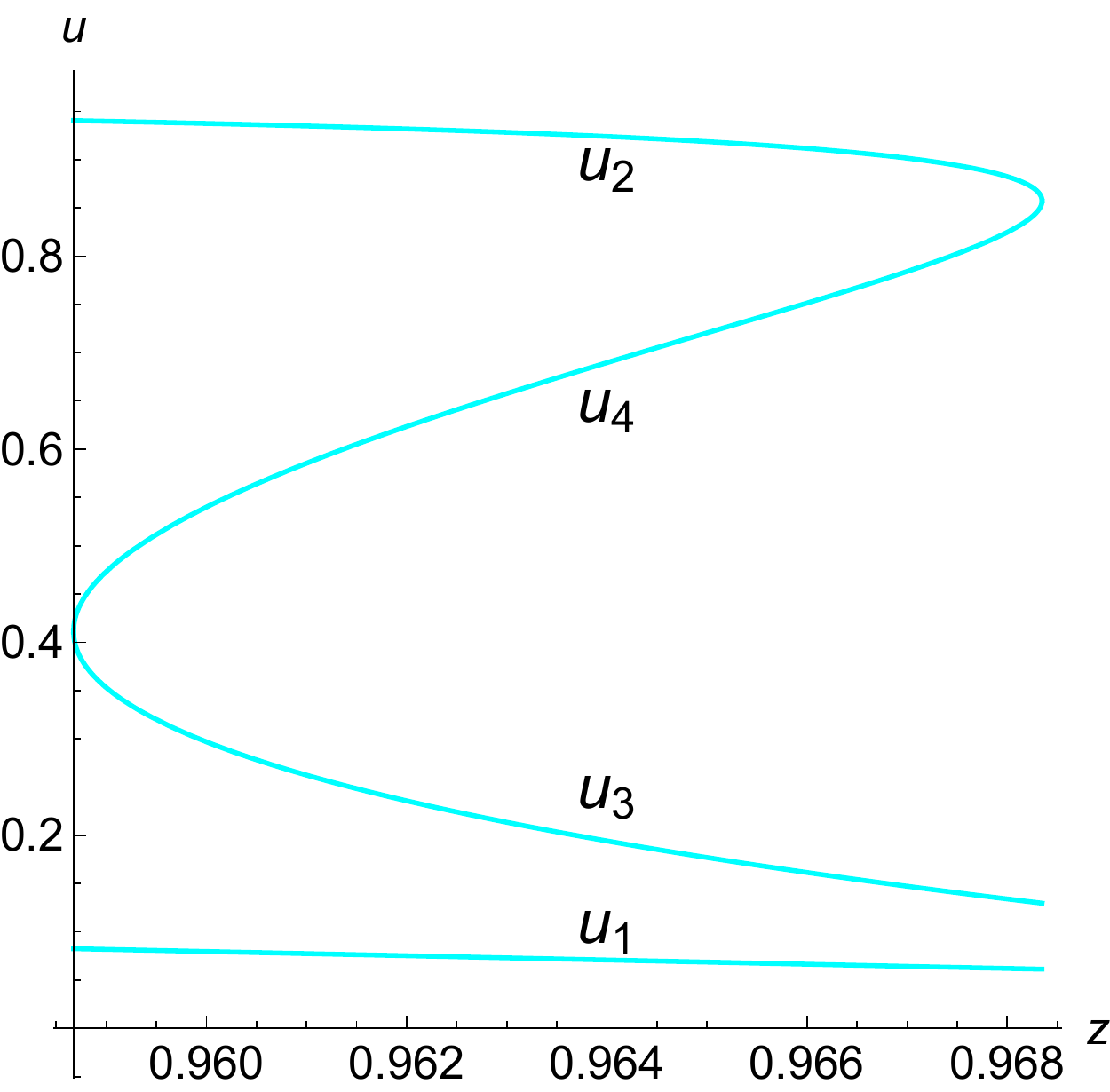}\\
(e): $(-1.25,-.15)$ {\bf (II)} and zoom, \qquad\qquad\qquad\qquad (f): $(-1,-.05)$ {\bf (III)} and zoom\\
 {\small \hspace*{-10mm}$z_{s_0}=0,\, z_{s_1}= 0.902896,\, z_{s_2}= 0.905097,\, z_{s_3}=0.905966$; \quad \qquad\quad $z_{s_1}=0.051008,\, z_{s_2}= 0.958672$,\qquad \qquad\qquad }
 {\small \hspace*{89mm} $ z_{s_3}=0.968347,\, z_{s_4}=0.974964$}
\\[3pt]
\includegraphics[width=8.7pc]{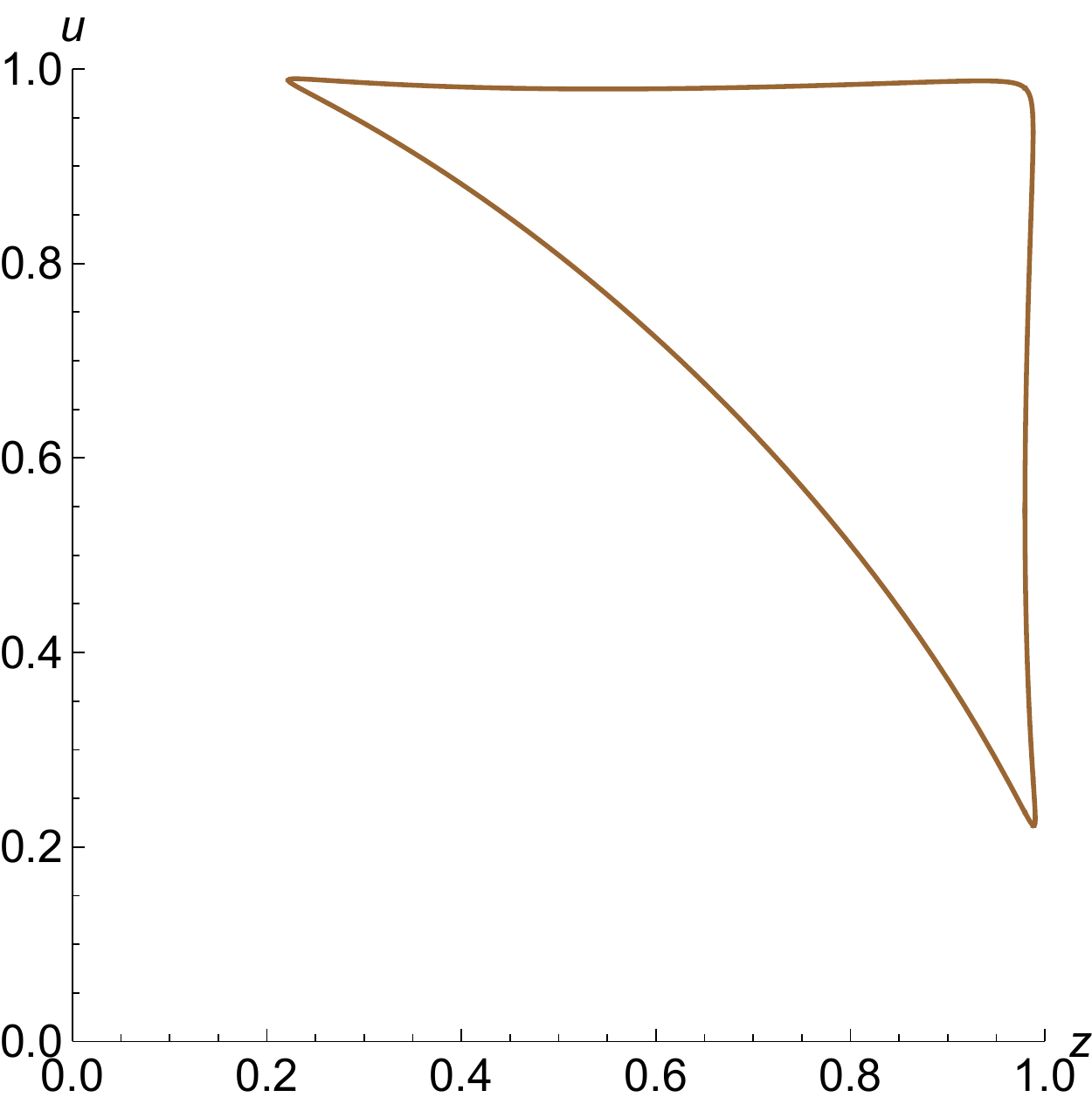}\quad \includegraphics[width=8.7pc]{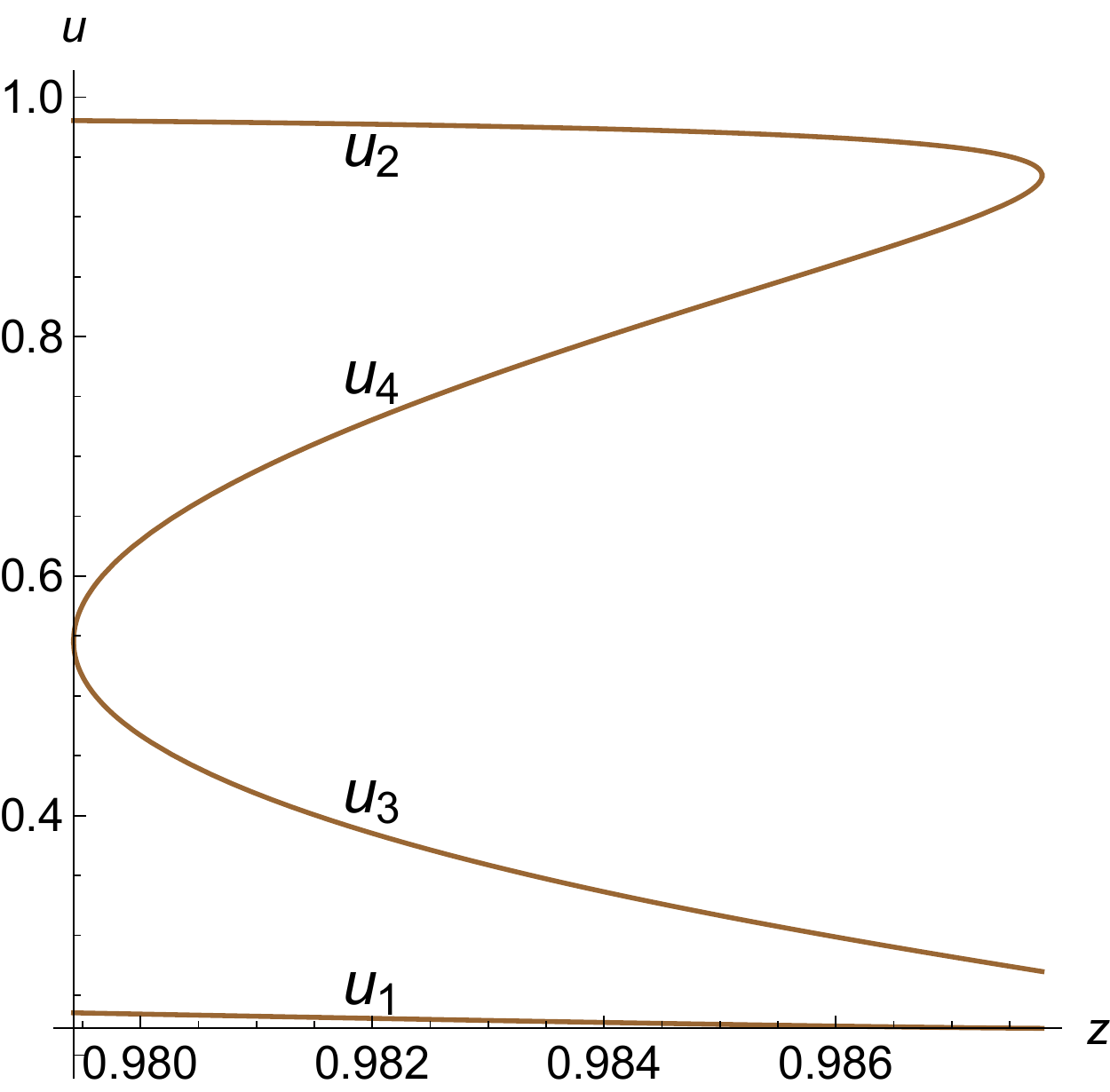}\quad
\includegraphics[width=8.7pc]{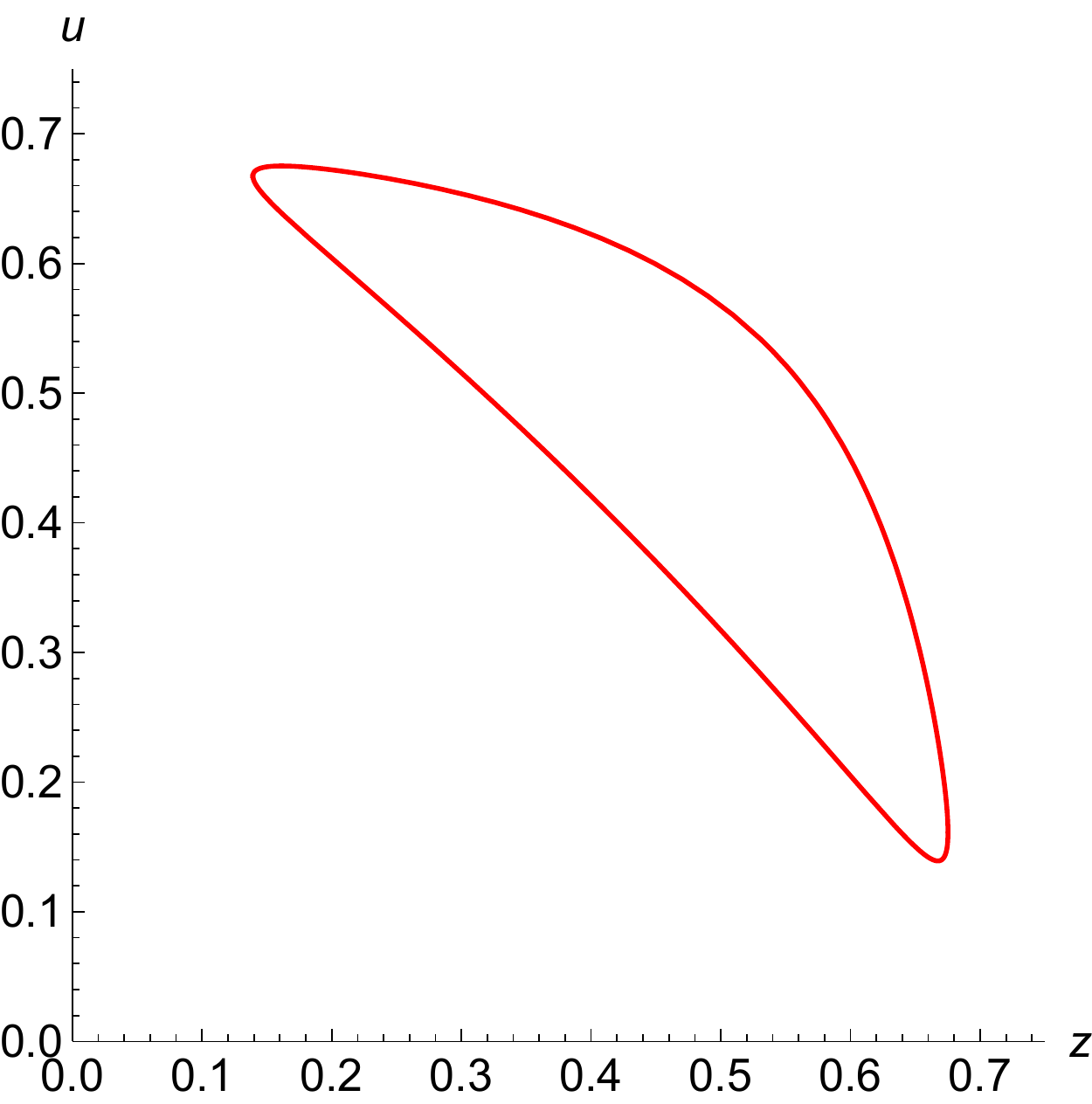}\quad \includegraphics[width=8.7pc]{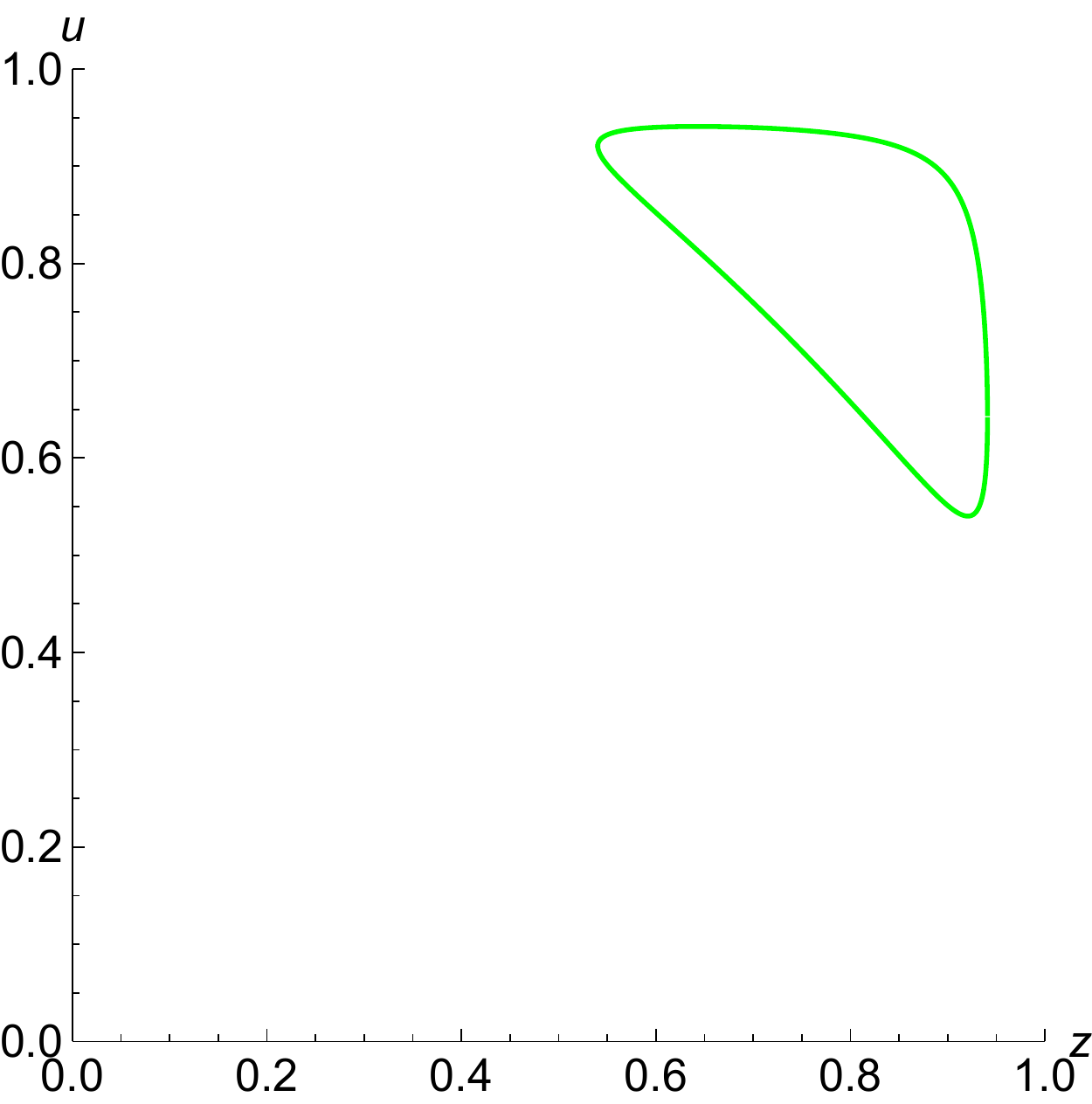}\\
(g): $(-0.8,-.02)$ {\bf (IV)} and zoom \qquad\qquad\quad (h): $(-1.5,-.6)$ {\bf (V)},\qquad (i): $(-.6, -.1)$ {\bf (VI)}\\
{\small \quad\quad $z_{s_2}= 0.22208, \,z_{s_3}= 0.97942$, \qquad \qquad \ $z_{s_1}=0.13922, \, z_{s_2}= 0.67523$; \ $z_{s_2}= 0.54021, \, z_{s_3}= 0.94092$,}\\
{\small \quad \quad $ z_{s_4}=0.98778, \, z_{s_5}=0.98997$ \hfill \null }
 \caption{A portrait gallery. At various points of the $(p,q)$ domain, 1st line: plot of the roots $u_i(z)$ of~$R$, which describe several branches of a closed curve; 2nd line: the value of $(p,q)$ and the region number; 3d line: the values of the {\it relevant} $z_s$.}\label{gallery}
\end{figure}

\subsection[The integrand $\varphi(z)$]{The integrand $\boldsymbol{\varphi(z)}$}
\begin{figure}[!tbp] \centering
 \includegraphics[width=36pc]{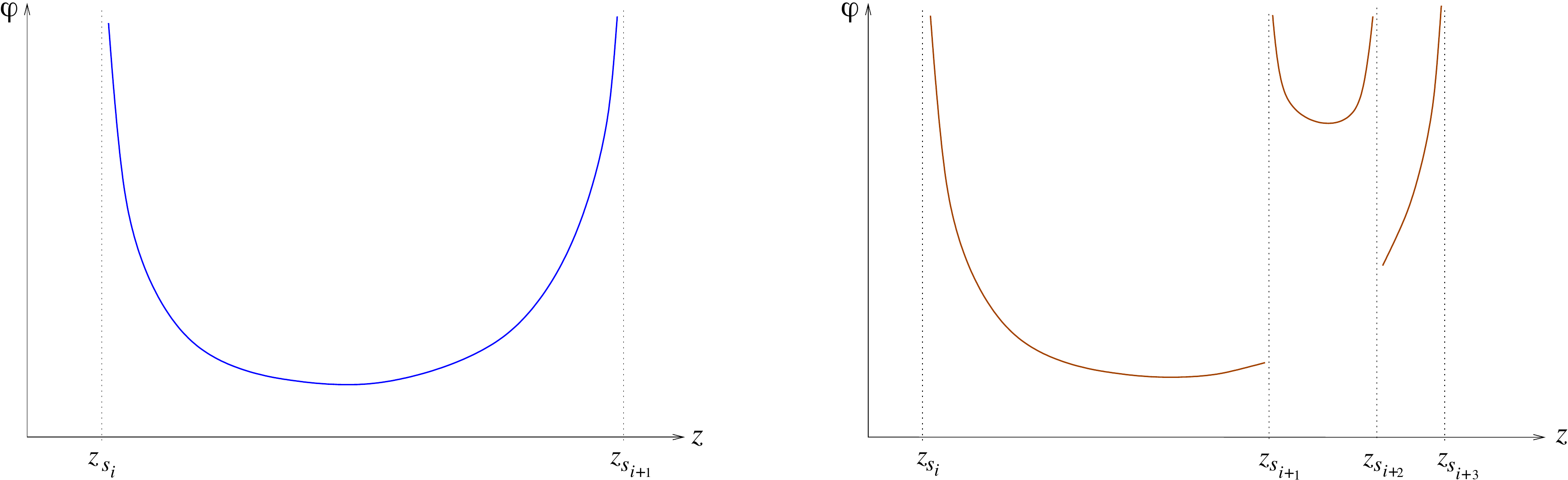}
\caption{Typical plots of $\varphi(z)$ in regions I, V or VI (left) or II, III, IV (right). In the latter, the middle and right intervals have been dilated for clarity. The discontinuity of $\varphi(z)$ at $z_{s_{i+1}}$, $z_{s_{i+2}}$ is due to the contribution of two new roots~$u_3$ and~$u_4$ in that interval.}\label{plots-of-phiI}
\end{figure}

Consider the integrand in~(\ref{pdf4})
\begin{gather}\label{phi} \varphi(z):= \sum_{{\rm roots}\ u_i \ {\rm of}\ R\atop |u_i|\le 1} \frac{(2+u+ z)}{|R'_u| }\Big|_{u=u_i}.\end{gather}
It has both integrable and non-integrable singularities, the latter where the integral diverges. Typical plots of $\varphi(z)$ in the various regions are displayed in Fig.~\ref{plots-of-phiI}.

 The singularities of $\varphi$ as a function of $z$ come either from singularities of $u_i(z)$, or from zeros of the denominator $|R'_u(u_i)|$. Both cases are associated with the merging of roots $u_i(z)$ of $R$, which occurs at some relevant root $z_s$ of its discriminant $\Delta_R$. The singularity of $u_i$ in the numerator is, however, at worst of square root type and gives rise to no divergence of the integral. We thus concentrate on the possible vanishing of the denominator $R'_u$.

 If we write the polynomial $R$ in a factorized form, $R = c \prod_{j} (u-u_j(z))$, its derivative at~$u_i$ reads $R'_u(u_i)= c\prod\limits_{j\ne i} (u_i(z)-u_j(z))$, and vanishes when~$u_i$ coalesces with some $u_j$, thus at some value $z=z_s$, a root of the discriminant $\Delta_R$.
\begin{enumerate}\itemsep=0pt
 \item[--] Either the pair of roots $(u_i,u_j)$ belonging to $[0,1]$ emerges from or disappears into the complex plane at $z_s$ with a square root behaviour, $u_{i,j}(z) \sim u_{i,j}(z_s) \pm \oh a |z-z_s|^{1/2} +\cdots$, hence $| u_i(z)-u_j(z)|\sim a |z-z_s|^{1/2}+\cdots$. Graphically, it manifests itself as a smooth curvature of the ``portraits'' of $u_i(z)$ and $u_j(z)$ at the points of vertical tangent in Fig.~\ref{gallery}. This is what happens at $z=z_{s_0}=0$ above the dashed line; and generically at the various relevant roots $z_s$ of $\Delta_R$. At such a point, the singularity of $\varphi(z)$ is integrable.
 \item[--] \looseness=-1 Or the two roots $u_i$, $u_j$ cross at a finite angle at $z_s$: $u_i(z) -u_j(z)\sim b(z-z_s) + o(z-z_s)$ with a finite coefficient $b$. This is what happens along the lines $p-q+1=0$, or $q=0$, and translates graphically into an angular point in Fig.~\ref{gallery}, see for example cases~(a),~(d),~(g). At such a point, the integral of $\varphi(z)$ diverges ``logarithmically'' at $z=z_s$, which supposes the introduction of some cut-off that measures the departure of $(p,q)$ from the singular point. This explains the growths of the PDF observed in Fig.~\ref{polygonJzJz} along the lines $q=0$, $p-q+1=0$ and (by symmetry $q\leftrightarrow -q$) $p+q+1=0$,
 as we discuss in the next subsection.
 \item[--] Or at exceptional points, the difference $u_i(z) -u_j(z)$ may vanish faster at $z_s$. This is what happens at the point $(-1,0)$, where it vanishes as $|z_s-z|^{3/2}$. Graphically, the two curves~$u_i(z)$ and $u_j(z)$ form a cusp, see for instance Fig.~\ref{gallery}{f}. In that case, the integral of~$\varphi(z)$ diverges as an inverse power of the cut-off, see below.
\end{enumerate}

But there is another source of divergence of $\rho$. The function $\varphi$ itself may diverge as $(p,q)$ approaches a~singular point. This is what happens at the three corners $(-4,0),\ (-3,-2),\ (0,0)$ of the domain, where we shall see that $R'_u(u_i)$ vanishes for all $z$ in the integration interval.

\subsection{The PDF, plots and divergences}\label{gen-sing}

\begin{figure}[!tbp]
\centering\includegraphics[width=16pc]{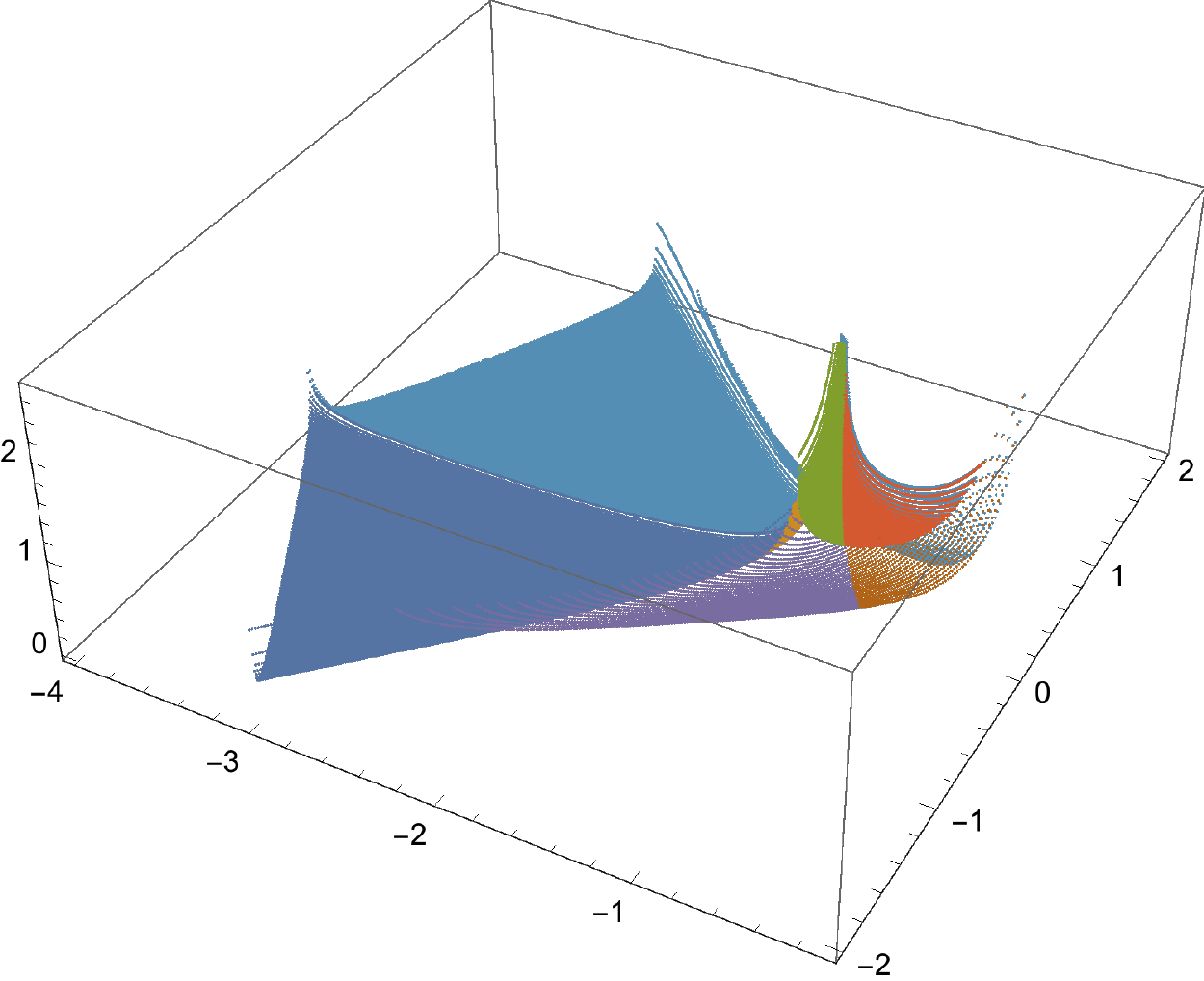}\qquad \includegraphics[width=16pc]{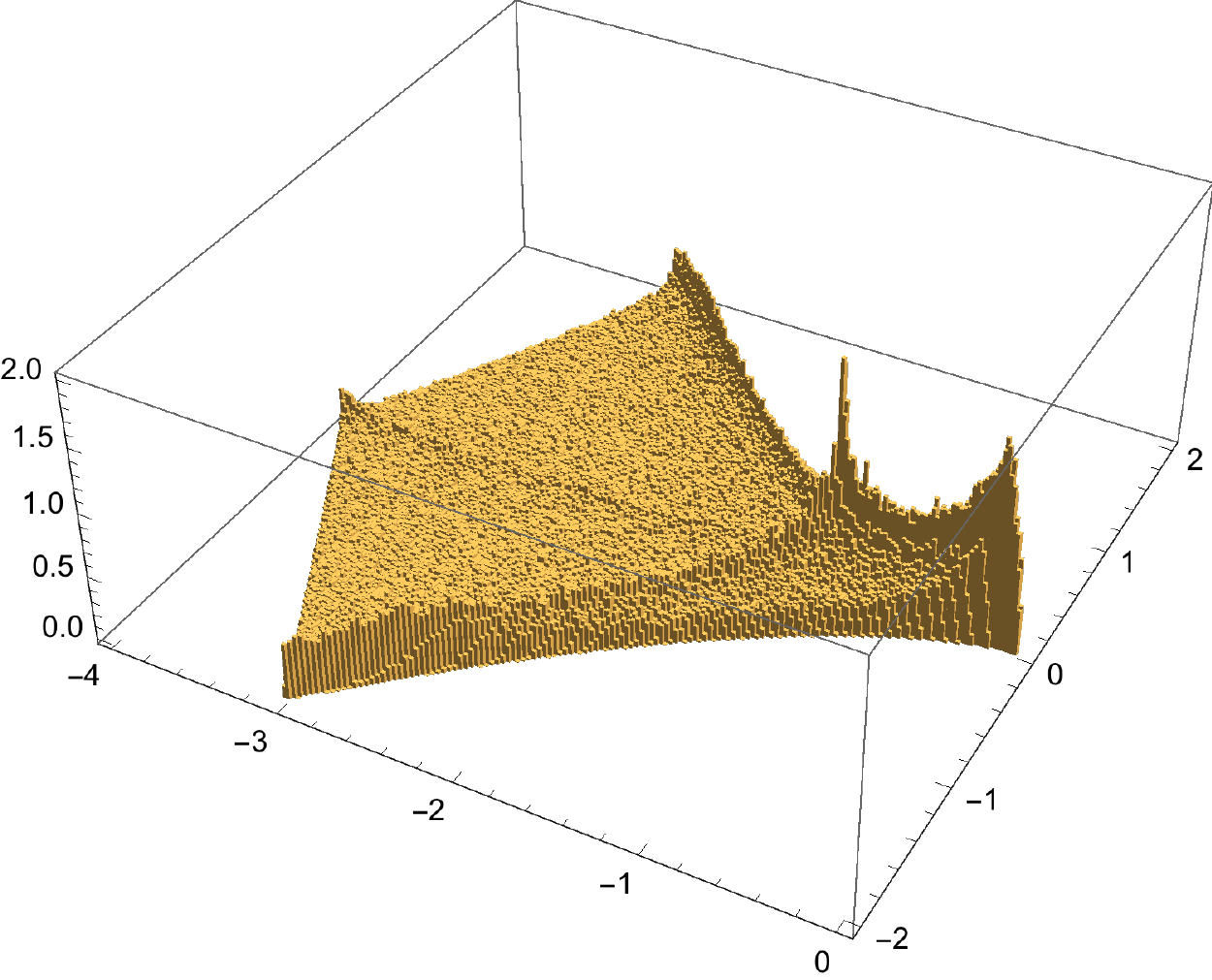}
\caption{Left: plot of the function $\rho$ computed on a grid of mesh $10^{-2}$; the regions I to VI, $q<0$ have been plotted in different colors (deep blue, orange, green, red, violet and brown), while their symmetric partners, for $q>0$, are all in light blue. Right: histogram of $10^6$ points in the $(p,q)$ plane.} \label{final-fig-pq}

\includegraphics[width=14pc]{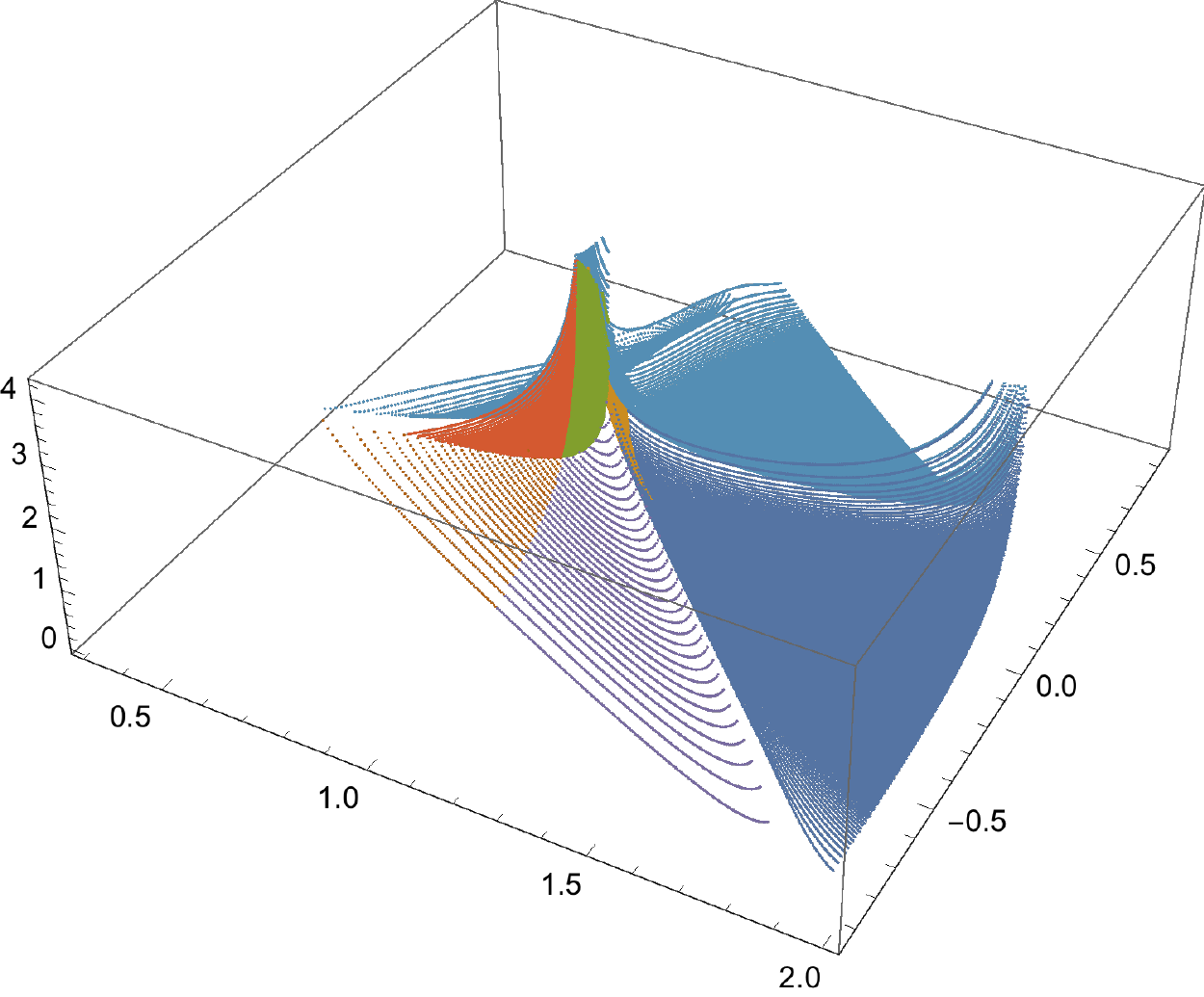}\qquad \includegraphics[width=14pc]{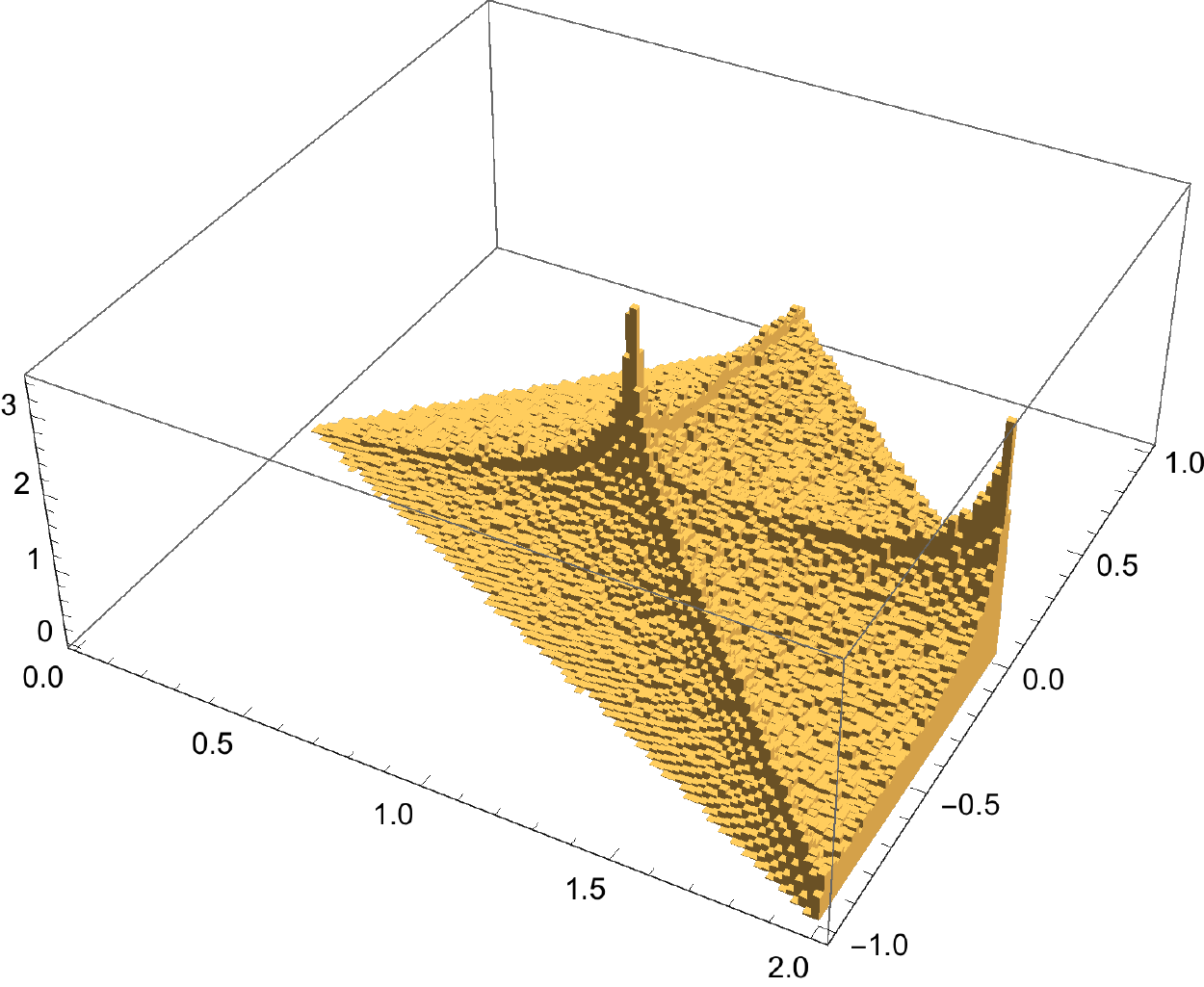}
\caption{The same, in the $(\gamma_1,\gamma_2)$ plane.}\label{final-fig-gamma}
 \end{figure}

We are now in position to draw the plot of the PDF $\rho(p,q)$ resulting from the integration in~(\ref{pdf4}) for $q\le 0$, supplemented by its mirror image by $q\to -q$, and to compare it with the histogram obtained by a simulation with $10^6$ points, see Figs.~\ref{final-fig-pq} and~\ref{final-fig-gamma}.

An important check consists in comparing the probability of occurrence of $(p,q)$ in a finite domain computed by integrating the PDF $\rho(p,q)$ over that domain to that estimated from a big sample of ``random events''. For example, $\Bbb{P}(-3.6\le p\le -3.5, q\le 0)\big|_{\rm computed} =0.04496$ while the estimate from a sample of size $10^6$ gives 0.044886.

As anticipated, the computed $\rho$ exhibits singularities along the lines $q=0$ and $p\pm q +1=0$. Note that the computation of $\rho$ is carried out point by point on a grid of mesh $10^{-2}$ in the $(p,q)$-plane, cutting off the vicinity of the singular lines, while the histogram uses bins of width $0.02$ throughout the domain. This explains the slight difference of appearance of the singularities.

By a long (and fairly tedious) case by case analysis, we may assert that the singularities are logarithmic in the approach of {\it generic} points of the singular lines $q=0$; $p-q+1=0$ for $q\le 0$; and $p+q+1=0$ for $q\ge 0$. At the end points and intersection of these lines, i.e., at the corners of the $(p,q)$-domain, as well as the point $(p,q)=(-1,0)$, the divergence is stronger, as an inverse power. This is summarized in the following Table, which gathers results obtained in the detailed discussion of the next subsection. The reader will also find in that subsection numerical verifications of the asserted divergences.

In most cases, we proceed as follows: as $(p,q)$ approaches a singular point, some $z_s$ approaches a limiting value $z_{s*}$ while the common value of a pair of coinciding roots, say, $u_s=u_1(z_s)=u_2(z_s)$ approaches $u_{s*} =u_{1,2}(z_{s*})$. Series expansions of $z_s$ and $u_s$ in powers of a ``distance'' $\epsilon\ll 1$ to the singularity may be computed. On the other hand, the roots $u_1(z)$ and $u_2(z)$ as well as the denominator $R'_u$ of $\varphi(z)$ in (\ref{pdf4}) may be expanded in powers of $\zeta=\sqrt{|z-z_s|} $. Finally, a~double series expansion in powers of $\epsilon$ and $\zeta$ is obtained for $R'_u(u_{1,2}(z))$, which upon integration in the vicinity of $z_s$, yields the singular contribution to $\rho$ of $z_{s*}$. This program is carried out in detail in the next subsection. Three particular singular points are treated separately.

\begin{table}[t]\centering \small \begin{tabular}{|c| c|c|c|c|c}
\hline
Position of singularity & Approach to singularity&{Divergent part} of $\rho$ \\
and $(z_{s*},u_{s*})$ && \\[3pt]
\hline && \\[-10pt]
$p-q+1 =0$, $q\le 0$ & & \\
$z_{s*}= 0$, $u_{s*}=(p+3)/2$ & $p-q+1\to 0$ &$
\rho_{\rm div}= \inv{ \pi^2 \sqrt{2} (p+3) \sqrt{-(1+p)} }\big|\log |p-q+1 | \big|$\\
$z_{s*}=(p+3)/2$, $u_{s*}=0$ &&\\[3pt]
\hline &&\\[-8pt]
$q= 0$, $-4< p<-1$ & $q\to 0$,\ $p$ fixed & $
\rho_{\rm div}= \inv{2 \pi^2 |p| \sqrt{p+4} }\big|\log |q | \big|$ \\
$z_{s*}= 1$, $u_{s*}=1$ &&\\[3pt]
\hline &&\\[-8pt]
$q= 0$, $-1< p<0 $ & & \\
$z_{s*}=p+1,\ u_{s*}=1$ &$q\to 0$,\ $p$ fixed &$
\rho_{\rm div}=\inv{2\pi^2 |p|}((1+p)^{-\oh} +(4+p)^{-\oh}) \big|\log |q | \big|$
 \\
$z_{s*}=1,\ u_{s*}=1$ & & \\
$z_{s*}=1,\ u_{s*}=p+1$ & & \\[3pt]

\hline &&\\[-8pt]
$(p,q)=(-1,0)$ & & \\
$z_{s*}=0,\ u_{s*}=1$&$q\to 0$, $p=-1$ &$
\rho_{\rm div}= \frac{1}{2\sqrt{6}\pi |q|^\oh} +O(\log|q|) $ \\
$z_{s*}=1,\ u_{s*}=0$ &&\\
$z_{s*}=u_{s*}=1$: subdominant &&
\\[3pt]
\hline &&\\[-8pt]
 $(p,q)=(-4,0)$ & &\\
$z_{s*}=0,\ u_{s*}=0$ & $\kappa (p+4) =-q \to 0$, &$
 \rho_{\rm div}=
C(\kappa)/ |q|^\oh$
\\
$z_{s*}=u_{s*}=\frac{2-\kappa }{3 \kappa +2}$ & $ 0<\kappa< 2$ &
\\[3pt]
\hline&&\\[-8pt]
$(p,q)=(0,0)$& $ \kappa p^3=-27q^2\to 0$ & $
 \rho_{\rm div}= {C''(\kappa)}/{ |q|^{2/3}} \sim 1/|p|$
 \\
$z_{s*}=u_{s*}=1$ & $0<\kappa<4$ & \\[3pt]
\hline&&\\[-8pt]
$(p,q)=(-3,-2)$ & $(p+3)=\kappa(q+2) \to 0$ & $
\rho_{\rm div}={C'(\kappa)}/{(q+2)^\oh}$ \\
$z_{s*}=u_{s*}=0$& $-\oh \le \kappa \le 1$ & \\[3pt]
\hline
\end{tabular}
\caption{Position and expression of the singularities of $\rho(p,q)$ in the $q<0$ part of the domain of Fig.~\ref{polygonJzJz}. The expressions of the coefficient functions $C(\kappa)$, etc., are given in {Appendix~\ref{detail-sing}}.} \label{Table2} \end{table}

Finally, note the PDF $\rho(p,q)$ does not vanish along the boundaries of the Horn domain in the $(p,q)$-plane. In particular, on the lower or upper sides of the domain, i.e., along the arcs of the cubic $q=\mp 2(-p/3)^{3/2}$, $\rho(p,q)$ has a finite limit. This is in no contradiction with the expected vanishing of the PDF~$\p(\gamma_1,\gamma_2)$ on the left boundaries $\gamma_1=\gamma_2$ and $\gamma_2=\gamma_3$ of the Horn domain in the $(\gamma_1,\gamma_2)$-plane, because of~(\ref{PDFgamma}) and of the vanishing of the Vandermonde determinant~$\Delta$ along those curves.

To make the discussion shorter, the detailed analysis of the singularities, leading to the conclusions of Table~\ref{Table2}, has been relegated to Appendix~\ref{detail-sing}.

\section{Zonal polynomials} \label{zonalpolynomials}

In the Hermitian case, it is known that the Horn problem discussed so far has a discrete counterpart, involving Littlewood--Richardson multiplicities, and may be regarded as a semi-classical limit of the latter. The PDF $\p(\gamma|\alpha,\beta)$, or rather the ``volume function'' ${\mathcal J}$ equal to the latter up to a Vandermonde factor, measures the distribution of (rescaled) Littlewood--Richardson multiplicities, i.e., structure constants of Schur polynomials, in the large scale limit, see \cite{CZ17, KT00}. In the real symmetric case, one expects similarly the PDF, or rather some ``volume like function'' ${\mathcal J :=\rho}$ proportional to it (see (\ref{PDFgamma})), to measure, at least in the generic case, the distribution of (rescaled) ``zonal multiplicities'', i.e., appropriate structure constants of zonal polynomials, see~\cite{FG2}. This motivates the discussion of the present section, where we pay special attention to the normalization and specializations of the Jack and zonal polynomials. We shall in particular use $\SU(n)$ reduction\footnote{I.e., eliminate integer partitions $\kappa$ when $\iota(\kappa) > n$, $\iota$ being the length of $\kappa$, and, when $\iota(\kappa) = n$, replace~$\kappa$ by $\kappa \operatorname{mod} \{1,1\ldots, 1\}$, the latter being therefore an {\sl extended partition} (see footnote~\ref{extendedpartition}) with a last part equal to~$0$.\label{sunreduction}} of zonal (or Jack) polynomials and introduce a notion of $\SU(n)$ zonal characters very similar to the usual Weyl characters.

\subsection{About Jack and zonal polynomials}

\subsubsection{Jack polynomials and their normalizations}

Zonal polynomials can be defined in many ways and we refer the reader to the abundant literature (see for instance \cite{Macdonald:book, zonalMathaiEtAl}). One possible approach is to start from Jack polynomials (themselves a particular case of the larger family called Macdonald polynomials). A~Jack polynomial with~$n$ variables is labelled by an integer partition~$\kappa$ and a real parameter $\alpha$. When one specializes the value of $\alpha$, Jack polynomials, in turn, give rise to various interesting families, in particular the Schur polynomials (case $\alpha =1$), the zonal polynomials (case $\alpha = 2$), and the quaternionic polynomials (case $\alpha = 1/2$).

Actually there are three variants of the Jack polynomials, denoted $J_X^\alpha$ with $X = P, C, J$\footnote{The capital subscripts $J$ and $C$, for both Jack and zonal polynomials, refer to the original papers of James \cite{JamesZonalPoly1,JamesZonalPoly3}, and Constantine~\cite{ConstantineZonalPoly}, see also~\cite{HualLoKen}.} differing by an overall normalization (an overall $\alpha$-dependent and partition-dependent numerical factor): for example, one writes $J_J=c_{PJ} J_P$, etc. When $\alpha = 2$ one has therefore also three kinds of zonal polynomials respectively denoted $Z_P$, $Z_J$ and $Z_C$. When studying zonal or Jack polynomials, many authors -- in particular in old papers -- use the James normalization (polynomials $Z_J$) without saying so explicitly. As we shall see, the most interesting family, for us, is the family of the zonal polynomials (and also the Jack polynomials) defined with the P normalization, the reason, that will be discussed below, is that this normalization is compatible with $\SU(n)$-reduction and with the conjugation of irreducible representations (irreps) of $\SU(n)$~-- the latter being described by integer partitions with at most $n-1$ parts. When expanded in terms of monomial symmetric polynomials, the $J$ normalization of the Jack polynomial defined by the partition $\kappa$ makes the coefficient of the lowest order monomial $[1^n]$ equal to $n!$, whereas, using the same expansion, the~$P$ normalization makes the coefficient of the monomial relative to the highest partition (i.e., $\kappa$) equal to~1. For a given partition $\kappa$, the normalization factor~$c_{PJ}$ is the lower $\alpha$-hook coefficient of~$\kappa$.

\emph{Note:} Zonal polynomials, with an un-specified normalization, were denoted $Z(\kappa)(x)$ in~(\ref{CI1}).

\subsubsection{Packages}
 Zonal polynomials and, more generally, Jack polynomials (variables are called $x_j$), are usually written in terms of monomial symmetric functions, or in terms of power sums, not very often in terms of the variables $x_j$ themselves because this would take too much space. To the authors' knowledge there are very few computer algebra packages devoted to the manipulation of those polynomials; we should certainly mention \cite{MOPS}, written for Mapple (that we did not use), and the small package \cite{Baratta} written for Mathematica~-- but it is slow, unstable (division by~0), and uses an obsolete version of the language. For those reasons we developed our own, using Mathematica: the definition chosen for Jack polynomials uses a recurrence algorithm in terms of skew Young diagrams and a modified Pieri's formula described by Macdonald in \cite{Macdonald:book}, see also \cite{DemmelKoev}; our code, which also contains commands to convert Jack, Schur and zonal polynomials to several other basis (elementary symmetric polynomials, power sums, monomial sums, complete sums), and commands for calculating structure constants in each basis, is freely available on the web site~\cite{MathematicaProgramsRC}; the same package contains commands giving the coefficients
$Z(\kappa)(I)$ and $c(\kappa)$ that appear in formula (\ref{CI1}), with various normalization choices.

\subsubsection{Structure constants} \label{zonal:structure constants}

Zonal polynomials form a basis of the space of the ring of symmetric polynomials in $n$ variables. Structure constants in this basis are only {\it rationals} (by way of contrast, the coefficients of Schur polynomials in the expansion of a product of two Schur polynomials are non-negative integers). For illustration, let us consider the zonal polynomial(s) for the extended\footnote{Given an integer partition $\kappa$, say of length $s$, of the integer $m$, it is convenient, in order to specify the number of variables in symmetric polynomials, to call ``extended partition'' of length $n$, assuming that $n\geq s$, the partition that is obtained from $\kappa$ by padding $n-s$ zeros to the right of $\kappa$. The length of the obtained partition (no longer an integer partition in the strict sense) is then equal to $n$, the chosen number of variables. Example: $Z_P(\{2,1,0,0\}) = Z_P(\{2,1\})(x_1,x_2,x_3,x_4)$.\label{extendedpartition}} partition $\{2,1,0\}$ (i.e., three variables $x_1$, $x_2$, $x_3$), and the decomposition of its square, using the three standard normalizations; we also give the decomposition obtained for the {\sl square} of the Schur polynomial $s(\{2,1,0\})$:
\begin{gather*}
Z_P(\{2,1,0\})^2= \frac{25}{12} {Z_P}(\{2,2,2\})+\frac{12}{5} {Z_P}(\{3,2,1\})\\
\hphantom{Z_P(\{2,1,0\})^2=}{} +\frac{4}{3} {Z_P}(\{3,3,0\})+\frac{4}{3} {Z_P}(\{4,1,1\})+{Z_P}(\{4,2,0\}),\\
Z_J(\{2,1,0\})^2= \frac{5}{54} Z_J(\{2,2,2\})+\frac{12}{35} Z_J(\{3,2,1\})+\frac{4}{135} Z_J(\{3,3,0\})\\
\hphantom{Z_J(\{2,1,0\})^2=}{} +\frac{32}{405} Z_J(\{4,1,1\})+\frac{1}{27} Z_J(\{4,2,0\}),\\
Z_C(\{2,1,0\})^2= \frac{3}{4} {Z_C}(\{2,2,2\})+\frac{12}{25} {Z_C}(\{3,2,1\})+\frac{21}{25} {Z_C}(\{3,3,0\})\\
\hphantom{Z_C(\{2,1,0\})^2=}{} +\frac{12}{25} {Z_C}(\{4,1,1\})+\frac{63}{125} {Z_C}(\{4,2,0\}), \\
 s(\{2,1,0\})^2= s(\{2,2,2\})+2 s(\{3,2,1\})+s(\{3,3,0\})+s(\{4,1,1\})+s(\{4,2,0\}).
\end{gather*}
Specifying the number of variables matters. In the previous example for instance, the square of $Z_P(\{2,1,0,0\})$, i.e., four variables, reads:
\begin{gather*}
\frac{125}{63} Z_ P(\{2,2,1,1\})+\frac{25}{12} Z_ P(\{2,2,2,0\})+\frac{25}{18} Z_ P(\{3,1,1,1\})+\frac{12}{5} Z_ P(\{3,2,1,0\})\\
\qquad{} +\frac{4}{3} Z_ P(\{3,3,0,0\})+\frac{4}{3} Z_ P(\{4,1,1,0\})+ Z_ P(\{4,2,0,0\}).
\end{gather*} Of course $Z_ P(\{\kappa_1,\kappa_2,\kappa_3,0\})$ restricts to $Z_ P(\{\kappa_1,\kappa_2,\kappa_3\})$ if $x_4=0$, but the last decomposition also contains new terms like $Z_ P(\{2,2,1,1\})$ {or $Z_ P(\{3,1,1,1\})$} that vanish when $x_4=0$.

\subsubsection[A particular feature of structure constants in the $Z_P$ basis]{A particular feature of structure constants in the $\boldsymbol{Z_P}$ basis}

We pause here to notice that the coefficients of ${Z_P}(\{3,3,0\})$ and of ${Z_P}(\{4,1,1\})$ are the same (both equal to $4/3$), with the same remark for the coefficients of ${s}(\{3,3,0\})$ and of ${s}(\{4,1,1\})$, which are both equal to $1$. This is not so for ${Z_J}$ and ${Z_C}$.

\looseness=1 In the Schur case this remark is not surprising: indeed, in terms of irreps of \SU(3), the Schur decomposition of $s(\{2,1,0\})^2$ corresponds to the tensor decomposition of the square of the adjoint representation labelled\footnote{The components of a highest weight $\lambda=[\lambda_1, \lambda_2]$ are written in the Dynkin basis (the basis of fundamental weights).} by its highest weight $[1,1]$; in other words one recovers the well-known decomposition $[1,1]^{\otimes 2} = [0,0] \oplus 2 [1,1] \oplus [0,3] \oplus [3,0] \oplus [2,2]$. The two integer partitions $\{3,3\}$ and $\{4,1,1\}$ determine, {\sl once they are reduced to \SU(3)} (see footnote~\ref{sunreduction}), two conjugate Young diagrams of respective shapes $\{3,3\}$ and $\{3\}$ describing the two complex conjugate irreducible representations $[0,3]$ and $[3,0]$ of \SU(3). The irrep $[1,1]$, associated with the partition $\{2,1,0\}$, is self-conjugate, actually real, and the coefficients of those complex conjugate irreps that appear in its square, in particular $[0,3]$ and $[3,0]$, are, of course, equal.

More generally we observe (the proof of this conjecture is left to the reader) that the following property holds in the zonal $P$ case\footnote{We remind the reader that Schur polynomials can be obtained from Jack polynomials, with the P normalization, just by setting $\alpha = 1$ (no pre-factors).}: Consider $\lambda$ and $\mu$, two irreducible representations of~$\SU(n)$ described by extended partitions $\ell(\lambda)$ and $\ell(\mu)$ of length $n$ (the last component \mbox{being~0}). Decompose the product $Z_P(\ell(\lambda)) \cdot Z_P(\ell(\mu))$ on zonal polynomials using the normalization P and call $f_{\ell(\lambda), \ell(\mu), \chi}$ the structure constants:
\begin{gather*}
Z_P(\ell(\lambda)) \cdot Z_P(\ell(\mu)) = \sum_{\chi^\prime} f_{\ell(\lambda), \ell(\mu), \chi^\prime} Z_P(\chi^\prime).\end{gather*} Calling $\bar\lambda$ and $\bar\mu$ the complex conjugate representations, we decompose in the same way the product
\begin{gather*} Z_P(\ell(\bar\lambda)) \cdot Z_P(\ell(\bar\mu)) = \sum_\chi^{\prime \prime} f_{\ell(\bar\lambda), \ell(\bar\mu), \chi^{\prime \prime}} Z_P(\chi^{\prime \prime}).
\end{gather*} Call $\chi|_{\SU(n)}$ the restriction of a partition $\chi$ to~$\SU(n)$~-- see footnote~\ref{sunreduction}. Then, if $\chi^{\prime \prime}|_{\SU(n)} = \bar\chi^\prime |_{ \SU(n)}$, we observe that we have $f_{\ell(\lambda), \ell(\mu), \chi^\prime} = f_{\ell(\bar\lambda), \ell(\bar\mu), \chi^{\prime \prime}}$. The case where we take $\mu = \bar\lambda$ and the sub-case where we also assume $\lambda = \bar\lambda$, are of particular interest. The latter gives:

If $\lambda$ is a self-conjugate irreducible representation of $\SU(n)$ described by a Young diagram of shape $\chi$ (an integer partition, or an extended partition of length~$n$), and if $\chi^\prime$ and $\chi^{\prime\prime}$ are two partitions appearing in the decomposition of the square of the zonal polynomial $Z_P(\chi)$ in the $Z_P$ basis that give rise, after $\SU(n)$ reduction, to complex conjugate representations, the coefficients (structure constants) of $Z_P(\chi^\prime)$ and of $Z_P(\chi^{\prime\prime})$ are equal.

Actually, the same property seems to hold for all values of the Jack parameter $\alpha$, when using the $P$ normalization. It does not hold for the normalizations $J$ and $C$.

Remember that the notion of complex conjugation on $\SU(n)$ irreps can be described in purely combinatorial terms: if $\kappa$ is the integer partition describing some irrep of $\SU(n)$, its length (number of parts) obeys $\iota(\kappa) < n$; then one obtains the partition describing the complex conjugate representation by taking the complement of (the Young diagram of) $\kappa$ in a rectangle which is $\kappa(1)$ units wide ($\kappa(1)$ being the largest part of $\kappa$) and $(\iota(\kappa) +1)$ units deep.

\subsection[From zonal polynomials to $\SU(n)$ zonal characters]{From zonal polynomials to $\boldsymbol{\SU(n)}$ zonal characters}

Although we have in mind applications to the zonal case $\alpha = 2$, or to the quaternionic (zonal) case $\alpha = 1/2$, most of our considerations, in the section that follows, apply to arbitrary values~$\alpha$ of the Jack parameter.

\subsubsection[$\SU(n)$-zonal characters]{$\boldsymbol{\SU(n)}$-zonal characters}

An irrep of $\SU(n)$, is characterized by its highest weight (hw)~$\lambda$. Its components in the basis of fundamental weights (Dynkin labels) are denoted $[\lambda_1, \ldots, \lambda_{n-1}]$. When considering irreps of~$\U(n)$ one adds a last index $\lambda_n$; here we only consider the case~$\SU(n)$, but it is often handy to keep this last index, while setting $\lambda_n=0$. One can also characterize the same irrep~$\lambda$ by the Young diagram defined by the extended partition $\alpha=\ell(\lambda)$ with components $\ell_i(\lambda)=\sum\limits_{j=i}^{n} \lambda_j$, $i=1,\ldots,n$, obeying the constraint $\alpha_i \geq \alpha_{i+1}$ for all~$i$. Conversely, given a partition~$\delta$ of length~$n$, extended or not, one obtains a~highest weight~$\lambda$ for $\SU(n)$ by setting $\lambda_i = \delta_i - \delta_{i+1}$; such a partition $\delta$ differs from $\alpha=\ell(\lambda)$ by a constant shift. We denote by $N_{\lambda\mu}^\nu$ the multiplicity\footnote{$N_{\lambda\mu}^\nu$ is sometimes called the Littlewood--Richardson (LR) multiplicity, although, strictly speaking, the latter refers to the coefficient of $\ell(\nu)$ in the decomposition in the Schur basis of the product of two Schur polynomials respectively defined by the partitions $\ell(\lambda)$ and $\ell(\mu)$. This decomposition often contains terms labelled by integer partitions of length larger than $n$, therefore not contributing to the tensor product of $\SU(n)$ representations~-- a~Young diagram of $\SU(n)$ cannot have more than $n-1$ lines.} of the irrep $\nu$ in the tensor product of the irreps of~$\SU(n)$ defined by $\lambda$ and $\mu$.

Given a dominant weight $\lambda$ of $\SU(n)$, i.e., a non-negative integer combination of the fundamental weights, call $\ell(\lambda)$ its associated partition of length $n$ (i.e., $\ell(\lambda)_n=0$) and take
the Jack-$P$ polynomial $J_P^\alpha (\ell(\lambda)) (x_1,\ldots,x_n)$ determined by the partition $\ell(\lambda)$. We then
consider the following Laurent polynomial in the variables $y_1,\ldots, y_{n-1}$:
\begin{gather*} J_P^\alpha (\ell(\lambda)) \left(x_1=y_1, x_2=\frac{y_2}{y_1}, \ldots, x_j = \frac{y_j}{y_{j-1}}, \ldots, x_{n-1} =\frac{y_{n-1}}{y_{n-2}}, x_n=\frac{1}{y_{n-1}}\right).\end{gather*}
If $\alpha=1$, $J_P^\alpha$ is a Schur polynomial and the previous Laurent polynomial is recognized as the Weyl character $\chi(\lambda)$ of the irrep $\lambda$, for the Lie group $\SU(n)$.\\
If $\alpha=2$, i.e., when $J_P^\alpha$ is a zonal polynomial $Z_P$ (with the normalization $P$), we introduce, by analogy, and for lack of a better name, the following notation and definition:\footnote{Another notion of ``zonal character'' can be found in the literature~\cite{FeraySniady}, but it is related to the symmetric group, not to irreducible representations of $\SU(n)$. It differs from the notion that we consider here.}
\begin{definition}
The zonal character $\chi_Z(\lambda)$ of $\SU(n)$ associated with the dominant weight $\lambda$ is defined as the Laurent polynomial
\begin{gather*} \chi_Z(\lambda) (y_1,\ldots, y_{n-1}) = Z_P(\ell(\lambda))\left(y_1, \frac{y_2}{y_1}, \ldots, \frac{y_j}{y_{j-1}}, \ldots, \frac{y_{n-1}}{y_{n-2}}, \frac{1}{y_{n-1}}\right).\end{gather*}
\end{definition}

Now, moving to the Lie algebra $\su(n)$ -- or, equivalently, to trigonometric characters, we start from the same $Z_P$ polynomial expressed in terms of $x_j$ variables but this time perform the following {consecutive} transformations on its arguments: $x_{i}\to {\rm e}^{i \big(a_i - \frac{1}{n}\sum\limits_{{j=1}}^{{n}} a _j\big)}$, then $a _{{j}}\to a _1-\sum\limits_{i=1}^{j-1} u_i$. The result is a trigonometrical expression in the variables $u_j$, that we call the {\sl Lie algebra zonal character} of $\su(n)$ associated with the hw $\lambda$, or the {\sl trigonometric zonal character} of $\SU(n)$ associated with the hw $\lambda$.

Let us give one example. Take $n=3$ and $\lambda = [1,1]$. The associated (extended) partition is $\ell(\lambda) =\{2,1,0\}$. The Jack polynomial $J_P^\alpha(\{2,1,0\})$ in terms of the variables $x_1$, $x_2$, $x_3$, the associated Laurent polynomial, and its trigonometric version are given below. The Schur polynomial $s(\{2,1,0\})$ and the Zonal-P polynomial $Z_P(\{2,1,0\})$ are obtained from the first expression by setting respectively $\alpha = 1$ and $\alpha = 2$. The corresponding $\SU({3})$ and $\su({3})$ zonal characters are obtained from the last two expressions by setting $\alpha = 2$
\begin{gather*}
J_P^\alpha(\{2,1,0\}) = x_1^2 x_2+x_1 x_2^2+x_2^2 x_3+x_2 x_3^2+x_1^2 x_3+x_1 x_3^2+\frac{6 x_1 x_2 x_3}{2+\alpha }, \\
\chi^{\SU({3}) }_{\alpha}([1,1]) = \frac{6}{\alpha +2}+\frac{y_1^2}{y_2}+\frac{y_2^2}{y_1}+ \frac{ y_1}{y_2^2} + \frac{y_2}{y_1^2}+y_1 y_2 + \frac{1}{y_1 y_2},\\
\chi^{\su({3}) }_{\alpha}([1,1]) = 2 \left(\frac{3}{\alpha +2}+\cos u_1 +\cos u_2+\cos (u_1+u_2 )\right).
\end{gather*}
Taking $\alpha=1$ in the second expression, one recognizes the Weyl character of the adjoint representation of~\SU(3), the powers (positive or negative) of the $y_j$ being, as usual, the components of the weights of the weight system of this representaton in the basis of fundamental weights. What happens in the zonal case, and more generally when $\alpha \neq 1${,} is that the ``multiplicities'' of the weights are no longer integers.
For this particular irrep, only the multiplicity of the weight at the origin of the weight system is modified by $\alpha$, its value being $2$ in the usual (Schur) case but $3/2$ in the zonal case.
Notice that the trigonometric expression is real~-- it is so because the hw $[1,1]$ is self-conjugate, otherwise the obtained expression would be complex. The arguments being specified (partitions or Dynkin labels), we shall denote $\chi^{\SU(n) }_{\alpha}$, in the cases $\alpha = 2$ and $\alpha = 1/2$, by $\chi_Z$ and~$\chi_Q$.

\subsubsection[Structure constants for ${\rm SU}(n)$-zonal characters]{Structure constants for $\boldsymbol{{\rm SU}(n)}$-zonal characters}

Given two $\SU(n)$ irreps, there are many ways to obtain the decomposition of their tensor products into a sum of irreps. The honeycomb technique, for instance,
is very fast\footnote{And the semi-magic square algorithm, valid for $\SU(3)$, is even faster, see~\cite{CZNuclPhys}.} but it is not available in the zonal case that we consider.
However, we can replace the multiplication of the associated $\SU(n)$ Weyl characters by the multiplication of the associated $\SU(n)$ zonal characters as defined above, which amounts to {using} the structure constants for the appropriate product of zonal polynomials. Let us illustrate this with our favorite example, the square of $[1,1]$. From the already given decomposition of the square of $Z_P(\{2,1,0\})$ we obtain immediately:
\begin{gather*}
\chi_Z([1,1])^2 = \frac{25}{12} \chi_Z([0,0])+\frac{12}{5} \chi_Z([1,1])+\frac{4}{3} \chi_Z([0,3])+\frac{4}{3} \chi_Z([3,0])+ \chi_Z([2,2]).
 \end{gather*}
The same decomposition can be obtained by using Laurent polynomials since the associated \SU(3) zonal characters are as follows (the reader can then check that the previous equality holds):
 \begin{gather*}
 \{2,1,0\} \!\rightarrow\! \chi_Z([1,1])= \frac{y_1^2}{y_2}+\left(y_2+\frac{1}{y_2^2}\right)
 y_1+\frac{y_2^2+\frac{1}{y_2}}{y_1}+\frac{y_2}{y_1^2}+\frac{3}{2} \\
 \hphantom{ \{2,1,0\} \!\rightarrow\! \chi_Z([1,1])}{} = \frac{2 y_2 y_1^4+2 y_2^3 y_1^3+2 y_1^3+3 y_2^2 y_1^2+2 y_2^4 y_1+2 y_2 y_1+2 y_2^3}{2
 y_1^2 y_2^2}, \\
 \{2,2,2\} \!\rightarrow\! \chi_Z([0,0])= 1, \\
 \{3,2,1\} \!\rightarrow\! \chi_Z([1,1])= \text{already given}, \\
 \{3,3,0\} \!\rightarrow\! \chi_Z([0,3])= \frac{1}{5} \!\left(\!\frac{5 y_1^3}{y_2^3}+\frac{3 y_1^2}{y_2}+3 y_2 y_1+\frac{3
 y_1}{y_2^2}+5 y_2^3+\frac{3 y_2^2}{y_1}+\frac{3}{y_2 y_1}+\frac{3
 y_2}{y_1^2}+\frac{5}{y_1^3}+2\!\right)\!, \\
 \{4,1,1\} \!\rightarrow\! \chi_Z([3,0])= \frac{1}{5} \!\left(\!5 y_1^3+\frac{3 y_1^2}{y_2}+3 y_2 y_1+\frac{3
 y_1}{y_2^2}+\frac{5}{y_2^3}+\frac{3 y_2^2}{y_1}+\frac{3}{y_2 y_1}+\frac{3
 y_2}{y_1^2}+\frac{5 y_2^3}{y_1^3}+2\!\right)\!, \\
 \{4,2,0\} \!\rightarrow\! \chi_Z([2,2])= \frac{1}{6} \left(\frac{6 y_1^4}{y_2^2}+\frac{4 y_1^3}{y_2^3}+4 y_1^3+6 y_2^2
 y_1^2+\frac{6 y_1^2}{y_2}+\frac{6 y_1^2}{y_2^4}+6 y_2 y_1+\frac{6 y_1}{y_2^2}+4
 y_2^3\right.\\
 \left. \hphantom{\{4,2,0\} \!\rightarrow\! \chi_Z([2,2])=}{} +\frac{4}{y_2^3}+\frac{6 y_2^2}{y_1}+\frac{6}{y_2 y_1}+\frac{6
 y_2^4}{y_1^2}+\frac{6 y_2}{y_1^2}+\frac{6}{y_2^2 y_1^2}+\frac{4 y_2^3}{y_1^3}+\frac{4}{y_1^3}+\frac{6 y_2^2}{y_1^4}+9\!\right)\!.
 \end{gather*}

The previous result -- and more generally any decomposition of a product of such characters -- can be checked by using a concept of dimension.
The dimension of an irreducible representation is the value taken by the associated Weyl $\SU(n)$ character at $y_j = 1$ (or the value taken by the $\su(n)$ character at $u_j = 0$), for all $j$.
In the same way one can define a ``zonal dimension'' for an irrep of hw $\lambda$ as the value taken by the $\SU(n)$ zonal character $\chi(\lambda)$ for $y_j=1$ (or the value taken by the $\su(n)$ character for $u_j=0$).
This is a slight terminological abuse since the obtained number is not an integer in general, but this dimension function is obviously compatible with addition and multiplication, as it should.
In the case\footnote{In the case of \SU(2) we find that the zonal dimension is given by
$\dim_Z^{}([\lambda]) = \frac{\sqrt{\pi } \Gamma (\lambda+1)}{\Gamma \left(\lambda+\frac{1}{2}\right)} = \frac{(2 \lambda)\text{!!}}{(2 \lambda-1)\text{!!}}$,
a formula which is the zonal analog of $\dim([\lambda]) = \frac{\Gamma(\lambda+2)}{\Gamma(\lambda+1)}=\lambda+1$.}
of \SU(3), we claim that
\begin{gather*} \dim_Z^{}([\lambda_1,\lambda_2])=\frac{{\lambda_1}! {\lambda_2}! (2 {\lambda_1}+2 {\lambda_2}+1){!!}}{(2 {\lambda_1}-1){!!} (2 {\lambda_2}-1){!!} ({\lambda_1}+{\lambda_2})!}.
\end{gather*}
The work relating the above to the quantity $u_0(P_\lambda)$, or $s_\lambda$ in the Schur case, defined by Macdonald in \cite[Section~6.11]{Macdonald:book} is left to the reader.
 Written $ \frac{2 \sqrt{\pi } \Gamma (\lambda_1+1) \Gamma (\lambda_2+1) \Gamma \left(\lambda_1+\lambda_2+\frac{3}{2}\right)}{\Gamma
 \left(\lambda_1+\frac{1}{2}\right) \Gamma \left(\lambda_2+\frac{1}{2}\right) \Gamma (\lambda_1+\lambda_2+1)}$ this expression looks very similar to the standard \SU(3) dimension
 \begin{gather*} \dim([\lambda_1,\lambda_2])=\frac{ (1+\lambda_1) (1+\lambda_2) (2+\lambda_1+\lambda_2)}{2} = \frac{\Gamma (\lambda_1+2) \Gamma (\lambda_2+2) \Gamma (\lambda_1+\lambda_2+3)}{2 \Gamma (\lambda_1+1) \Gamma (\lambda_2+1) \Gamma (\lambda_1+\lambda_2+2)}.
 \end{gather*}
This zonal dimension is (non-surprisingly) related to the normalization coefficient $Z_J(p)(I)$ that enters~(\ref{CI1}), and was introduced by A.T.~James in~\cite{JamesZonalPoly2}:
\begin{gather*} \dim_Z(\lambda) = \frac{{Z_J(\ell(\lambda))(I)} } {c_{PJ}(\ell(\lambda),\alpha=2)}.\end{gather*}
It may be interesting to notice that both the numerator and the denominator of this formula are not invariant under a global shift (translation of the partition $\ell(\lambda)$ by an arbitrary integer), but their ratio is invariant~-- this can be interpreted as a kind of ``gauge freedom'' in the writing of the $\SU(n)$ highest weight $\lambda$ as a partition.

Going back to our favorite \SU(3) example, we can check the consistency of the obtained result for $\chi_Z([1,1])^2$ in terms of dimensions. The (usual) dimensions of irreps labelled [0,0], [1,1], [0,3], [3,0], [2,2] are 1, 8, 10, 10, 27, and
the usual decomposition of $\chi([1,1])^2$ is compatible with the identity $8 \times 8 = ( 1 + 2 \times 8 + 1 \times 10 + 1 \times 10 + 1 \times 27$). The zonal dimensions of the same irreps are 1, 15/2, 7, 7, 35/2, and the decomposition of $\chi_Z([1,1])^2$ implies
\begin{gather*} \left(\frac{15}{2}\right)^2 = \frac{25}{12} +\frac{12}{5} \times\frac{15}{2}+\frac{4}{3} \times 7 +\frac{4}{3} \times 7 +1 \times \frac{35}{2}.\end{gather*}

Finally, still another way to calculate a structure constant is to use the ($\alpha$-dependent) Hall inner product, see for instance \cite{Macdonald:book} or \cite{Stanley:JackSymFun}, for which Jack polynomials -- and in particular zonal polynomials -- are orthogonal. Using it, one can see for instance that the ``multiplicity'' of $[1,1]$ in the decomposition of the square of $[1,1]$, or of $\{3,2,1\}$ in the square of $\{2,1,0\}$, which is $2$ in the usual case ($\alpha = 1$), and $12/5$ in the zonal-P case ($\alpha=2$) is more generally (i.e., in the Jack-P case) equal to
\begin{gather*}
\frac{\langle J_P(\{2, 1, 0\}) J_P(\{2, 1, 0\}), J_P(\{3, 2, 1\}) \rangle} {\langle J_P(\{3, 2, 1\}), J_P(\{3, 2, 1\}) \rangle}
= \frac{6 \alpha \left(2 \alpha ^2+11 \alpha +2\right)}{(\alpha +2) (2 \alpha +1) (3 \alpha +2)}.
\end{gather*}
In particular we can also find in this way the decomposition of $J_P[\{2,1,0\}]^2$:
\begin{gather*}
\frac{3 \alpha ^2 (\alpha +3) (2 \alpha +1)}{(\alpha +1)^2 (\alpha +2)^2} J_P[{2,2,2}] +
\frac{6 \alpha \big(2 \alpha ^2+11 \alpha +2\big)}{(\alpha +2) (2 \alpha +1) (3 \alpha +2)}J_P[{3,2,1}]\\
\qquad{} + \frac{2 \alpha }{\alpha +1} J_P[{3,3,0}]+\frac{2 \alpha }{\alpha +1} J_P[{4,1,1}]+J_P[{4,2,0}]
\end{gather*}
and we recover the first and last decompositions already given in Section~\ref{zonal:structure constants} by setting $\alpha = 2$ or~$1$.

\subsection{Back to the PDF (symmetric case) and to the ``volume function''}

In the Hermitian case it is known that one may associate with a given admissible triple $(\lambda, \mu, \nu)$ a convex polytope $\CH_{\lambda\mu}^\nu$, {the ``polytope of honeycombs''}, in a $d\le (n-1)(n-2)/2$-dimensional space~\cite{KT}. As recalled above, it is known that the function called ${\mathcal J}$ in \cite{CZ17}, which differs from the PDF $\p(\gamma)$ mainly by a Vandermonde factor $\Delta$, measures the volume of~$\CH_{\lambda\mu}^\nu$ and that it is also a good approximation of the LR multiplicity $N_{\lambda\mu}^\nu$ of that triple\footnote{Provided the triple is generic, i.e., provided ${\mathcal J}$ does not vanish.}. More precisely ${\mathcal J}$ is equal to the highest degree coefficient of the stretching (or LR) polynomial that gives the multiplicity when the triple $(\lambda, \mu, \nu)$ is scaled by a factor $s$, i.e., ${\mathcal J}$ is also the dominant coefficient of the Ehrhart polynomial of the polytope $\CH_{\lambda\mu}^\nu$. Since this multiplicity is known to be given by the number of {\it integral} points inside the polytope \cite{KT}, this property is just expressing that in the large $s$ limit, a semi-classical picture approximates well this number of points by the volume\footnote{For $n\times n$ matrices, and in the generic case, this is its $(n-1)(n-2)/2$ volume.} of the polytope.

In the symmetric case multiplicities are not integral, there are no honeycombs (at least the concept was not (yet?) generalized to cover this case), no polytope $\CH_{\lambda\mu}^\nu$, and no volume function either. However, as already mentioned, one expects the PDF $\p(\gamma)$, or rather the ``volume like function'' ${\mathcal J}=\p(\gamma)/\Delta=\rho$ (see~(\ref{pdfsym}),~(\ref{PDFgamma})), to measure, at least in the {\it generic case} where $\CJ$ does not vanish, the behavior of ``zonal multiplicities'' (i.e., appropriate zonal structure constants) under scaling.

In \cite[Section~4]{OO}\footnote{We thank Vadim Gorin for bringing our attention to that reference.} the function ${\mathcal I}_1$ is called a ``generalized Bessel function" and it is considered as a spherical function on the space $\O(n)\ltimes H(n, \R)$, $ H(n, \R)$ being the space of real symmetric matrices. This non-compact symmetric space is associated, by a standard limiting procedure (contraction) to the compact symmetric space~$\U(n)/\O(n)$. In turn, the spherical functions on the latter are essentially normalized zonal polynomials from which the generalized Bessel functions can be obtained by the same limiting procedure. This shows that the structure constants in the expansion of a product of two spherical functions ${\mathcal I}_1$ on spherical functions, i.e., our function ${\mathcal J}$, is the limit of the structure constants of the corresponding zonal polynomials. Some claims in that direction can also be found in~\cite{FG2}.

\looseness=-1 We are happy, in the present paper, to show ``experimentally'' that the overall features of the PDF $\rho$ computed in Section~\ref{sectionparticularcase} are consistent with the values obtained for zonal multiplicities, when the argument (written as a highest weight), is scaled. Remember that the list of eigenvalues $(1, 0, -1)$ chosen for the example studied in Sections~\ref{quaternion},~\ref{PDF} and~\ref{sectionparticularcase} differs from the partition $\{2,1,0\}$ (aka $[1,1]$ if reinterpreted in terms of \SU(3) highest weights) only by a constant shift\footnote{In the same way, the partitions that appear in the decomposition of the square of $Z_P(\{2,1,0\})$, namely $\{2,2,2\}$, $\{3,2,1\}$, $\{3,3,0\}$, $\{4,1,1\}$, $\{4,2,0\}$, whose reduction to $\SU(3)$ are $\{0,0,0\}$, $\{2,1,0\}$, $\{3,3,0\}$, $\{3,0,0\}$, $\{4,2,0\}$, and read $[0,0]$, $[1,1]$, $[0,3]$, $[3,0]$, $[2,2]$ in terms of Dynkin labels, differ by a constant shift from the following lists of eigenvalues (traceless condition for $\gamma_1,\gamma_2,\gamma_3$ see Section~\ref{PandQ}): $\{0,0,0\}$, $\{1,0,-1\}$, $\{1,1,-2\}$, $\{2,-1,-1\}$, $\{2,0,-2\}$, and appear in the Horn polygon of Fig.~\ref{polygonJzJz} as (special) points with respective coordinates $(p,q)$ given by $(0,0), (-1,0), (-3,2), (-3,-2), (4,0)$.}; we are therefore led to consider the behavior of the zonal structure constants that appear in the reduction of the square of $s\{2,1,0\}$, with a scaling factor $s=1,2,\ldots$, with the plot of the function~$\rho$ in terms of Dynkin labels. Using~\cite{MathematicaProgramsRC} we could perform exact calculations of multiplicities, i.e., obtain the decomposition of the square of $\chi_Z(s[1,1])$ up to the value $s=8$ of the scaling factor\footnote{The explicit reduction of $\chi_Z([s,s]) \times \chi_Z([s,s])$, for $s=1,2,\ldots$ can be obtained from the web site \url{http://www.cpt.univ-mrs.fr/~coque/Varia/JackZonalSchurResults.html}.\label{webpageRC}}.
The same considerations and calculations extend to the quaternionic case, where we compare multiplicities in the decomposition of $\chi_Q([8,8])^2$ with the function $\CJ_3$ computed in Section~\ref{quaternion}.

\begin{figure}[!tbp] \centering
 \includegraphics[width=11.7pc]{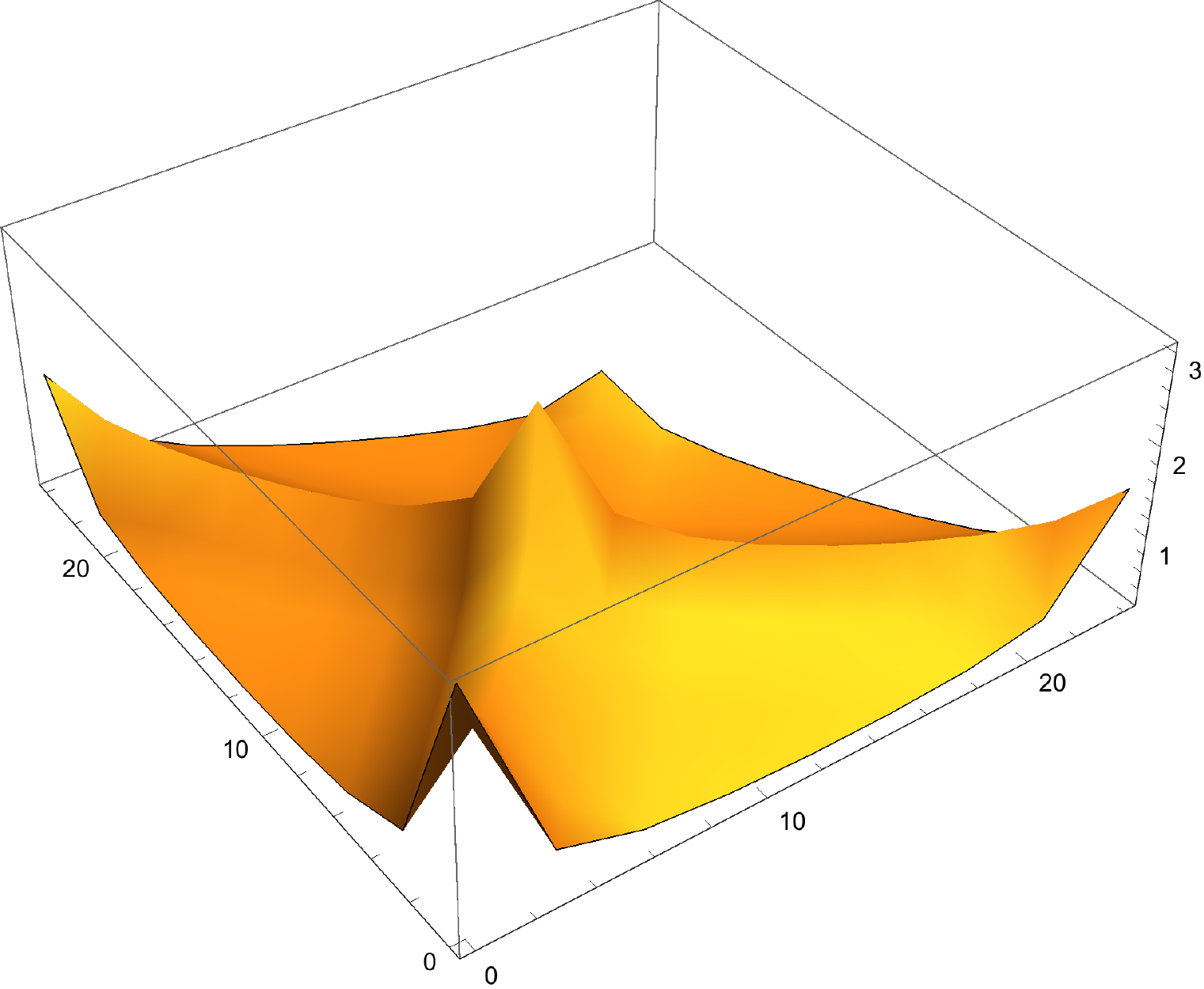}
 \includegraphics[width=11.7pc]{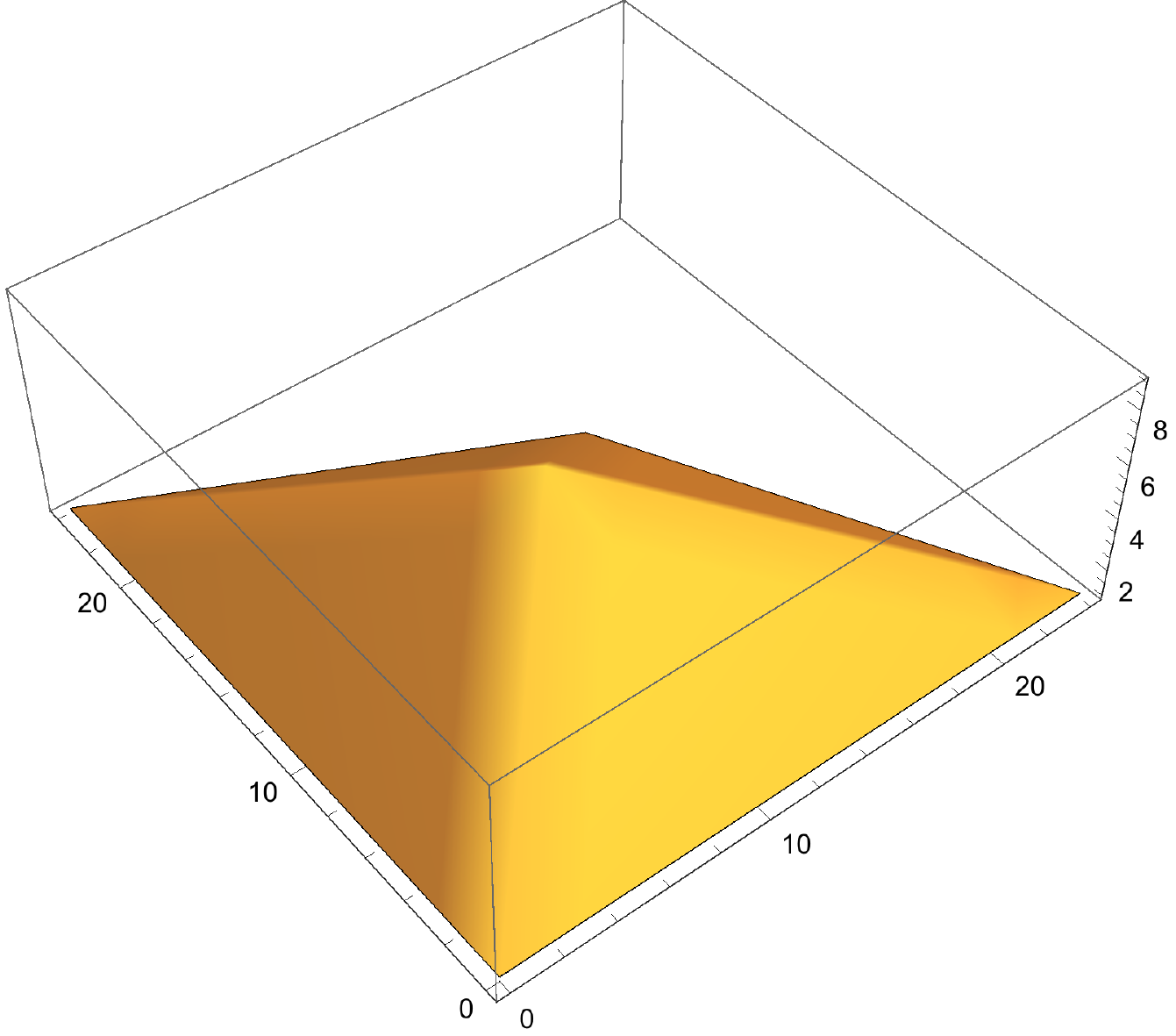}
 \includegraphics[width=11.7pc]{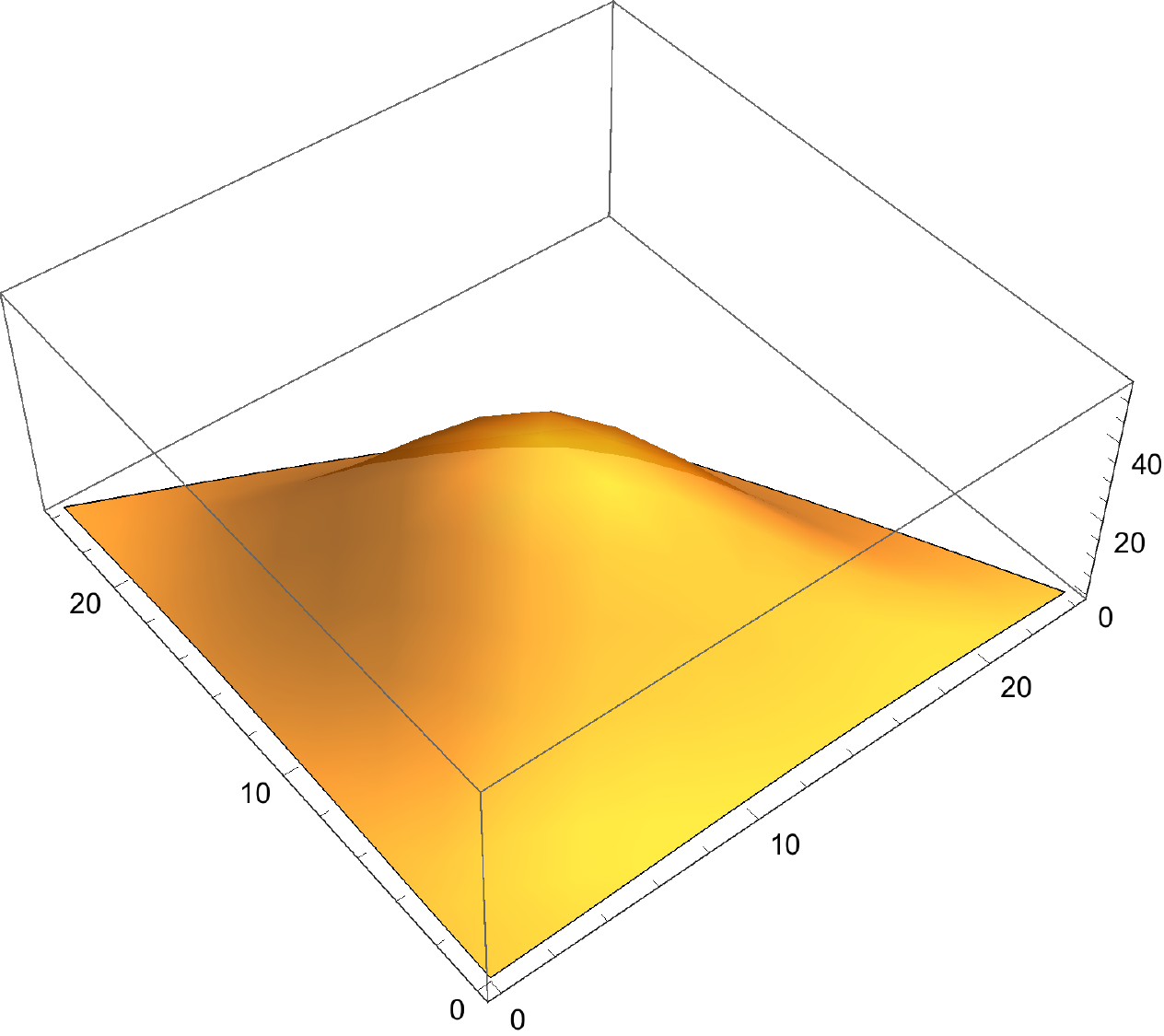}\\
 \includegraphics[width=11.7pc]{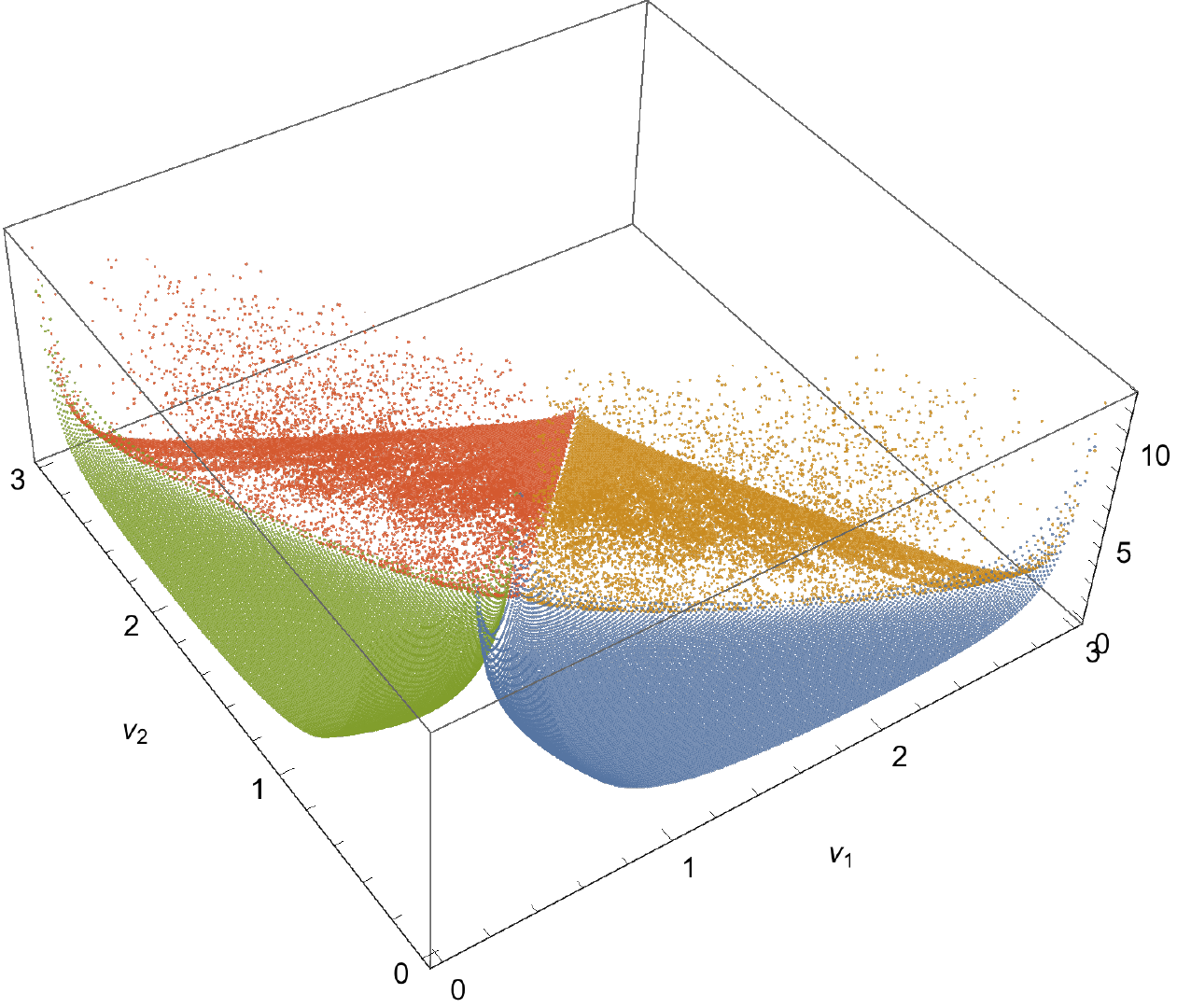} \includegraphics[width=11.7pc]{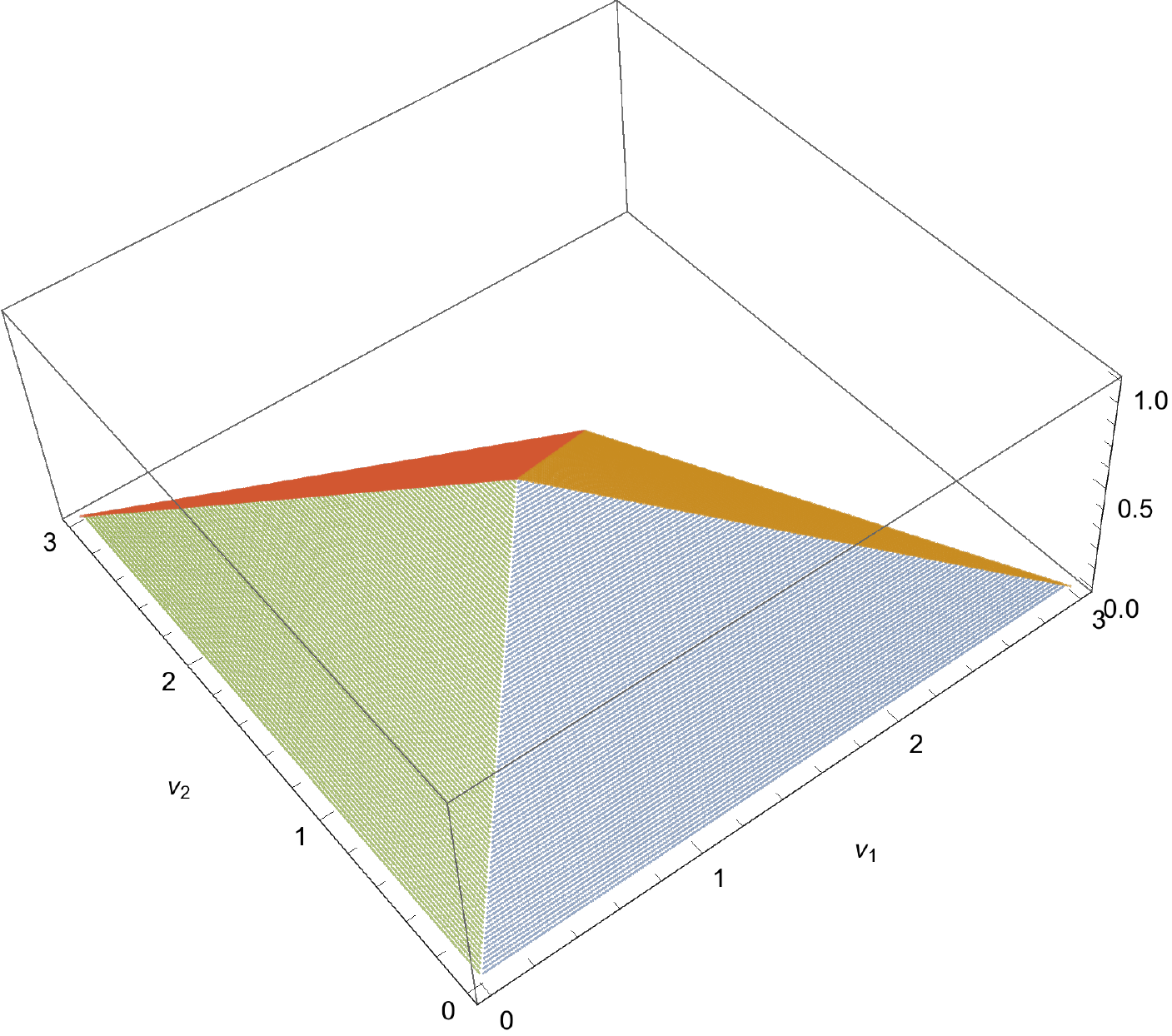} \includegraphics[width=11.7pc]{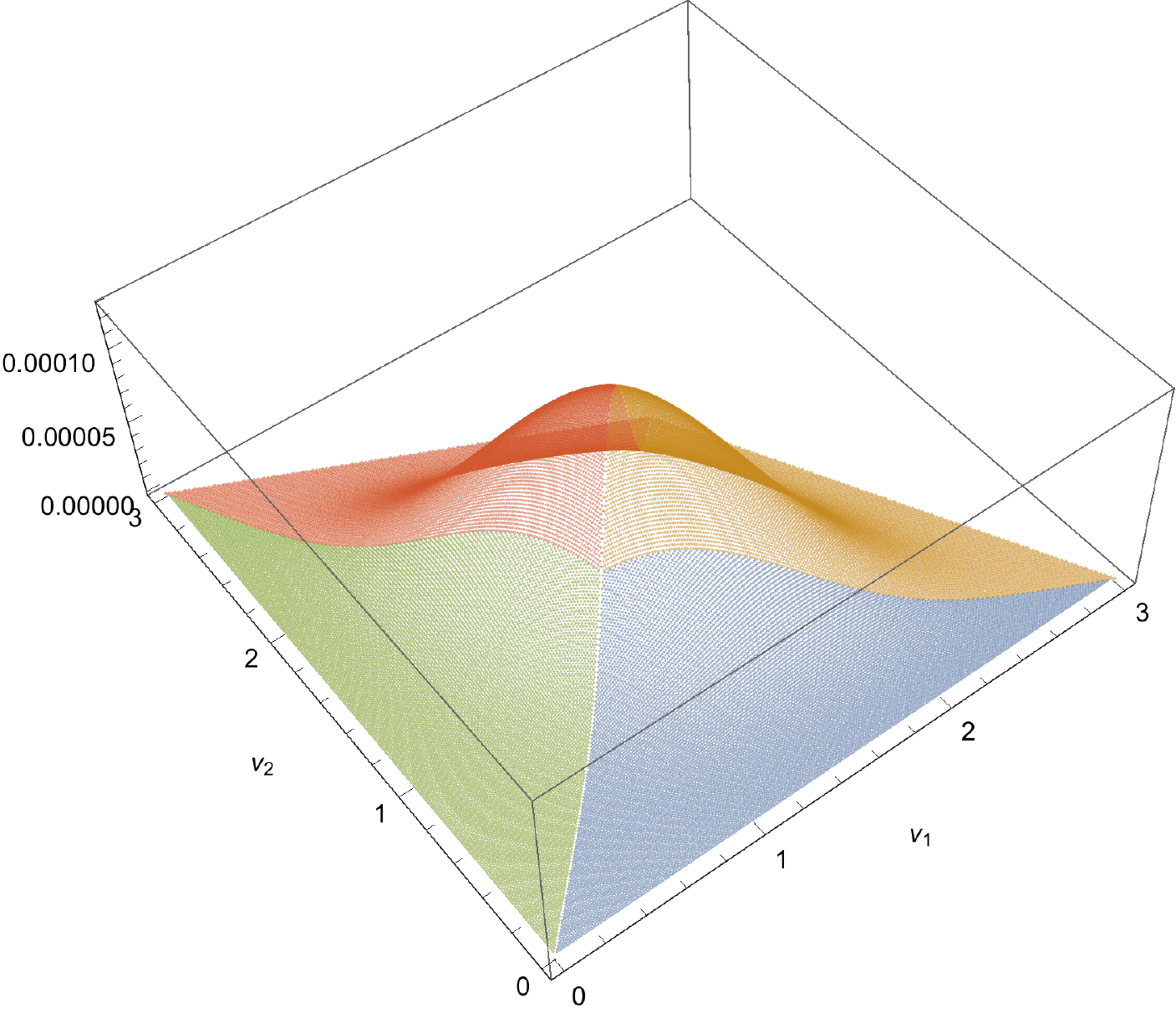}
 \caption{First line: Interpolating function surfaces representing an array of ``multiplicities''.
 Left: Decomposition of $\chi_Z([8,8])^2$ on $\SU(3)$ zonal characters.
 Middle: Decomposition of $\chi([8,8])^2$ on $\SU(3)$ Weyl characters.
 Right: Decomposition of $\chi_Q([8,8])^2$ on $\SU(3)$ (zonal) quaternionic characters.
 Second line: Results of our analytical calculations.
 }\label{multiplicitiesSquareOfZonal88}
\end{figure}

In Fig.~\ref{multiplicitiesSquareOfZonal88}, we compare the multiplicities obtained for the decomposition of $\chi_Z([8,8])^2$, $\chi([8,8])^2$ and $\chi_Q([8,8])^2$ with the plot of the volume functions $J$, ie $J = J_3$ calculated from~(\ref{CJ3}) in the quaternionic case, and $J = \rho$ in the symmetric case, calculated from~(\ref{pdf4}). After rescaling the coordinates $\nu_1$ and $\nu_2$ by a factor $s=8$, these plots give a good approximation of the volume function displayed on the second line, for the three cases. The four colors refer to the four sectors of Fig.~\ref{HornPolygonQuaternionic}, but here in Dynkin labels. The divergence at $\nu_1=\nu_2=1$ in the left figure, being a~very narrow and very high peak, is hardly visible.

\begin{figure}[!tbp] \centering
 \includegraphics[width=10pc]{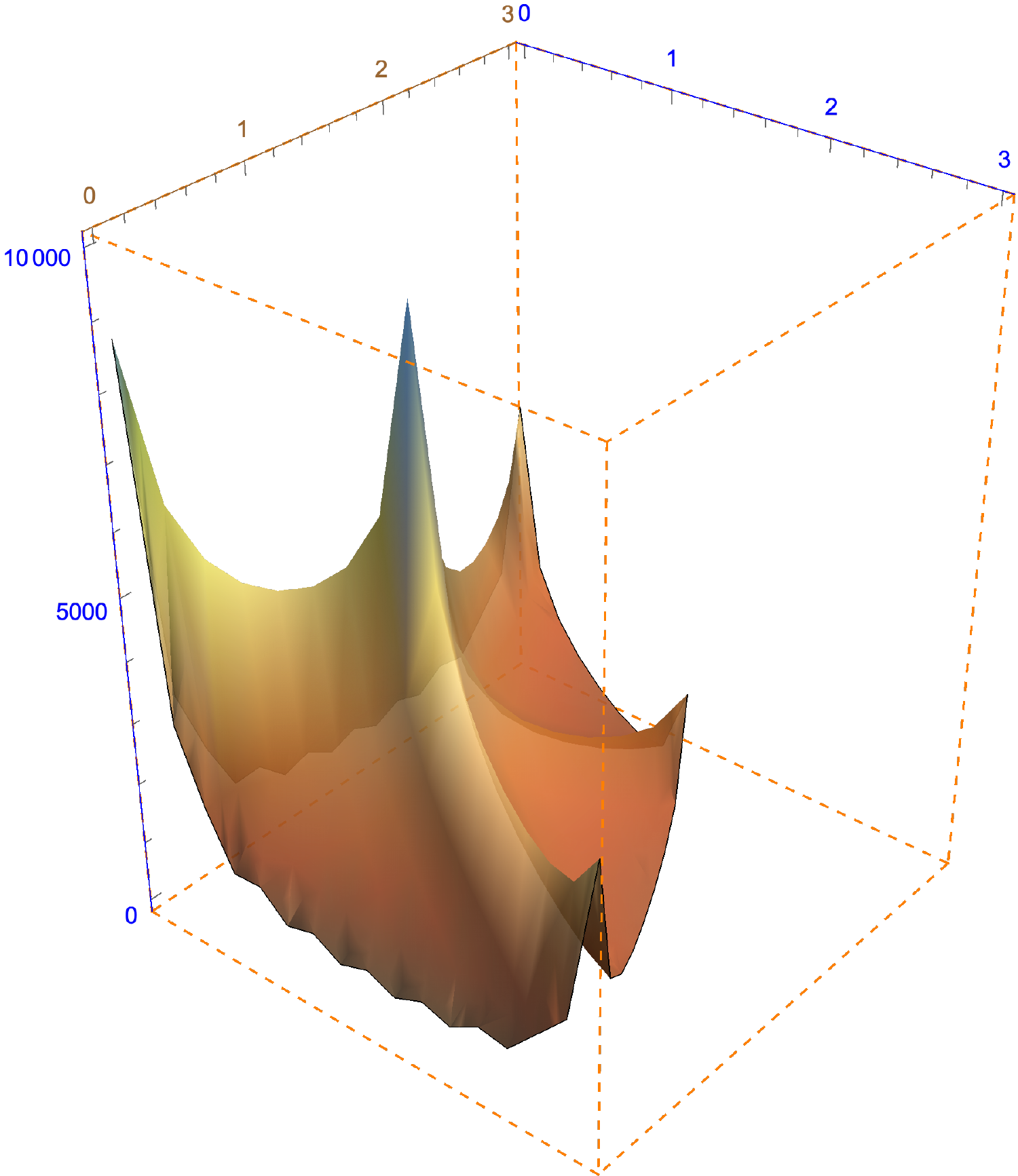}
 \includegraphics[width=10pc]{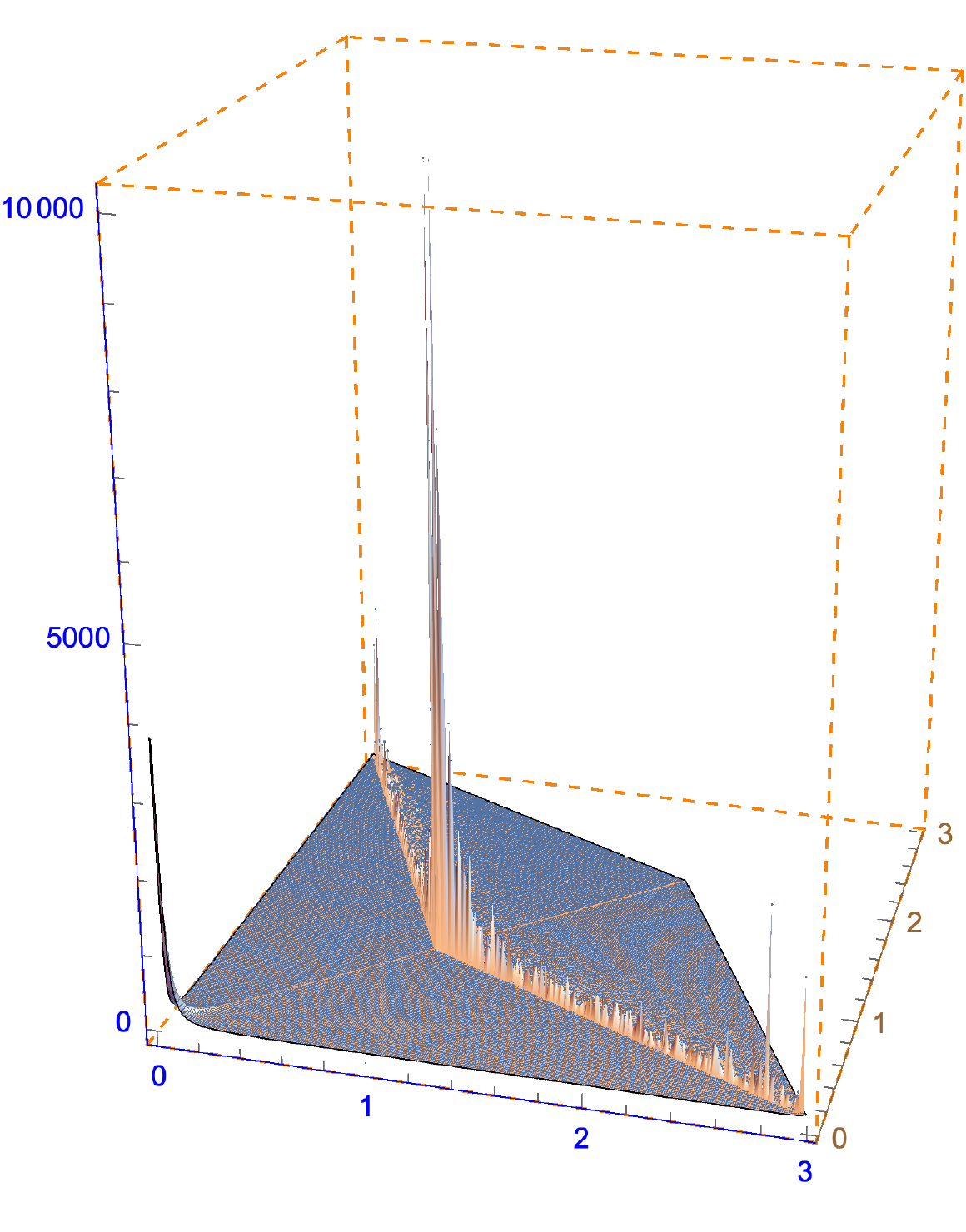}
 \includegraphics[width=10pc]{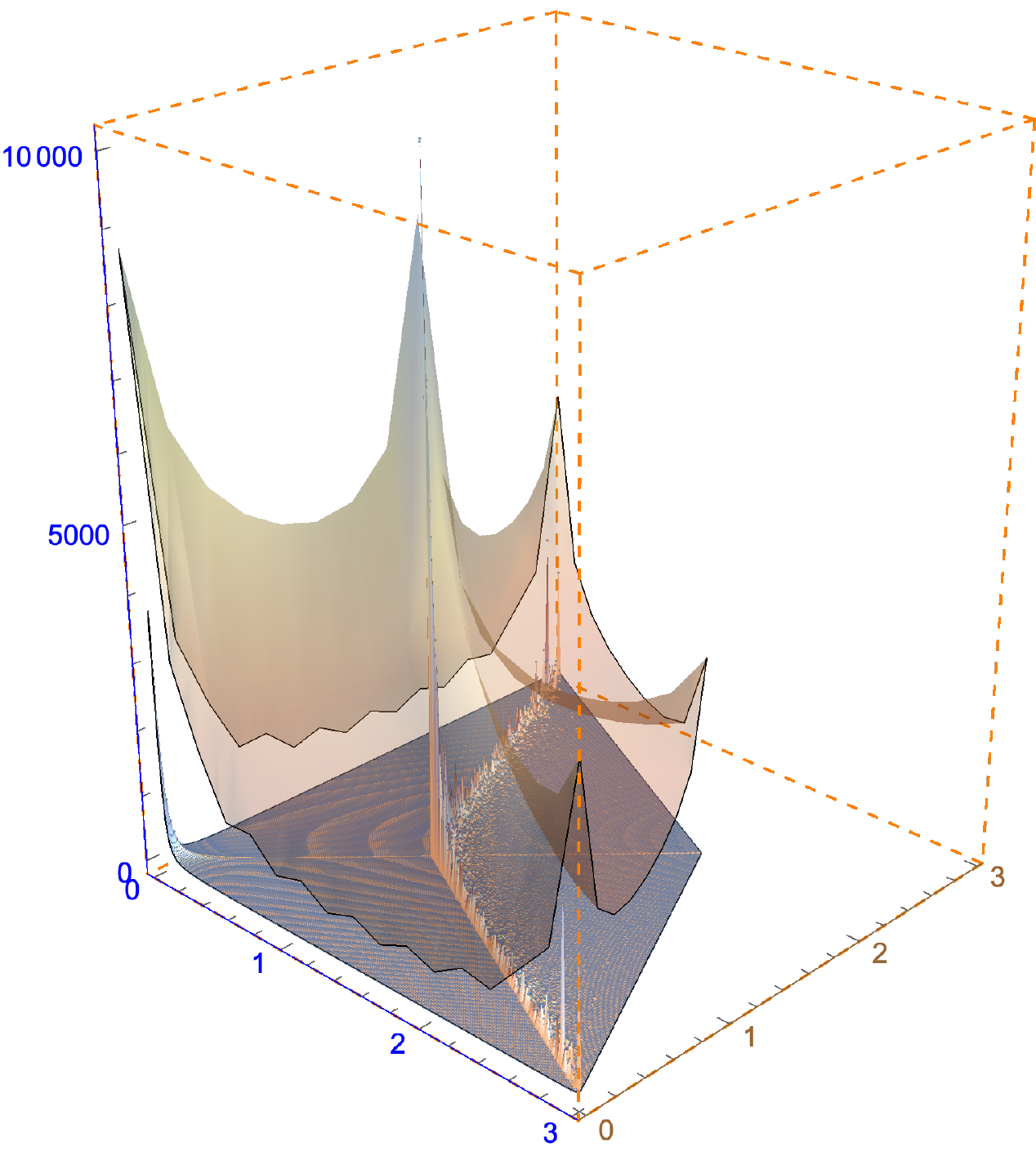}
 \caption{a)~Semi-transparent surface approximating the vertically rescaled zonal structure constants of $\chi_Z([8,8])^2$. b)~Point plot of the volume function $J=\rho$ --- because of the large vertical coordinate range, those parts of the surface lying far from the singularities look essentially flat and one only sees the singularities themselves (high values of $J$). c)~Superposition of a) and b). The horizontal coordinates (variables $\nu=(\nu_1, \nu_2)$ running between $0$ and $3$) are the Dynkin coordinates corresponding to $\alpha = \beta = \{1,0,-1\}$, obtained from the partition $\{2,1,0\}$ after a global shift by $-1$.}\label{FigXX}
 \end{figure}

In the symmetric case only, the structure of singularities is illustrated in the superimposed pictures of Fig.~\ref{FigXX} that display the volume function $J=\rho$, calculated from the integral (cloud of $59000$ points, from equation~(\ref{pdf4})), together with a vertically scaled version of the surface approximating the zonal structure constants $N_{\lambda \mu}^\nu$ (with $\lambda = \mu = (8,8)$).

We conclude that already with $s=8$, the classical limit provided by the ``volume'' approximates very well the distribution of multiplicities.

\section{Conclusion}
 To summarize:
 \begin{enumerate}\itemsep=0pt
\item[--] we have reproduced the main features of the Horn problem for symmetric matrices and understood the analytic origin of the singularities, at least for $n=3$, and in detail for the particular case of $A\sim B\sim J_z$;
\item[--] we have confirmed, at least in that particular case, that the divergences of the PDF occur on the same locus of non-analyticities as in the Hermitian and quaternionic cases;
\item[--] we have also confirmed numerically the connection between the ``volume functions'' and the (asymptotic) distribution of multiplicities in the product of zonal/Schur/quaternionic polynomials.
\end{enumerate}

This leaves, however, open several issues and room for further progress:
\begin{enumerate}\itemsep=0pt
\item[--] a more synthetic and general discussion of the singularities in the symmetric case would clearly be desirable. Can one understand their origin from a geometric point of view and assert a priori their location and nature without detailed calculations?
\item[--] what happens for higher $n$ and/or for generic $\bbeta$ (or $\alpha=2/\bbeta$)? The methods developed recently in \cite{Gorin2, Gorin1} should be helpful in that respect.
\end{enumerate}

The points we find most challenging are the following:
\begin{enumerate}\itemsep=0pt
\item[--] There is an enhancement of particular eigenvalues in the Horn spectrum of real symmetric matrices,
due to the divergences of the PDF. Is this enhancement observable in some physical process?
\item[--] The discussion of Section~\ref{zonalpolynomials} has pointed to an analogue of the volume function for real symmetric or quaternionic matrices: is there an underlying geometric interpretation to this ``volume''? is there a~geometric object generalizing the polytope $\CH_{\lambda\mu}^\nu$ of the Hermitian case, whose volume is computed there? on a representation theoretic side, what is the origin of the enhancement of certain multiplicities?
\end{enumerate}

We leave these questions to the sagacity of our readers \dots

\appendix

\section{Analysis of the singularities}\label{detail-sing}

In this subsection, we proceed to a detailed -- and lengthy~-- case-by-case analysis of the divergent singularities of $\rho$.

A preliminary observation is that, due to the $u\leftrightarrow z$ symmetry of our particular case, the $u<z$ and $u>z$ sectors will contribute equally:
\begin{gather*} \int\int {\rm d}u {\rm d}z \delta(R(u,z)) (2+u+z) \\
\qquad {}= \int\int_{u<z} {\rm d}u {\rm d}z \delta(R(u,z)) (2+u+z)+ \int\int_{u>z} {\rm d}u {\rm d}z \delta(R(u,z)) (2+u+z)\\
 \qquad{} = 2 \int\int_{u>z} {\rm d}u {\rm d}z \delta(R(u,z)) (2+u+z).
 \end{gather*}

This will be apparent in the following.
As explained in the previous subsection, we proceed heuristically, making appropriate Taylor expansions close to the singularities. We do it in detail in the first case (singularity along the dashed line), and are then more sketchy.
\begin{enumerate}\itemsep=0pt
\item Along the dashed line $p+1-q=0$, the integrand $\varphi$ has two non-integrable singularities: $\varphi\sim 1/z$ for $z\to z_{s_0}=0$, corresponding to $u_{1,2}(z)\to u_{s_0}=(p+3)/2$, and, by symmetry between $u$ and $z$, $\varphi\sim 1/(z_s-z)$ for $z\to z_s=(p+3)/2$. Thus the integral diverges and the PDF is infinite along the line.
\begin{itemize}\itemsep=0pt
\item For $q=p+1+\epsilon$ (i.e., close to and above the line), $z_s$ remains equal to zero, and we determine the common value of the two roots at that point by plugging a series expansion of the form $u_s= \frac{p+3}{2} +\alpha \epsilon+\cdots$ in the equation $R(u_s, z_s=0)=0$, whence
\begin{gather*} u_{1, 2}(0) = u_s= \frac{p+3}{2} +\frac{ 5+p}{2 (9+p)} \epsilon+\cdots.\end{gather*}
\item For $ z$ close to $z_{s_0}=0$, we write $z=\zeta^2$, approximate
\begin{gather*} u_{1,2}(z)= u_s +\beta_{1,2} \zeta +\gamma_{1,2} \zeta^2 +O\big(\zeta^3\big),\end{gather*} and determine the coefficients of that expansion by plugging it again in the equation $R\big(u_{1,2}(z), z=\zeta^2\big)=0$. $\beta_1$ (resp.~$\beta_2$) is the negative (positive) root of an equation which, for $q=p+1+\epsilon$ reduces to
 \begin{gather} \label{beta2}\beta^2=-\frac{2 (3 + p) (1 + p) (7 + p)^2} {(9 + p)^3}\epsilon +O\big(\epsilon^2\big)\end{gather}
 and the coefficient $\gamma$ is given by
 \begin{gather} \label{gamma12} \gamma_{1,2}= -\frac{p^3+5 p^2-13 p +15\pm 4\sqrt{2} (p+3) (p+7) \sqrt{-p-1} }{2 (p+9)^2} +O(\epsilon).\end{gather}
 \item We then expand the denominator $R'_u$ of~(\ref{phi}), for $ z$ close to $z_{s_0}=0$ and $\epsilon$ small and of order $\zeta^2$, as
 \begin{gather} R'_u\big(u_{1,2},z=\zeta^2\big) = 2 \sqrt{2} \sqrt{-p-1} (p+3) (p+7) \left(\sqrt{\frac{p+9}{p+3}} \sqrt{\epsilon }\mp 2
 \zeta \right) \zeta \nonumber\\
 \hphantom{R'_u\big(u_{1,2},z=\zeta^2\big) =}{} + o\big(\epsilon, \zeta\sqrt{\epsilon}, \zeta^2\big).\label{Rpu}\end{gather}
 \item and we finally derive the divergence of $\int_0^{z_s} {\rm d}z\varphi $ at the lower end point $0$ as
\begin{gather*} \int_{0}{\rm d}z \varphi \Big|_{\rm div} =\sum_{u_1,u_2}\int_0 \frac{2 \zeta {\rm d}\zeta (p+7)/2}{|R'_u(u_{1,2},z)|}\Big|_{\rm div}=\inv{4\sqrt{-2(p+1)}(p+3) }
 |\log\epsilon|.\end{gather*}
\end{itemize}

As explained above, the divergence of the integral at its other end point, obtained by the symmetry $u\leftrightarrow z$, contributes the same amount.

Thus the total divergence of $\rho(p,q)$ as the dashed line is approached from above is
\begin{gather*} \rho(p,q)\big|_{p-q+1 \to 0} \approx \rho_{\rm div}:= \inv{ \pi^2 \sqrt{2} (p+3) \sqrt{-(1+p)} }\big|\log |p-q+1 | \big|.\end{gather*}
This formula is well verified on numerical data, see Fig.~\ref{logdashedline}. Note that the logarithmic behavior is enhanced in the approach to $p=-1$, $q=0$ {and to $p=-3$, $q=-2$, at the southern end of the plot. } This will be reconsidered in detail in items~4 and~7 below.

\begin{figure}[!tbp] \centering
\includegraphics[width=25pc]{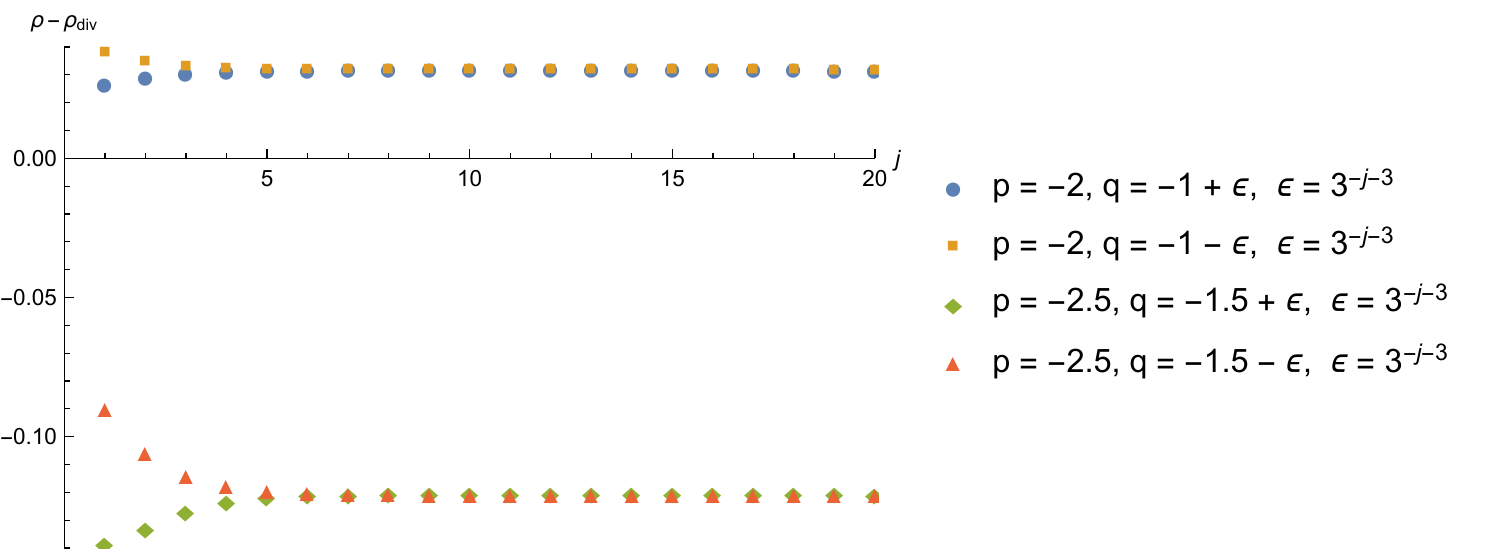}\qquad
\caption{$\rho-\rho_{\rm div}$ above (blue disks) or below (yellow squares) the dashed line for $p=-2.5$, $q=p+1\pm\epsilon$, $\epsilon=3^{-j-3}$, plotted against $j=1, \dots,20$.}\label{logdashedline}
\end{figure}

One also checks that the same formula applies to the approach to the dashed line from below ($\epsilon <0$). The relevant expressions are
\begin{gather*} z_{s_1}= \frac{(p+9) |\epsilon| }{4 (p+3)}, \qquad u_{s_1}=\frac{p+3}{2}-\frac{(p+5) |\epsilon| }{2 (p+9)}, \\
u_{1,2}\big(z=z_{s_1}+ \zeta^2\big)= u_{s_1}\mp \beta \zeta +\gamma_{1,2}\zeta^2+\cdots\end{gather*}
with $\beta>0$ and $\gamma_{1,2}$ as given in (\ref{beta2}) and (\ref{gamma12}), and the same expression for $R'_u$ as in~(\ref{Rpu}):
\begin{gather*} R'_u\big(u_{1,2}, z_{s_1}+\zeta^2\big)
\approx \pm 2\zeta\sqrt{-2(1+p)} (p+7)(p+3) \big({-}\sqrt{(p+9)/(3+p)}\sqrt{|\epsilon|} + 2 \zeta\big)\end{gather*}
and hence the same divergence as above the line, see Fig.~\ref{logdashedline} for an illustration at $p=-2.5$.

\item By a similar discussion, one finds that when $q\to 0$ with $-4\le p\le p_0=-1.21891$ (i.e., in region~I), the singularity at $z_{s_0}=0$ is integrable while that at $z_{s_1}\to 1$ gives rise to a~divergence of the integral. We write for short $z_s=z_{s_1}$. For $|q|$ small,
\begin{gather*} z_s = 1 +2q/(4+p) +O\big(q^2\big), \qquad\! u_s =u_1(z_s)=u_2(z_s)= 1 -8q/(p (4+p))+O\big(q^2\big) ; \end{gather*}
 as $z\to z_s$, we write $z= z_s-\zeta^2$, and the two roots $u_{1,2} \approx u_s \mp \beta_1 \zeta - \gamma_{1,2} \zeta^2$, for some computable coefficients $\beta_1$, $\gamma_{1,2} $, so that
 \begin{gather*} R'_u\big(u_{1,2},z_s-\zeta^2\big) \approx \mp 32 |p| \big(\sqrt{4+p} \zeta + \sqrt{-2q} \big) \zeta,
 \end{gather*} whence a divergence of the integral
 \begin{gather*} \int^{z_{s}} {\rm d}z \varphi(z) \sim \inv{ 4 |p| \sqrt{4+p}} |\log |q||.\end{gather*}

The same applies in region II, $ p_0=-1.21891\le p < 1$, where now four values of $z_s$ exist (see Section~\ref{subsectionroots} and Table~\ref{Table1}), but the singularities of~$R'_u$ at the points $z_{s_0}=0$, $z_{s_1}$ and~$z_{s_2}$ are of inverse square root type, hence integrable, and only the linear vanishing of~$R'_u$ at~$z_{s_3}$ matters.
Hence
\begin{gather*}
 \rho(p,q)\big|_{q \to 0\atop -4<p<-1} \approx \rho_{\rm div}:= \inv{2 \pi^2 |p| \sqrt{p+4} }\big|\log |q | \big|.\end{gather*}
See Fig.~\ref{q=0-line} for comparison with numerical data.

As $p\to -1$, however, the singularities at $z_{s_0}=0$ and $z_{s_2}\to 1$ become sharper and sharper, resulting in a stronger divergence at $z=1$, see below item~4.

\begin{figure}[!tbp]\centering
\includegraphics[width=20pc]{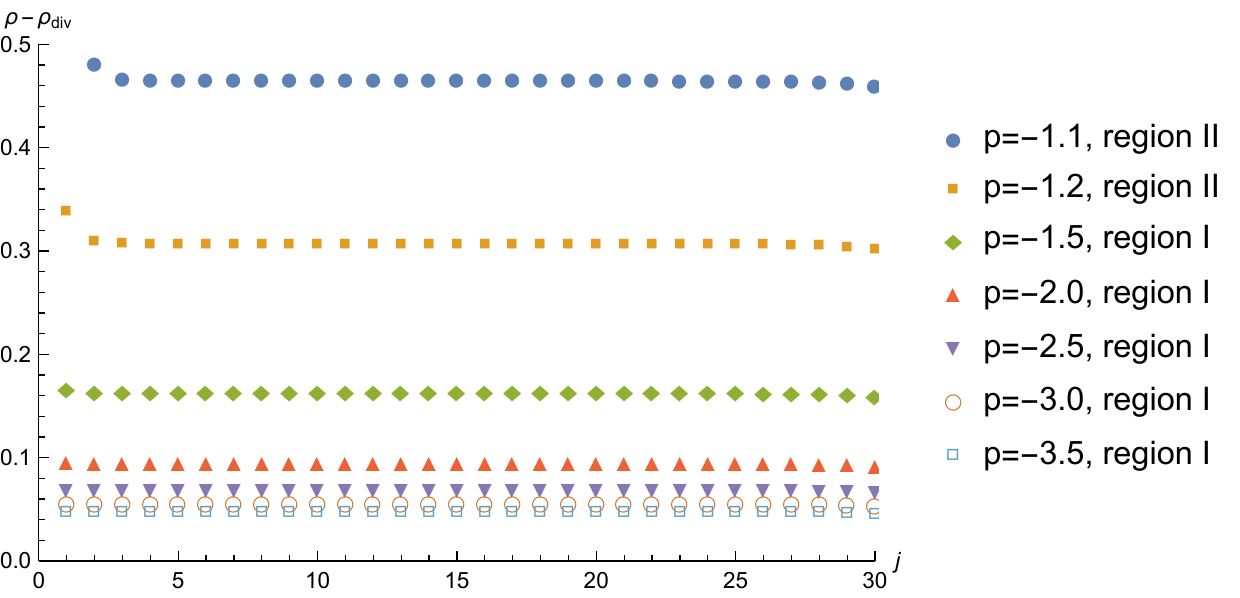}
\caption{$\rho- \rho_{\rm div}$ for $p=-3.5,\dots, -1.1$, $q=-3^{-j-1}$, plotted against $j=1, \dots,30$.}\label{q=0-line}
\end{figure}

 \item In region IV, ($-1<p<0$), $q$ is small ($|q|<1/10$), and the relevant values of $(z_s,u_s)$ are well approximated, up to higher powers of $q$, by
\begin{gather*}
 z_{s_2} \approx p+1+\alpha_{z_2} q \quad \text{with} \ \alpha_{z_2} \ \text{the largest root of} \ -4+p+p^2+2\alpha p (p+1) -4p\alpha^2, \\
 \text{i.e.}, \ \alpha_{z_2}=\big(p + p^2 - \sqrt{-16 p + 5 p^2 + 6 p^3 + p^4}\big)/(4 p);\\
 u_{s_2} \approx 1+ \alpha_{u_2} q,\qquad \alpha_{u_2}=\frac{1}{4}\big(1+\big(8+5 p+p^2\big)/\sqrt{-16 p+5 p^2+6 p^3+p^4}\big) ;\\
z_{s_4}\approx 1+ \frac{2q}{p+4},\qquad u_{s_4}\approx 1-\frac{8}{p (4+p)}q ;\\
z_{s_5} \approx 1+\oh q, \qquad u_{s_5}\approx p+1+\oh q \left(p+3+\frac{4}{p}\right).
 \end{gather*}
 (The value of $z_{s_3}$ where $u_3$ and $u_4$ merge, see Fig.~\ref{gallery}(g), is of no concern to us here, as $\varphi$ is integrable there.)

 As $q\to 0$ three divergences occur:
 \begin{enumerate}\itemsep=0pt
 \item[--] the pair of roots $(u_1,u_2)$ merge toward $u_{s_2}$ as $z\to z_{s_2}$; setting $\zeta^2=z-z_{s_2}$, as $\zeta\to 0$, $u_{1\atop 2}=u_{s_2} \mp \beta \zeta + \gamma_{1\atop 2} \zeta^2 +\cdots$ with $\beta^2$ and $\gamma_{1,2}$ determined by a series expansion of the equation $R(u,z)=0$. One finds that $\beta=O\big( |q|^\oh\big)$, while $\gamma_{1,2}=O(1)$. Consequently
\begin{gather*} R'_u\big(u_{1\atop 2}, z_{s_2}+\zeta^2\big) \approx \mp 8(p+4) (p+1)^\oh \\
\hphantom{R'_u\big(u_{1\atop 2}, z_{s_2}+\zeta^2\big) \approx}{} \times \big( \sqrt{-2p} |q|^\oh \big(p^4 + 6 p^3 + 5 p^2 - 16 p\big)^{\inv{4}} +2 p \zeta ) \big)\zeta,\end{gather*}
whence a divergence of the integral at its lower bound
\begin{gather*} \int_{z_{s_2}} \varphi(z) {\rm d}z \Big|_{\rm div }= \inv{8 |p| (p+1)^\oh} |\log |q| |; \end{gather*}
 \item[--] the pair of roots $(u_1,u_3)$ converges to $u_{s_5}$ as $z\to z_{s_5}$, and
 \begin{gather*} R'_u(u_{1,3},z) \approx 16|p| \big( (p+4) \sqrt{1+p} (1-z) +|q|\big);
 \end{gather*} this is just the image by the $u\leftrightarrow z$ symmetry of the latter, whence another contribution $ |\log|q|| / \big(8 |p| \sqrt{p+1}\big) $ to the integral of $\varphi$;
 \item[--] the pair of roots $(u_2,u_4)$ converges to $u_{s_4}$ as $z\to z_{s_4}$. By similar expansions near $z_{s_4}$ one finds that
\begin{gather*} R'_u\big(u_{2\atop4},z_{s_4}-\zeta^2\big)= 16 |p| \big({\mp} 2\sqrt{-2 q}+ (-(p+8) \mp 2(p+4)^\oh \zeta\big) \zeta,\end{gather*}
 whence a divergence of the integral
 \begin{gather*} \int^{z_{s_4}} {\rm d}z \varphi(z)\Big|_{\rm div} = |\log|q|| / \big(4 |p| \sqrt{p+4}\big).\end{gather*}
Thus as above in regions I and II, for $q\approx 0$, $\rho \sim C\log |q|$, but with a larger value of the coefficient, $C =\inv{2\pi^2 |p|}\big((1+p)^{-\oh} +(4+p)^{-\oh}\big)$. The agreement between numerical results and that coefficient $C$ is illustrated in Fig.~\ref{rhoIV-rhodiv} for $p=-0.8$. It deteriorates at small values of $|q|$ where the convergence of the integral is bad.
\end{enumerate}

\begin{figure}[!tbp] \centering
\includegraphics[width=16pc]{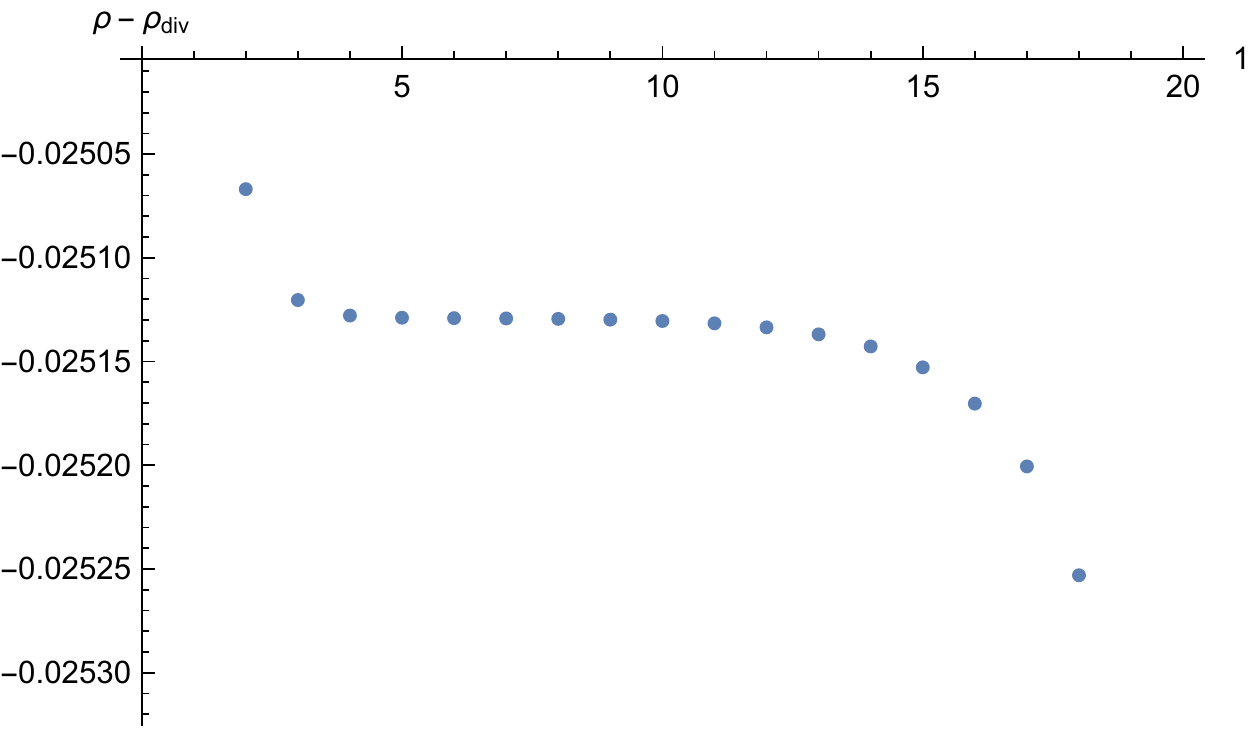}
\caption{$\rho- \rho_{\rm div}$ for $p=-0.8$, $q=-3^{-j-3}$, plotted against $j=1, \dots,18$. }\label{rhoIV-rhodiv}
\end{figure}

\item At $p=-1$, $q=0$, one can see that, as $z\to 0$, $u_1$ and $u_2$ approach~1 with the same slope, so that we have now $|u_1-u_2| \sim \alpha z^{3/2}$, causing a strong divergence of $\rho$. More precisely, for $p=-1$ and $q$ small, $z_{s_1}= -q+O\big(q^2\big)$ and if $z-z_{s_1}=\zeta^2\to 0$, one finds
\begin{gather*} R'_u| \approx 36 \sqrt{2} \left(\zeta^2 + \frac{4}{3} |q|\right)\zeta,\end{gather*}
whence a contribution to the divergent part of $\int_{z_{s_1}} {\rm d}z\varphi(z) $ equal to $\frac{\pi}{4\sqrt{6} |q|^\oh}$. The coalescence of roots $u_1$ and $u_3$ towards $u_{s_4}\approx 0$ at $z_{s_4}\approx 1$ gives rise to the same divergence (by the $u-z$ symmetry once again), while there is a weaker (logarithmic) divergence coming from~$u_2 $ and~$u_4$ merging to $u_{s_3}\approx 1$ as $z\to z_{s_3}\approx 1$ (see Fig.~\ref{gallery}(f)). In total, we have
\begin{gather*}
 \rho(p, q)\Big|_{p=-1\atop q\to 0_-} \approx \rho_{\rm div}:= \frac{1}{\sqrt{6}\pi |q|^\oh} \end{gather*}
up to a subdominant $\log|q|$ term. See Fig.~\ref{approachto-10-40}(left) for a numerical plot of $\rho/\rho_{\rm div}$ converging to~1.

 \begin{figure}[t] \centering
 \includegraphics[width=16pc]{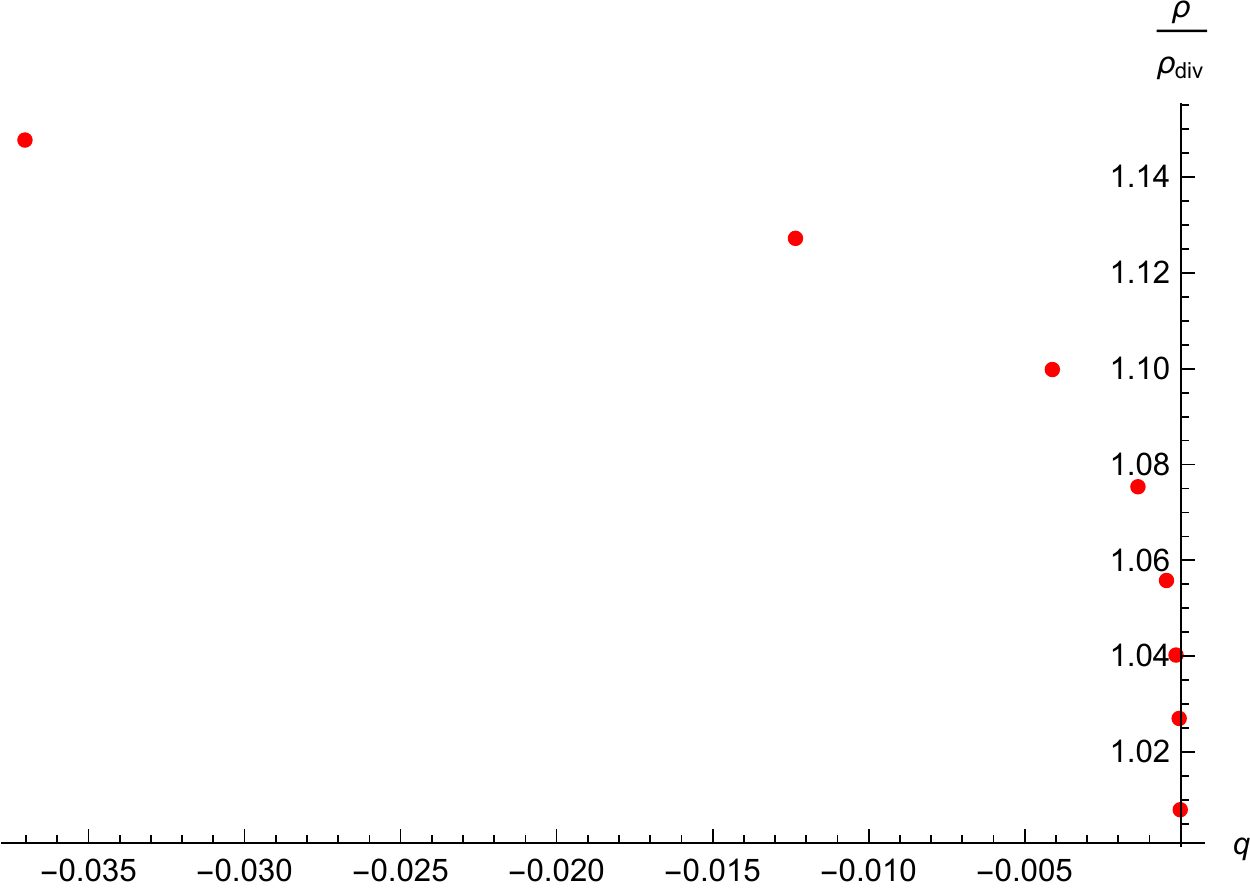}\qquad\qquad \centering\includegraphics[width=18pc]{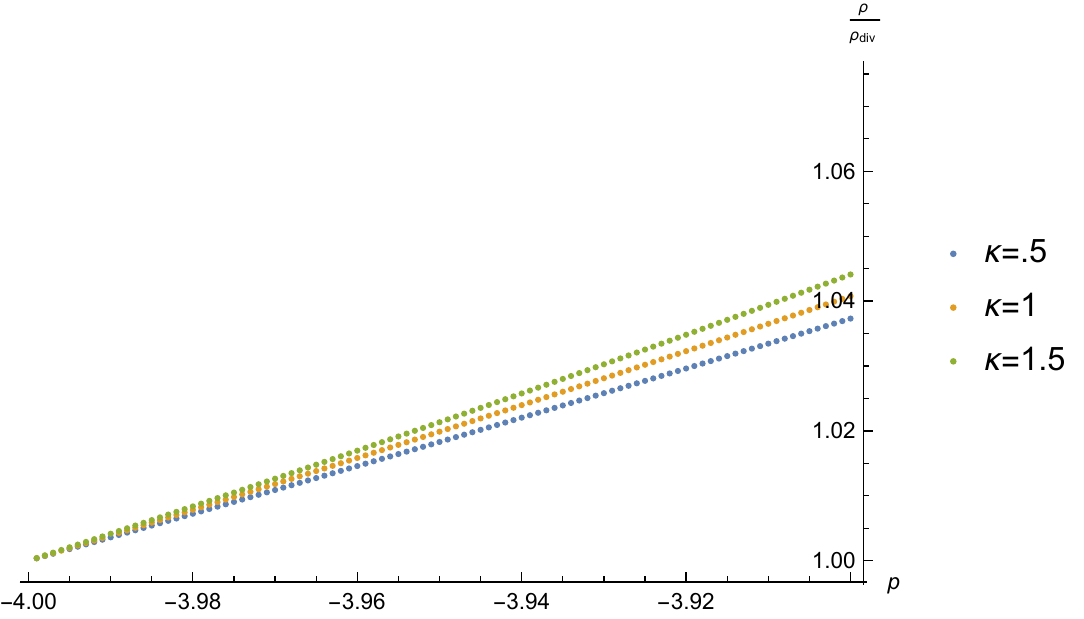}
\caption{Left: $\rho/ \rho_{\rm div}$ for $p=-1,\ q= -3^{-j}\to 0$, $j=3,\dots,10$, plotted against $q$. Right: Plot of $\rho/ \rho_{\rm div}$ as a function of $p$, in the approach of $(-4,0)$ along the line $p+q/\kappa=-4$ for three values of $\kappa$.}\label{approachto-10-40}
\end{figure}

\item Divergence at $(p,q)=(-4,0)$. Let us approach that singular point in a linear way, letting $p+4 = \epsilon$, $q = - \kappa \epsilon$, $\epsilon\to 0$, with $ 0<\kappa< 2$ (so as to remain within the Horn domain); then to leading order in $\epsilon$, the $z$-integration runs between $z_{s_0}=0$ and $z_{s_1}=(2 - \kappa)/(2 + 3 \kappa)$. Solving $R=0$ in that limit, one finds that the portrait $(u_1(z),u_2(z))$ forms a vanishingly thin
ellipse-like curve along the diagonal:
\begin{gather*} u_{1\atop 2}(z) = z \mp \sqrt{2+3\kappa} \left(\frac{(1-z) z (z_{s_1}-z)}{1+z} \right)^\oh\epsilon^\oh
\end{gather*} to the first non-trivial order in $\epsilon$, and plugging this expression into $R'_u$ yields
\begin{gather*} |R'_u(u_{1,2})|=32 \sqrt{2+3\kappa} (1+z)^\frac{3}{2} \big((1-z) z (z_{s_1}-z)\big)^\oh \epsilon^\oh +O(\epsilon),\end{gather*}
whence a divergence of $\rho$ as $|q|^{-\oh}$
\begin{gather*}
\rho_{\rm div}:= \inv{4\pi^2}\sqrt{\frac{\kappa}{2+3\kappa}} \inv{|q|^\oh} \int_0^{z_{s_1}} \frac{{\rm d}z}{\sqrt{z(1-z)(1+z)(z_{s_1}-s)}} \end{gather*}
in very good agreement with numerical data, see Fig.~\ref{approachto-10-40}(right).

\item Divergence at $p=q=0$. We let $(p,q)$ approach $(0,0)$ in region VI, for example along the cubic $ \kappa p^3+27q^2=0$, with $\kappa<4$; then the two end points $z_{s_2}$ and $z_{s_3}$ of the integral go to 1 as $p=-\epsilon^2$ goes to~0, i.e., $z_{s_i}\approx 1 - \alpha_{z,i} \epsilon^2 +O\big(\epsilon^3\big)$ where $\alpha_{z,i}$, $i=2,3$, are the second and third largest roots of $27 \alpha_z(1-\alpha_z)^2 -\kappa$. Both $z$ and $u$ are thus confined in an interval of size $\epsilon^2$ near 1, and solving the equation $R=0$ in the rescale variable $z=1-\zeta \epsilon^2$ and plugging into $R'_u$, one finds that the latter has a limiting shape described by the elliptic curve
\begin{gather*} R'_u(u_{2,3})=\pm\frac{2^6}{3\sqrt{3}} \epsilon^4 \sqrt{\zeta \big(27 \zeta(1-\zeta)^2-\kappa\big)}\end{gather*}
so that
 \begin{gather*} 
 \rho_{\rm div}:=\frac{\sqrt{3} \kappa^\inv{3}}{4\pi^2 |q|^\frac{2}{3}} \int_{\alpha_{z,3}}^{\alpha_{z,2}}\frac{{\rm d}\zeta}{\sqrt{\zeta \big(27 \zeta (\zeta -1)^2-\kappa\big)}}.
 \end{gather*}
 This behavior is again well supported by numerical calculations, see Fig.~\ref{limit0032}(left).

\begin{figure}[!tbp]\centering
 \includegraphics[width=17pc]{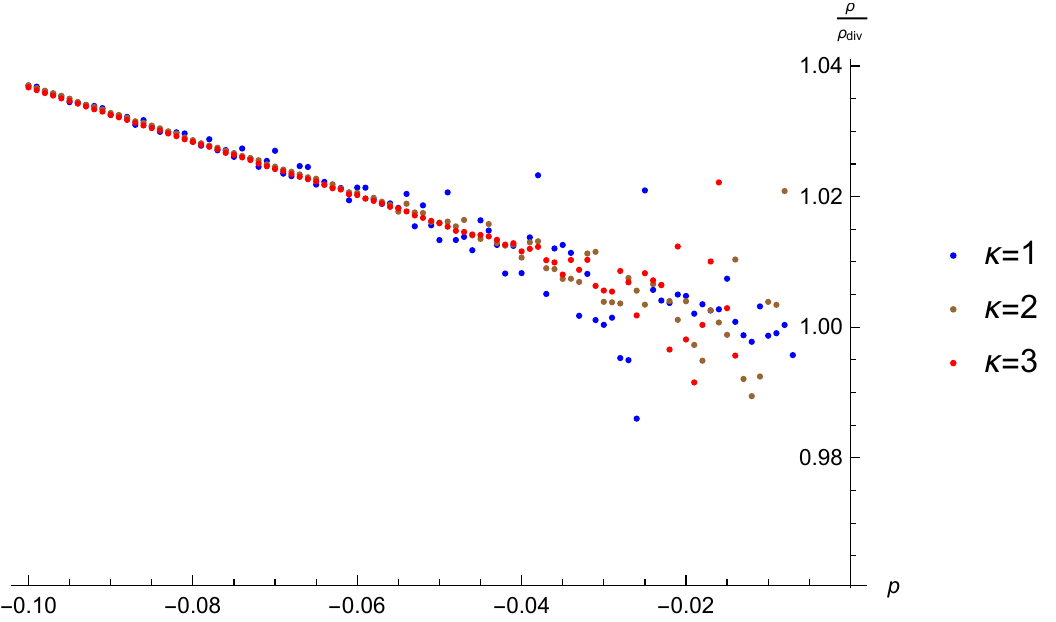} \qquad \includegraphics[width=17pc]{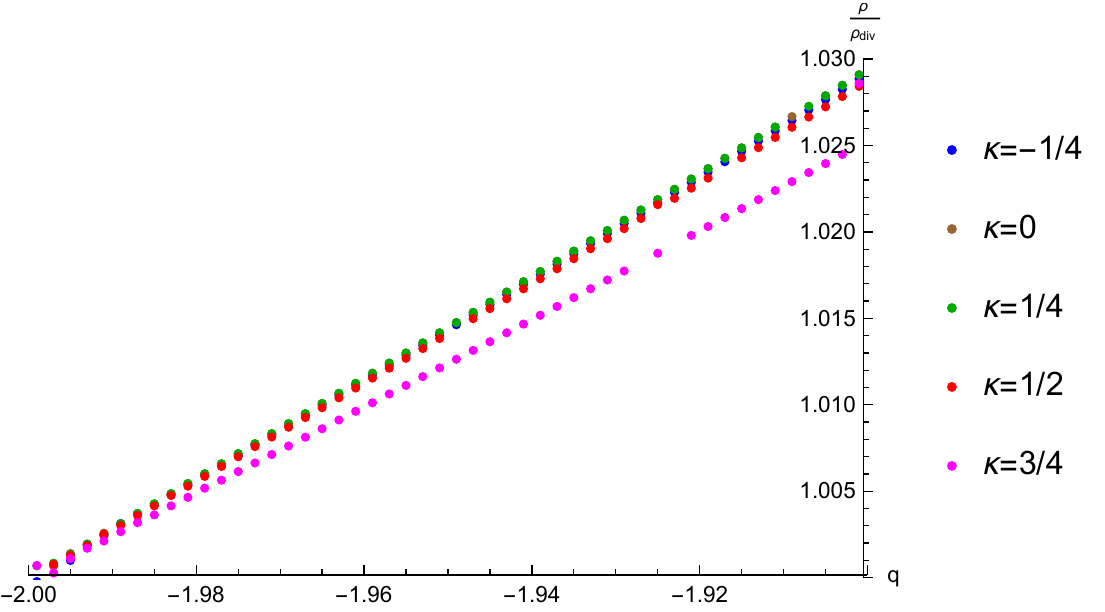}
\caption{Left: Plot of $\rho(p,q)/\rho_{\rm div}$ in the approach of $(0,0)$ along the cubic $\kappa p^3+27 q^2=0$ for $\kappa=1, 2, 3$ (in region VI). Convergence of the integral deteriorates for
small values of $p$. Right: Plot of $\rho(p,q)/\rho_{\rm div}$ in the approach of $(-3,-2)$ along the line $p-q \kappa=-3+2\kappa$ for five values of $\kappa$.}\label{limit0032}
 \end{figure}

\item Divergence at $p=-3$, $q=-2$. If one approaches that corner of the domain along lines $p=-3+\kappa \epsilon$, $q=-2+\epsilon$ with $-\oh<\kappa< 1$ so as to remain in region~I, one finds that $z_{s_0}=0$ and $z_{s_1}= \inv{6} (1+2\kappa)\epsilon$ to the lowest order in $\epsilon$, so that both $z$ and $u$ remain small of order $\epsilon$. Solving the equation $R=0$ to order $\epsilon^3$, one finds that
 \begin{gather*} u_{1\atop 2}(z)= z_{s_1}- z \pm \frac{2}{9} \big(12\epsilon (1-\kappa)\big(z(z_{s_1}-z)\big)\big)^\oh\end{gather*} so that
 \begin{gather*} |R'_u(u_{1,2})|\approx 32 \sqrt{3(1-\kappa)} \big(z(z_{s_1}-z)\big)^\oh \epsilon^\oh,\end{gather*} and the $z$-integration may be carried out, leading to
\begin{gather*}
\rho_{\rm div}:= \inv{\pi\sqrt{48 (1-\kappa)} } \inv{|q+2|^\oh}. \end{gather*}
 This is corroborated by the numerical calculation at various values of $\kappa$, see Fig.~\ref{limit0032}(right).
\end{enumerate}

The alert reader may wonder why the singularity along the line $p+q+1=0$ of the {\it upper} half-plane, (a reflection of the singularity along the dashed line of the $q<0$ half-plane) does not manifest itself along the {\it dotted} line of the lower half-plane. The reason is that, in that lower half-plane, the two $z_s$ that merge there are in fact irrelevant for $q<0$.

\subsection*{Acknowledgements}
Many thanks to Mich\`ele Vergne for challenging us to carry out the $n=3$ calculation {and for her patient explanations on the location of the singularities occurring in orbital integrals}, and to Michel Bauer for his encouragement and in particular for his constructive criticism of identity~(\ref{deltaR2}). Stimulating conversations with O.~Babelon, \'E.~Br\'ezin, P.~Di Francesco, J.~Faraut, V.~Gorin, S.~Majumdar and G.~Schehr are also acknowledged. We thank the referees for suggesting several editorial improvements.

\pdfbookmark[1]{References}{ref}
\LastPageEnding


\begin{thebibliography}{99}
\footnotesize\itemsep=0pt

\bibitem{Baratta}
Baratta W., SpecialFunctions.nb, {T}he University of Melbourne, 2008,
 \url{https://researchers.ms.unimelb.edu.au/~wbaratta/}.

\bibitem{BZ19}
Bauer M., Zuber J.-B., On products of delta distributions and resultants, in
 preparation.

\bibitem{BH1}
Br\'{e}zin E., Hikami S., An extension of the
 {H}arish-{C}handra--{I}tzykson--{Z}uber integral, \href{https://doi.org/10.1007/s00220-003-0804-x}{\textit{Comm. Math. Phys.}}
 \textbf{235} (2003), 125--137, \href{https://arxiv.org/abs/math-ph/0208002}{arXiv:math-ph/0208002}.

\bibitem{Gorin2}
Bufetov A., Gorin V., Fourier transform on high-dimensional unitary groups with
 applications to random tilings, \href{https://arxiv.org/abs/1712.09925}{arXiv:1712.09925}.

\bibitem{ConstantineZonalPoly}
Constantine A.G., Some non-central distribution problems in multivariate
 analysis, \href{https://doi.org/10.1214/aoms/1177703863}{\textit{Ann. Math. Statist.}} \textbf{34} (1963), 1270--1285.

\bibitem{MathematicaProgramsRC}
Coquereaux R., Mathematica package ``SymPol\$Package'', available at
 \url{http://www.cpt.univ-mrs.fr/~coque/Computer_programs/index.html} and
 \url{https://github.com/RobertCoquereaux/Lie-AffineLie-Representations}.

\bibitem{CMcSZ}
Coquereaux R., McSwiggen C., Zuber J.-B., On {H}orn's problem and its volume
 function, \href{https://arxiv.org/abs/1904.00752}{arXiv:1904.00752}.

\bibitem{CZNuclPhys}
Coquereaux R., Zuber J.-B., On some properties of {$\rm SU(3)$} fusion
 coefficients, \href{https://doi.org/10.1016/j.nuclphysb.2016.05.029}{\textit{Nuclear Phys.~B}} \textbf{912} (2016), 119--150,
 \href{https://arxiv.org/abs/1605.05864}{arXiv:1605.05864}.

\bibitem{CZ17}
Coquereaux R., Zuber J.-B., From orbital measures to {L}ittlewood--{R}ichardson
 coefficients and hive polytopes, \href{https://doi.org/10.4171/AIHPD/57}{\textit{Ann. Inst. Henri Poincar\'{e}~D}}
 \textbf{5} (2018), 339--386, \href{https://arxiv.org/abs/1706.02793}{arXiv:1706.02793}.

\bibitem{DemmelKoev}
Demmel J., Koev P., Accurate and efficient evaluation of {S}chur and {J}ack
 functions, \href{https://doi.org/10.1090/S0025-5718-05-01780-1}{\textit{Math. Comp.}} \textbf{75} (2006), 223--239.

\bibitem{MOPS}
Dumitriu I., Edelman A., Shuman G., M{OPS}: multivariate orthogonal polynomials
 (symbolically), \href{https://doi.org/10.1016/j.jsc.2007.01.005}{\textit{J.~Symbolic Comput.}} \textbf{42} (2007), 587--620,
 \href{https://arxiv.org/abs/math-ph/0409066}{arXiv:math-ph/0409066}.

\bibitem{Fa}
Faraut J., Horn's problem and {F}ourier analysis, \href{https://doi.org/10.2140/tunis.2019.1.585}{\textit{Tunis.~J. Math.}}
 \textbf{1} (2019), 585--606.

\bibitem{FeraySniady}
F\'{e}ray V., \'{S}niady P., Zonal polynomials via {S}tanley's coordinates and
 free cumulants, \href{https://doi.org/10.1016/j.jalgebra.2011.03.008}{\textit{J.~Algebra}} \textbf{334} (2011), 338--373,
 \href{https://arxiv.org/abs/1005.0316}{arXiv:1005.0316}.

\bibitem{FG2}
Frumkin A., Goldberger A., Diagonals of real symmetric matrices of given
 spectra as a measure space, \href{https://arxiv.org/abs/1505.06418}{arXiv:1505.06418}.

\bibitem{FG}
Frumkin A., Goldberger A., On the distribution of the spectrum of the sum of
 two {H}ermitian or real symmetric matrices, \href{https://doi.org/10.1016/j.aam.2005.12.007}{\textit{Adv. in Appl. Math.}}
 \textbf{37} (2006), 268--286.

\bibitem{Fu}
Fulton W., Eigenvalues, invariant factors, highest weights, and {S}chubert
 calculus, \href{https://doi.org/10.1090/S0273-0979-00-00865-X}{\textit{Bull. Amer. Math. Soc.}} \textbf{37} (2000), 209--249,
 \href{https://arxiv.org/abs/math.AG/9908012}{arXiv:math.AG/9908012}.

\bibitem{GelfandShilov}
Gel'fand I.M., Shilov G.E., Generalized functions, {V}ol.~{I}, {P}roperties and
 operations, Academic Press, New York~-- London, 1964.

\bibitem{Gorin1}
Gorin V., Marcus A.W., Crystallization of random matrix orbits,
 \href{https://arxiv.org/abs/1706.07393}{arXiv:1706.07393}.

\bibitem{HC}
Harish-Chandra, Differential operators on a semisimple {L}ie algebra,
 \href{https://doi.org/10.2307/2372387}{\textit{Amer.~J. Math.}} \textbf{79} (1957), 87--120.

\bibitem{BH2}
Hikami S., Br\'{e}zin E., W{KB}-expansion of the
 {H}arish-{C}handra--{I}tzykson--{Z}uber integral for arbitrary~{$\beta$},
 \href{https://doi.org/10.1143/PTP.116.441}{\textit{Progr. Theoret. Phys.}} \textbf{116} (2006), 441--502,
 \href{https://arxiv.org/abs/math-ph/0604041}{arXiv:math-ph/0604041}.

\bibitem{HualLoKen}
Hua L.-K., Harmonic analysis of functions of several complex variables in the
 classical domains, Amer. Math. Soc., Providence, R.I., 1963.

\bibitem{IZ}
Itzykson C., Zuber J.-B., The planar approximation.~{II}, \href{https://doi.org/10.1063/1.524438}{\textit{J.~Math.
 Phys.}} \textbf{21} (1980), 411--421.

\bibitem{JamesZonalPoly1}
James A.T., Normal multivariate analysis and the orthogonal group, \href{https://doi.org/10.1214/aoms/1177728846}{\textit{Ann.
 Math. Statist.}} \textbf{25} (1954), 40--75.

\bibitem{JamesZonalPoly2}
James A.T., The distribution of the latent roots of the covariance matrix,
 \href{https://doi.org/10.1214/aoms/1177705994}{\textit{Ann. Math. Statist.}} \textbf{31} (1960), 151--158.

\bibitem{JamesZonalPoly3}
James A.T., Zonal polynomials of the real positive definite symmetric matrices,
 \href{https://doi.org/10.2307/1970291}{\textit{Ann. of Math.}} \textbf{74} (1961), 456--469.

\bibitem{Kl}
Klyachko A.A., Stable bundles, representation theory and {H}ermitian operators,
 \href{https://doi.org/10.1007/s000290050037}{\textit{Selecta Math. (N.S.)}} \textbf{4} (1998), 419--445.

\bibitem{KT}
Knutson A., Tao T., The honeycomb model of {${\rm GL}_n({\mathbb C})$} tensor
 products. {I}.~{P}roof of the saturation conjecture, \href{https://doi.org/10.1090/S0894-0347-99-00299-4}{\textit{J.~Amer. Math.
 Soc.}} \textbf{12} (1999), 1055--1090, \href{https://arxiv.org/abs/math.RT/9807160}{arXiv:math.RT/9807160}.

\bibitem{KT00}
Knutson A., Tao T., Honeycombs and sums of {H}ermitian matrices,
 \textit{Notices Amer. Math. Soc.} \textbf{48} (2001), 175--186,
 \href{https://arxiv.org/abs/math.RT/0009048}{arXiv:math.RT/0009048}.

\bibitem{KTW01}
Knutson A., Tao T., Woodward C., The honeycomb model of {${\rm GL}_n({\mathbb
 C})$} tensor products. {II}.~{P}uzzles determine facets of the
 {L}ittlewood--{R}ichardson cone, \href{https://doi.org/10.1090/S0894-0347-03-00441-7}{\textit{J.~Amer. Math. Soc.}} \textbf{17}
 (2004), 19--48, \href{https://arxiv.org/abs/math.CO/0107011}{arXiv:math.CO/0107011}.

\bibitem{Macdonald:book}
Macdonald I.G., Symmetric functions and {H}all polynomials, \textit{Oxford Mathematical
 Monographs}, The Clarendon Press, Oxford University Press, New York, 1979.

\bibitem{zonalMathaiEtAl}
Mathai A.M., Provost S.B., Hayakawa T., Bilinear forms and zonal polynomials,
 \textit{Lecture Notes in Statistics}, Vol.~102, \href{https://doi.org/10.1007/978-1-4612-4242-0}{Springer-Verlag}, New York,
 1995.

\bibitem{OO}
Okounkov A., Olshanski G., Shifted {J}ack polynomials, binomial formula, and
 applications, \href{https://doi.org/10.4310/MRL.1997.v4.n1.a7}{\textit{Math. Res. Lett.}} \textbf{4} (1997), 69--78,
 \href{https://arxiv.org/abs/q-alg/9608020}{arXiv:q-alg/9608020}.

\bibitem{Stanley:JackSymFun}
Stanley R.P., Some combinatorial properties of {J}ack symmetric functions,
 \href{https://doi.org/10.1016/0001-8708(89)90015-7}{\textit{Adv. Math.}} \textbf{77} (1989), 76--115.

\bibitem{Z1}
Zuber J.-B., Horn's problem and {H}arish-{C}handra's integrals. {P}robability
 density functions, \href{https://doi.org/10.4171/AIHPD/56}{\textit{Ann. Inst. Henri Poincar\'{e}~D}} \textbf{5}
 (2018), 309--338, \href{https://arxiv.org/abs/1705.01186}{arXiv:1705.01186}.

\end{thebibliography}
\end{document}